%% file: master.tex
\documentclass[11pt]{article}

\usepackage[utf8]{inputenc}
\usepackage{geometry}
\usepackage{amssymb,amsmath,amsthm,mathtools,amsbsy,amsfonts}
\usepackage{subcaption}
\usepackage{graphicx,wasysym}
\usepackage{stmaryrd}
\usepackage{cases}
\usepackage{empheq}
\usepackage{paralist}
\usepackage{hyperref}
\hypersetup{colorlinks=true,linkcolor=red,urlcolor=black,citecolor=green}
\usepackage{multimedia}
\usepackage{tikz}
\usepackage{enumitem}
\usepackage{cite}
\usepackage{multirow}
\usepackage{hhline}
\usepackage{tabularx}
\usepackage{makecell}
\usepackage[ruled,vlined]{algorithm2e}
\SetKwRepeat{Do}{do}{while}
\usepackage{threeparttable}
\usepackage{siunitx}
\usepackage{bbm}
\graphicspath{{./fig/}}

\usepackage{fullpage}
\usepackage[explicit]{titlesec}
\titleformat{\section}{\large\bfseries\center\raggedright}{\thesection}{0.5em}{{#1}}[]
\titleformat{\subsection}[runin]{\bfseries}{\thesubsection}{0.5em}{{#1}}[.]
\titleformat{\subsubsection}[runin]{\bfseries\itshape}{\thesubsubsection}{0.5em}{{#1}}[.]
\titleformat{\paragraph}[runin]{\itshape}{\theparagraph}{0.5em}{{#1}}[.---]
\titlespacing*{\section}{0pt}{0.8\baselineskip}{0.6\baselineskip}
\titlespacing*{\subsection}{0pt}{0.6\baselineskip}{0.4\baselineskip}
\titlespacing*{\subsubsection}{0pt}{0.4\baselineskip}{0.4\baselineskip}
\titlespacing*{\paragraph}{0pt}{0.2\baselineskip}{0.2\baselineskip}

\numberwithin{equation}{section}

\setenumerate{itemsep=-0.2em,topsep=0.3em}
\setitemize{itemsep=-0.2em,topsep=0.3em}

\newcommand\namefont{\normalfont\scshape}
\newcommand\numberfont{\normalfont\scshape}
\newcommand\notefont{\normalfont}
\newtheoremstyle{mystyle} %Name
    {0.3em} %Space above (default, \topsep)
    {0.3em} %Space below (default, \topsep)
    {\itshape} %Body font
    {} %Indent amount 1
    {\normalfont} %Theorem head font
    {.} %Punctuation after theorem head
    {.5em} %Space after theorem head 2
    {{\namefont\thmname{#1}}~{\numberfont\thmnumber{#2}}{\notefont\thmnote{ (#3)}}} %Theorem head spec
\theoremstyle{plain}
\newtheorem{thm}{Theorem}[section]

\newtheorem{rem}[thm]{Remark}

\newtheorem{pb}[thm]{Problem}

\newtheorem{ex}[thm]{Example}

\newtheoremstyle{namedthmstyle} %Name
        {} %Space above (default, \topsep)
	{} %Space below (default, \topsep)
	{\itshape} %Body font
	{} %Indent amount 1
	{\bfseries} %Theorem head font
	{} %Punctuation after theorem head
	{ } %Space after theorem head 2
        {}
\theoremstyle{namedthmstyle}
\newcommand{\thistheoremname}{}
\newtheorem*{genericthm}{\thistheoremname}

\makeatletter
\def\namedlabel#1#2{\begingroup
   \def\@currentlabel{#2}%
   \label{#1}\endgroup
}
\makeatother

\allowdisplaybreaks
\mathtoolsset{showonlyrefs}
\definecolor{darkred}{rgb}{0.6,0.1,0.1}
\definecolor{darkgreen}{rgb}{0.1,0.6,0.1}
\definecolor{darkblue}{rgb}{0.1,0.1,0.6}

\newcommand{\mt}[1]{\mathrm{#1}}
\def\st{\, \left|\right. \,} %"such that" in set notation
\def\:{\colon} %colon for function definition
\newcommand{\abs}[1]{\left\vert#1\right\vert}
\newcommand{\norm}[1]{\left\Vert#1\right\Vert}

\def\bgrad{\boldsymbol{\nabla}}

\DeclareMathOperator{\dive}{div}

\def\R{\mathbb{R}} %Reals
\def\N{\mathbb{N}} %Integers
\def\I{\mathcal{I}} %Index set

\def\e{\varepsilon}

\def\d{\,\mathrm{d}}

\def\p{\partial}

\def\bx{\boldsymbol{x}}

\def\bu{\boldsymbol{u}}
\def\bv{\boldsymbol{v}}

\def\bc{\boldsymbol{c}}

\def\bff{\boldsymbol{f}}
\def\bn{\boldsymbol{n}}

\def\bvarphi{\boldsymbol{\varphi}}

\def\blambda{\boldsymbol{\lambda}}

\def\E{\mathcal{E}} %Energy
\def\K{\mathbb{K}}
\def\k{\mathbbm{k}}

\def\d{\mathbbm{d}}

\newcolumntype{C}[1]{>{\centering\arraybackslash}m{#1}}
\newcolumntype{L}[1]{>{\raggedright\arraybackslash}m{#1}}

\renewcommand{\leq}{\leqslant}
\renewcommand{\geq}{\geqslant}
\newcommand{\sib}[1]{[\si{#1}]}

\begin{document}

\title{Numerical validation of an adaptive model for the determination of nonlinear-flow regions
        in highly heterogeneous porous media}

\author{Alessio Fumagalli$^1$ \and Francesco Saverio Patacchini$^2$}

\date{
        $^1$ Department of Mathematics, Politecnico di Milano, p.za Leonardo da Vinci 32, Milano 20133, Italy\\%
        $^2$ IFP Energies nouvelles, 1 et 4 avenue de Bois-Pr\'eau, 92852 Rueil-Malmaison, France
}

\maketitle

\input{abstract}

%\tableofcontents

\input{todo}
\input{introduction}
\input{model}
\input{validation}
\input{conclusion}

%Journals
%Advances in Water Resources
%Mathematics and Computers in Simulation

\end{document}

%% file: abstract.tex
\begin{abstract}
\noindent An adaptive model for the description of flows in highly heterogeneous porous media is developed in~\cite{FP21,FP23}. There, depending on the magnitude of the fluid's velocity, the constitutive law linking velocity and pressure gradient is selected between two possible options, one better adapted to slow motion and the other to fast motion. We propose here to validate further this adaptive approach by means of more extensive numerical experiments, including a three-dimensional case, as well as to use such approach to determine a partition of the domain into slow- and fast-flow regions. 
%This partitioning, combined with machine-learning and domain-decomposition techniques, will be useful in a future work to provide a fast and accurate flow model in highly heterogeneous porous media. 

\end{abstract}

%%% Local Variables:
%%% mode: latex
%%% TeX-master: "master"
%%% End:

%% file: introduction.tex
\section{Introduction}
\label{sec:introduction}

In the context of flow modeling in porous media, Darcy's law is omnipresent. It is a linear relationship between velocity and pressure gradient, and allows for tensor permeabilities. Its validity is supported both theoretically, via a derivation from Navier--Stokes equations, and empirically, through countless experiments. When velocities deviate from those typically found in porous media, as a result of particularly heterogeneous permeabilities, Darcy's law can, however, be challenged. This was noted experimentally in many studies~\cite{ST14,ST14_2,ZZGYZCC18}. Indeed, in highly permeable areas of a porous medium, such as in fractures, macropores or conduits, the velocity may get high enough to generate nonnegligeable inertial and frictional effects which translate into a nonlinear relation with the pressure gradient. Experimentally, this deviation is noticed as an overestimation of the velocity by Darcy's law compared to that empirically measured. Indeed, from an energetic viewpoint, the effect of the nonlinearity is to penalize high velocities by giving them a high energetic cost. The nonlinearity can also be retrieved theoretically, from averaging and homogenization methods, without, however, providing an explicit formulation of this nonlinearity~\cite{MM00,Neuman77,RM92,Whitaker86,Whitaker96}. Instead, these methods yield the nonlinearity as a perturbation, or correction, of Darcy's law which needs to be fitted empirically. One possible expression for the nonlinearity, which we assume in the present work, is polynomial, meaning that the pressure gradient is a polynomial function of the velocity. This approach includes the widely used Darcy--Forchheimer law, which is a quadratic correction.

Numerically, the nonlinearity thus introduced increases precision but slows down simulations. Nevertheless, the portions of the porous medium where the nonlinear correction has a significant impact are typically small compared to the size of the whole domain. Therefore, the use of the nonlinear perturbation of Darcy's law can be limited to a small subdomain, while the classical Darcy law can be kept in most of the medium. As a result, a-priori knowledge of where these nonlinear regions lie should help improve accuracy without significant computational cost. Indeed, if a domain partition can be devised subdividing the medium into a slow region and a fast region, then a domain-decomposition algorithm~\cite{EHE98,SKN13,AFB19} can be run over this partition where Darcy's law is solved in the slow region and a nonlinear law, such as the Darcy--Forchheimer law, is solved in the fast region. We expect this approach to be more precise than solving Darcy's law globally and faster than solving the nonlinear law globally. The crucial point here is to get the knowledge \emph{a priori} of where each region lies; we aim to answer this point in the present paper, while we leave the domain-decomposition step for future investigation. Methods \emph{a posteriori} are also envisageable~\cite{FRV24,DVY15}, though not explored here.

To locate the linear and nonlinear regions, we propose to use the adaptive model introduced in~\cite{FP21,FP23}. There, the authors derive a model where the constitutive law linking velocity and pressure gradient switches from linear to nonlinear according to a fixed threshold on the seepage-flux magnitude: when the magnitude is below the threshold, the linear law is solved, and when it is above, the nonlinear law is solved instead. As a result, a law which naturally adapts itself given the flux speed is obtained. Although well-posedness of the adaptive model is proved in~\cite{FP21,FP23}, the adaptive law is discontinuous in speed, leading to a jump of pressure gradient at the transition zones between the two laws. Consequently, the adaptive model is rather unsuited for numerical simulation. To circumvent this, the authors in~\cite{FP23} tackle this discontinuity via a regularization argument. More precisely, an energetic, or variational, formulation of the discontinuous problem is derived and smoothed by convolving it with a Gaussian kernel. Then, as an optimization problem, the smooth constitutive model ensues and is shown to converge to the original adaptive model as the regularization parameter vanishes. This regularized law is precisely that which we propose to use in this paper to determine the partition of the medium into slow and fast regions. The result of this task is a labeling of each mesh cell as belonging to one or the other subdomain and strongly depends on the tolerance we give ourselves on the error allowed when solving Darcy's law rather than a nonlinear one; this tolerance is encoded in the very definition of the flux threshold.

This approach may seem to yield a slower model than a globally nonlinear one, since indeed the determination of the regions is carried out by solving another nonlinear, albeit smooth, problem, that is, the regularized problem. Overall then, solving the regularized model to obtain the regions and subsequently implementing a domain decomposition according to the obtained regions should be slower than simply solving the nonlinear model globally. However, we hope that our two-step approach (solving the regularized model followed by a domain decomposition) opens the possibility in a future work to introduce machine-learning techniques to determine the regions in a computationally efficient way. 
%Indeed, one could use the regularized model to generate a sensible dataset through the classification of each cell as belonging to the slow or the fast region, and then train a neural network on these data so that one may then directly infer the partition on which to apply the domain decomposition, thus bypassing the resolution of the regularized problem in practice. 
%\marginpar{non so se lasciare tutta questa parte sul ML, qualcosa si ma qui stiamo dicendo un po' tanto..}

The paper is organized as follows: in Section 2, we summarize the various nonlinear laws considered by our model, we derive an expression for the flux threshold depending on the prescribed tolerance and we recall the adaptive and regularized models of~\cite{FP21,FP23}; in Section 3, we show the validity of the approach by experimenting on four different test cases in two and three space dimensions, the first two cases being inspired by an application in landfill management~\cite{WVCF22} and the remaining two cases being taken from the SPE10 benchmark reservoir scenario~\cite{Christie2001}. The codes used to obtain our results are written in Python and can be found at~\cite{ValAdaptGit}, while the underlying discretization and equation resolution are performed using the open source Python simulation tool PorePy~\cite{KBFSSVB20,PorePyGit}.

%%% Local Variables:
%%% mode: latex
%%% TeX-master: "master"
%%% End:

%% file: model.tex
\section{Model}
\label{sec:model}

We place ourselves in the modeling context of~\cite{FP23}; for completeness, we summarize it here. In addition, we introduce some notions and quantities of particular interest in the present paper.

\subsection{Notation}
\label{sec:notation}

We denote by $\Omega\subset\R^d$ the porous medium, which we assume to be open, bounded and with Lipschitz boundary $\p\Omega$; we write $\bn$ the outward normal unit vector of $\p\Omega$. Also, we write $\R_+ = [0,\infty)$ and $\N_0 = \N \cup \{0\}$.

The unknowns of the problems discussed throughout the paper are the fluid's pressure $p\:\Omega\to\R$ and the seepage flux $\bu\:\Omega\to\R^d$ defined by $\bu = \rho\phi\bv$, where $\phi$ is the medium's porosity and $\rho$ and $\bv$ are the fluid's density and velocity. Furthermore, we write $\nu$ and $\k$ the fluid's kinematic viscosity and the medium's scalar permeability. We suppose that $\rho, \nu, \k\:\Omega\to(0,\infty)$ and $\phi\:\Omega\to(0,1)$ are space-dependent knowns of the problem.

\subsection{Conservation equations}
\label{sec:cons-equat}

We study the stationary flow of the fluid through the porous medium. The conservation of mass then reads
\begin{equation}
  \label{eq:cons-mass}
  \dive\bu = q \quad \text{in $\Omega$},
\end{equation}
where $q\:\Omega\to\R$ is a known fluid mass source, and the conservation of momentum is given by
\begin{equation}
 \label{eq:cons-mom}
 -\bgrad p + \bff = \blambda[\bu] \quad \text{in $\Omega$},
\end{equation}
where $\bff\:\Omega\to\R^d$ is a known vector of external body forces, possibly including gravity. The conservation of momentum shows that the pressure gradient and the external forces balance the \emph{drag force} $-\blambda[\bu]$ undergone by the fluid. We refer to \eqref{eq:cons-mom} as the \emph{seepage}, or \emph{constitutive}, \emph{law} and to $\blambda$ as the \emph{drag operator}.

\subsection{Drag operator}
\label{sec:drag-operator}

We let $\blambda$ be of the form
\begin{equation}
  \label{eq:drag-op}
  \blambda[\bu]\: \Omega \ni \bx \mapsto \alpha(\bx, \norm{\bu(\bx)}) \d(\bx)\bu(\bx),
\end{equation}
where $\d \: \Omega \to (0,\infty)$ is the medium's \emph{drag}, that is, $\d = \nu/ \k$, and $\alpha \: \Omega \times \R_+ \to \R_+$ is a space- and flux-dependent function (the product $\alpha\nu$ possibly being interpreted as the kinematic viscosity perceived by the fluid in motion). This assumption physically means that $\blambda[\bu]$ is colinear with $\bu$ and has magnitude directly dependent on the position in $\Omega$ and on the Euclidean norm $\norm{\bu}$. Moreover, it does not exhibit any anisotropy, be it from the permeability or other quantities. Extensions to tensor permeabilities are envisageable~\cite{FP23,SV01,KL95} and one possibility is in fact mentioned in the context of our three-dimensional test case (cf. Section~\ref{subsec:case3}), but further investigation is required. Regarding anisotropies coming from physical quantities other than permeability, there does not seem to be any consensus on how to incorporate them~\cite{SV01}, so that we do not consider this possibility in this paper.

Volume-averaging and homogenization methods yield relevant information on the expression of the function $\alpha$ directly from the Navier--Stokes equations~\cite{MM00,RM92,Neuman77,Whitaker86,Whitaker96}. By considering the interaction between micro- and macroscopic scales, these methods can be used to motivate the following choice:
\begin{equation*}
  \alpha = 1 + \eta,
\end{equation*}
where $\eta \: \Omega \times \R_+ \to \R_+$ is referred to as the \emph{inertial correction} and behaves linearly about $0$ with respect to its second variable~\cite{Whitaker96}; in particular, $\eta(\bx,0) = 0$ for all $\bx\in\Omega$. With this information, we take $\blambda$ in \eqref{eq:drag-op} as
\begin{equation}
  \label{eq:drag-op-inertial}
  \blambda[\bu]\: \Omega \ni \bx \mapsto \left( 1 + \eta(\bx, \norm{\bu(\bx)}) \right) \d(\bx)\bu(\bx).
\end{equation}
Additional theoretical arguments require that $\eta$ be monotonously nondecreasing and coercive in its second variable~\cite{FP21,FP23}, which intuitively implies that the viscosity perceived by the fluid cannot decrease with its speed.

\subsection{Inertial correction}
\label{sec:inertial-correction}

Unfortunately, volume-averaging and homogenization techniques do not give a unique and explicit expression of the inertial correction, so that assumptions on its form and experimental fitting of coefficients are needed.

\subsubsection{Assumption}
\label{sec:assumption}

We assume that the inertial correction in \eqref{eq:drag-op-inertial} is such that $\eta(\bx,0) = 0$ and $\eta(\bx,\cdot)$ is entire for all $\bx\in\Omega$, i.e., that it can be expanded as
\begin{equation}
  \label{eq:inertial-correction}
  \eta(\bx,a) = \sum_{i=1}^\infty \eta_i(\bx) a^i \quad \text{for all $\bx\in\Omega$ and $a\geq0$},
\end{equation}
where, for all $i\in \N$, the empirical function $\eta_i:\Omega\to\R_+$ is defined by
\begin{equation*}
  \eta_i(\bx) = \frac{1}{i!}\frac{\p^i \eta}{\p a^i}(\bx, 0).
\end{equation*}
A less general expression for the inertial correction is found in~\cite[Equation (5)]{SV01}. 

Note that the entireness of $\eta$ sets the constraint that
\begin{equation*}
  \lim_{i\to\infty} \sqrt[i]{\eta_i(\bx)} = 0 \quad \text{for all $\bx\in\Omega$}.
\end{equation*}
Further note that the nonnegativity of the functions $\eta_i$ is not a mathematically, nor physically, necessary condition but is assumed here for ease of presentation.

\subsubsection{Coefficients}
\label{sec:coefficients}

For all $i \in\N$, we introduce the following dimensionless function $c_i\:\Omega \to \R_+$, referred to as \emph{inertial coefficient}:
\begin{equation*}
  c_i = \left( \frac{\mu}{\sqrt{\k}} \right)^i\eta_i,
\end{equation*}
where $\mu \coloneq \rho\nu$ is the fluid's dynamic viscosity. With this notation, the inertial correction from \eqref{eq:inertial-correction} satisfies
\begin{equation}
  \label{eq:inertial-coeff}
  \eta(\bx,a) = \sum_{i=1}^\infty c_i(\bx) \left( \left(\sqrt{\k}/\mu\right)(\bx)\, a \right)^i \quad \text{for all and $\bx\in\Omega$ and $a\geq0$}.
\end{equation}

The nonzero inertial coefficients in \eqref{eq:inertial-coeff} are precisely the quantities to be chosen experimentally to get an expression for the inertial correction.

\subsubsection{Reynolds number}
In the present porous-medium context~\cite{SV01}, given a flux $\bu$, the \emph{Reynolds number} $\mt{Re}[\bu]\colon \Omega \to \R_+$ is set by
\begin{equation}
  \label{eq:reynolds-num}
  \mt{Re}[\bu] = \frac{\sqrt{\k}}{\mu}\norm{\bu}.
\end{equation}
With this notation, together with \eqref{eq:inertial-coeff}, one can rewrite the seepage law \eqref{eq:cons-mom} as
\begin{equation}
  \label{eq:general-seepage-law-rey}
  -\bgrad p + \bff = \left( 1 + \sum_{i=1}^\infty c_i\, \mt{Re}[\bu]^i \right)\d\,\bu  \quad \text{in $\Omega$}. 
\end{equation}

\subsection{Truncated seepage laws}
\label{sec:trunc-seepage-laws}

In practice, one cannot experimentally tune all the nonzero inertial coefficients since, a priori, these are infinitely many. Therefore, we propose to truncate the law in \eqref{eq:general-seepage-law-rey}, i.e., to keep only a finite number of coefficients. Another way of reducing the number of coefficients is presented in Section \ref{sec:remark-expon-laws}.

We let $\I\subset\N$ be the set of indices we want to keep in the truncated law, so that \eqref{eq:general-seepage-law-rey} rewrites as
\begin{equation}
  \label{eq:general-seepage-law-trunc}
  -\bgrad p + \bff = \left( 1 + \sum_{i\in\I} c_i\, \mt{Re}[\bu]^i \right)\d\,\bu  \quad \text{in $\Omega$}. 
\end{equation}
In the limit case $\I = \N$, we recover the untruncated law \eqref{eq:general-seepage-law-rey}. By considering finite index sets, a wide range of seepage laws can be obtained, many of which are already proposed in various theoretical and numerical studies~\cite{ABHI09,HI11,Kovtunenko23,AFM18,FGL97}. There does not seem to be a general agreement on which is best.

\subsubsection{Darcy's law}
Taking $\mathcal{I} = \emptyset$ in \eqref{eq:general-seepage-law-trunc} gives the widely used \emph{Darcy law}:
\begin{equation*}
  -\bgrad p + \bff = \d\, \bu \quad \text{in $\Omega$}.
\end{equation*}
Note that, here, no coefficient needs fitting so that this law is in fact not empirical.

\subsubsection{Single-power law}
For some $m\in\N$, we take $\mathcal{I} = \{m\}$ in \eqref{eq:general-seepage-law-trunc} to get
\begin{equation*}
  -\bgrad p + \bff = \left( 1 + c_{\mt{F}}\, \mt{Re}[\bu]^m \right) \d\, \bu \quad \text{in $\Omega$},
\end{equation*}
where the coefficient $c_m$ is preferably denoted by $c_{\mt{F}}\:\Omega\to\R$, which is the \emph{Forchheimer coefficient} (or \emph{Ergün number}~\cite{BSD20,BGHA22}), containing information on porosity and tortuosity~\cite{RM92}. When $m=1$, this is called the \emph{Darcy--Forchheimer law}. Note that this formulation stays meaningful for any positive, noninteger value of $m$.

\subsubsection{Power-expansion law}
By setting $\mathcal{I} = \{1, 2, \dots,m\}$ in \eqref{eq:general-seepage-law-trunc} for some $m\in\N$, we yield
\begin{equation*}
  -\bgrad p + \bff = \left( 1 + c_1\, \mt{Re}[\bu] + c_2\, \mt{Re}[\bu]^2 + \cdots + c_m\, \mt{Re}[\bu]^m \right) \d\, \bu \quad \text{in $\Omega$},
\end{equation*}
which, as with the single-power law, includes the Darcy--Forchheimer law. Single-power and power-expansion laws may thus be called \emph{generalized Forchheimer laws}.

\subsection{Comparing seepage laws}
\label{sec:comp-seep-laws}

Given a flux $\bu$ and two drag operators $\blambda$ and $\hat\blambda$, the resulting forces $-\bgrad p + \bff$ and $-\bgrad \hat p + \bff$ are obtained by the laws (cf. \eqref{eq:cons-mom})
\begin{equation}
  \label{eq:two-laws}
  \begin{cases}
    -\bgrad p + \bff = \blambda[\bu],\\
    -\bgrad \hat p + \bff = \hat\blambda[\bu],
  \end{cases}
  \quad \text{in $\Omega$}.
\end{equation}
The drag operator $\hat\blambda$ is referred to as the \emph{reference} operator.

\subsubsection{Untruncated error}
\label{sec:untrunc-error}

Subtracting the two equations in \eqref{eq:two-laws} and using the untruncated law \eqref{eq:general-seepage-law-rey} with inertial coefficients $\bc \coloneq (c_i)_{i\in\N}$ and $\hat\bc \coloneq (\hat c_i)_{i\in\N}$, the \emph{local error} $d_{\bc, \hat\bc}[\bu] \colon \Omega \to \R_+$ of the first law relative to the second is defined by
\begin{equation*}
  d_{\bc, \hat\bc}[\bu] = \frac{\abs{ \sum_{i=1}^\infty \left(c_i - \hat c_i\right) \mt{Re}[\bu]^i}}{1 + \sum_{i=1}^\infty \hat c_i\, \mt{Re}[\bu]^i}.
\end{equation*}

\subsubsection{Truncated error}
\label{sec:truncated-error}

For truncated laws, subtracting again the two equations in \eqref{eq:two-laws} and using this time \eqref{eq:general-seepage-law-trunc} with coefficients $\bc \coloneq (c_i)_{i\in\N}$ and $\hat\bc \coloneq (\hat c_i)_{i\in\N}$ and index sets $\I$ and $\hat\I$, we define the local error $d_{\bc, \I; \hat\bc, \hat\I}[\bu] \colon \Omega \to \R_+$ as
\begin{equation*}
  d_{\bc, \I; \hat\bc, \hat\I}[\bu] = \frac{\abs{ \sum_{i\in\I} c_i\, \mt{Re}[\bu]^i - \sum_{i\in\hat\I}\hat c_i\, \mt{Re}[\bu]^i}}{1 + \sum_{i\in\hat\I} \hat c_i\, \mt{Re}[\bu]^i}.
\end{equation*}

\begin{ex}[equal coefficients and Darcy--Forchheimer]
  \label{ex:DF}
  Of particular interest to us is the case $\bc = \hat \bc$ and $\I = \emptyset$. In this setting, the above error becomes
\begin{equation*}
  d_{\hat\bc, \hat\I}[\bu] \coloneq d_{\hat\bc, \emptyset; \hat\bc, \hat\I}[\bu] = \frac{\sum_{i\in\hat\I} \hat c_i\, \mt{Re}[\bu]^i}{1 + \sum_{i\in\hat\I} \hat c_i\, \mt{Re}[\bu]^i} < 1.
\end{equation*}
More specifically, we shall consider $\hat\I = \{m\}$, with $m\in\N$ and $\hat c_m = c_{\mt{F}}$, that is, the case when the first law is Darcy's and the reference law is a single-power law; then,
\begin{equation}
  \label{eq:error-DF}
  d_{c_{\mt{F}}, m}[\bu] \coloneq d_{\hat\bc, \emptyset; \hat\bc, \{m\}}[\bu] = \frac{c_{\mt{F}}\, \mt{Re}[\bu]^{m}}{1 + c_{\mt{F}}\, \mt{Re}[\bu]^{m}} \eqcolon \frac{\mt{Fo}_{m}[\bu]^{m}}{1 + \mt{Fo}_{m}[\bu]^{m}},
\end{equation}
where the \emph{$m$th Forchheimer number} $\mt{Fo}_{m}[\bu]\colon \Omega \to \R_+$ is defined by
\begin{equation*}
  \mt{Fo}_{m}[\bu] = c_{\mt{F}}^{1/m} \frac{\sqrt{\k}}{\mu}\norm{\bu} = c_{\mt{F}}^{1/m} \mt{Re}[\bu].
\end{equation*}
The dimensionless number $\mt{Fo}_{m}[\bu]^{m}$ quantifies the weight of the nonlinear term in the single-power law relative to $1$, i.e., relative to the linear term. The first Forchheimer number, resulting from considering the Darcy--Forchheimer law as reference, is simply referred to as the Forchheimer number in the literature~\cite{MMV11,ZG06}.
\end{ex}

\subsubsection{Flux threshold and subdomains}
\label{sec:flux-thresh-subd}

In the setting given by Example \ref{ex:DF}, there exists a function $\Delta_{\hat\bc,\hat\I} \colon \Omega\times \R_+ \to [0,1)$ such that, for all $\bx\in\Omega$, there holds
\begin{equation*}
  d_{\hat\bc, \hat\I}[\bu](\bx) = \Delta_{\hat\bc, \hat\I}(\bx, \norm{\bu(\bx)}),
\end{equation*}
and $\Delta_{\hat\bc, \hat\I}(\bx, \cdot)$ is an increasing bijection. Let $\delta\in[0,1)$, referred to as \emph{error tolerance}.

\paragraph{Threshold}

By bijectivity, there is a unique $w_\delta\colon \Omega \to \R_+$ so that $\Delta_{\hat\bc, \hat\I}(\bx, w_{\delta}(\bx)) = \delta$ for all $\bx\in\Omega$. We then set the \emph{flux threshold} $\bar u_\delta$ as
\begin{equation}
  \label{eq:thresh}
  \bar u_{\delta} = \inf_{\Omega} w_\delta \in \R_+.
\end{equation}

When the reference law is a single-power law with exponent $\hat m\in\N$, then, from \eqref{eq:error-DF}, one finds
\begin{equation*}
  w_\delta = \mt{Fo}_\delta\, \frac{\mu}{c_{\mt{F}}^{1/m}\sqrt{\k}},
\end{equation*}
where $\mt{Fo}_\delta \coloneq \sqrt[m]{\delta/(1 - \delta)}$ is a \emph{critical} Forchheimer number, which yields
\begin{equation*}
  \bar u_\delta = \mt{Fo}_\delta\, \inf_\Omega \frac{\mu }{c_{\mt{F}}^{1/m}\sqrt{\k}}.
\end{equation*}
In particular, when $m = 1$, i.e., the reference law is Darcy--Forchheimer, we get
\begin{equation*}
  w_\delta = \mt{Fo}_\delta\, \frac{\mu}{c_{\mt{F}} \sqrt{\k}} \quad \text{and} \quad \bar u_\delta = \mt{Fo}_\delta\, \inf_\Omega \frac{\mu}{c_{\mt{F}}\sqrt{\k}} = \frac{\delta}{1-\delta} \inf_\Omega \frac{\mu}{c_{\mt{F}}\sqrt{\k}}.
\end{equation*}

\paragraph{Subdomains}

We then define the subdomains
\begin{equation*}
  \Omega_{\delta}[\bu] = \{ \bx \in \Omega \st \norm{\bu(\bx)} < \bar u_{\delta} \} \quad \text{and} \quad \hat\Omega_{\delta}[\bu] = \{ \bx \in \Omega \st \norm{\bu(\bx)} > \bar u_{\delta} \},
\end{equation*}
which we respectively refer to as the \emph{slow} and \emph{fast} subdomains.

If $\bx\in\Omega_{\delta}[\bu]$, we know that solving at $\bx$ the law given by $\blambda$ (cf. first equation in \eqref{eq:two-laws}) does not generate a local error greater than $\delta$ relative to solving at $\bx$ the law stemming from $\hat\blambda$ (cf. second equation in \eqref{eq:two-laws}). Indeed, suppose $\bx\in \Omega_\delta[\bu]$; then,
\begin{equation*}
  \norm{\bu(\bx)} < \bar u_\delta \leq w_\delta(\bx), 
\end{equation*}
so that 
\begin{equation*}
  d_{\hat\bc, \hat\I}[\bu](\bx) = \Delta_{\hat\bc, \hat\I}(\bx, \norm{\bu(\bx)}) < \Delta_{\hat\bc, \hat\I}(\bx, w_\delta(\bx)) = \delta.
\end{equation*}
Hence, by solving the seepage law for $\blambda$ in $\Omega_\delta[\bu]$ and that for $\hat\blambda$ in $\hat\Omega_\delta[\bu]$, we make sure that the local error stays below the error tolerance $\delta$. Also note that if $\bx\in\Omega$ and $\mathrm{Fo}_m[\bu](\bx) > \mathrm{Fo}_\delta$, then $\bx\in\hat\Omega_\delta[\bu]$; this is why $\mathrm{Fo}_\delta$ is seen as a critical Forchheimer number above which we are certain that nonlinear effects are nonnegligeable.

\subsection{Final model}
\label{sec:final-model}

Suppose we are in the framework of Section \ref{sec:comp-seep-laws} with $\delta \in [0, 1)$. 

Furthermore, let $\Sigma_{\mt{v}},\Sigma_{\mt{p}} \subset \p\Omega$ be relatively open in $\p\Omega$ (i.e., $\Sigma_{\mt{v}}$ and $\Sigma_{\mt{p}}$ are each the intersection of an open subset of $\R^d$ with $\p\Omega$) and such that $\p\Omega = \overline{\Sigma_{\mt{v}} \cup \Sigma_{\mt{p}}}$ and $\Sigma_{\mt{v}} \cap \Sigma_{\mt{p}} = \emptyset$. We also let $u_0\: \Sigma_{\mt{v}} \rightarrow \R$ and $p_0\:\Sigma_{\mt{p}}\to\R$ be functions setting the conditions on the boundary for the flux and pressure. For simplicity, a map on $\Omega$ and its trace on $\p\Omega$ are denoted by the same symbol.

\subsubsection{Adaptive setting}
\label{sec:discont-setting}

We consider the problem where the seepage law for $\blambda$ applies in the slow subdomain and that for $\hat\blambda$ in the fast one. In practice, the reference drag operator $\hat\blambda$ is numerically more costly to handle, so that by only solving it in the fast subdomain, rather than the entire domain, we reduce the overall computational cost. This discontinuous problem, where the discontinuity in flux arises from the jump occurring between the seepage laws, is referred to as the \emph{adaptive problem} and is as follows:
\begin{pb}[adaptive]
  \label{pb:discont}
  Find $\bu\:\Omega\to\R^d$ and $p\:\Omega\to\R$ such that
  \begin{equation}
    \label{eq:discont}
    \begin{cases}
      \dive\bu = q & \text{in $\Omega$},\\
      -\bgrad p + \bff = \blambda[\bu] & \text{in $\Omega_\delta[\bu]$},\\
      -\bgrad p + \bff = \hat\blambda[\bu] & \text{in $\hat\Omega_\delta[\bu]$},\\
      \bu \cdot \bn = u_0 & \text{on $\Sigma_{\mt{v}}$},\\
      p = p_0 & \text{on $\Sigma_{\mt{p}}$}.
    \end{cases}
  \end{equation}
\end{pb}

The rigorous setting in which to solve such problem is described in~\cite{FP23}, from where well-posedness can be deduced via a straightforward generalization of~\cite[Assumption~4.1 and Corollary~4.4]{FP23} to space-dependent
physical parameters, such as $\mu$, $\nu$, $\k$ and the inertial coefficients, which we further assume to be bounded below and above by strictly positive constants in $\Omega$. In~\cite{FP23}, the appropriate functional framework is developed and a multivalued formulation is derived to tackle the flux discontinuity at the transition zone between the laws.

\begin{rem}[average pressure]
  \label{rem:press-ave}
  If the boundary piece $\Sigma_{\mt{p}}$ verifies $\mt{Vol}^{d-1}(\Sigma_{\mt{p}}) = 0$, where $\mt{Vol}^{d-1}$ is the $(d-1)$-dimensional Lebesgue measure, to ensure uniqueness of the pressure satisfying Problem~\ref{pb:discont}, one imposes a constraint on the average of $p$:
\begin{equation}
  \label{eq:pressure-average}
  \frac{1}{\abs{\Omega}}\int_\Omega p = \bar p,
\end{equation}
for a given $\bar p\in\R$. Tacitly, we thus require \eqref{eq:pressure-average} in \eqref{eq:discont} whenever $\mt{Vol}^{d-1}(\Sigma_{\mt{p}}) = 0$.
\end{rem}

\subsubsection{Regularized setting}
\label{sec:regul-setting}

Numerically, it is not diserable to solve a discontinuous problem such as Problem \ref{pb:discont}. Indeed, since we use an explicit fixed-point algorithm, convergence is not guaranteed without continuity. We therefore follow the approach developed once again in~\cite{FP23}, which consists in averaging the two drag operators $\blambda$ and $\hat \blambda$ via a convolution to obtain a smooth law. This introduces an additional parameter, $\e > 0$, which vanishes as the regularization weakens. The details of this regularization are found in~\cite{FP23} and the corresponding Python code at~\cite{ValAdaptGit}. Still, for completeness, we give here an informal derivation of this regularized law and the respective regularized problem. 

The authors in~\cite{FP23} exploit the variational, or energetic, formulation of the adaptive model to regularize the law via a convolution of the underlying energy functional. It is shown there that a pair $(\bu, p)$ is solution to Problem \ref{pb:discont} if and only if it is a saddle point of the energy functional $\E$ defined by
\begin{equation*}
    \E(\bvarphi, \psi) = \int_\Omega \bff \cdot \bvarphi- \int_\Omega q \psi + \int_{\Sigma_{\mt{v}}} u_0 \psi - \int_\Omega \bgrad \psi \cdot \bvarphi - \mathcal{D}(\bvarphi),
\end{equation*}
where $\mathcal{D}$ is called the \emph{dissipation} and is minimized on a certain restricted function space by a certain component of the solution $\bu$ (cf.~\cite{FP21,FP23}). It encodes the information of the two chosen drag operators $\blambda$ and $\hat \blambda$ and is of the form
\begin{equation*}
  \mathcal{D}(\bvarphi) = \frac12 \int_\Omega \d \left( \Psi(\norm{\bvarphi}^2) \mathbbm{1}_{\Omega_\delta[\bvarphi]} + \hat\Psi(\norm{\bvarphi}^2) \mathbbm{1}_{\hat \Omega_\delta[\bvarphi]} \right),
\end{equation*}
where $\mathbbm{1}_A$ stands for the indicator function of set $A$, and $\Psi, \hat \Psi \: \R_+ \to \R_+$ satisfy
\begin{equation*}
  \blambda[\bvarphi] = \Psi'(\norm{\bvarphi}^2)\d\, \bvarphi \quad \text{and} \quad \hat \blambda[\bvarphi] = \hat \Psi'(\norm{\bvarphi}^2)\d\, \bvarphi.
\end{equation*}
We note that the discontinuity of the adaptive law is reflected in the discontinuity of the integrand of $\mathcal{D}$ with respect to the variable $\norm{\bvarphi}$. In fact, $\mathcal{D}$ is not differentiable and merely has a subdifferential, which is equal to the right-hand side of the adaptive law in Problem \ref{pb:discont}. To make the dissipation smooth, the integrand of $\mathcal{D}$ is convolved with a Gaussian kernel $G_\e$ of mean $0$ and standard deviation some fixed $\e>0$. The convolution is applied on the space of flux magnitudes and yields the regularized dissipation $\mathcal{D}_\e$, defined as
\begin{equation*}
  \mathcal{D}_\e(\bvarphi) = \frac12 \int_\Omega  \d\, G_\e * \left[ \left( \Psi(\norm{\bvarphi}^2) \mathbbm{1}_{\Omega_\delta[\bvarphi]} + \hat\Psi(\norm{\bvarphi}^2) \mathbbm{1}_{\hat \Omega_\delta[\bvarphi]} \right) \right],
\end{equation*}
and the corresponding energy $\mathcal{E}_\e$, given by
\begin{equation*}
    \E_\e(\bvarphi, \psi) = \int_\Omega \bff \cdot \bvarphi - \int_\Omega q \psi + \int_{\Sigma_{\mt{v}}} u_0 \psi - \int_\Omega \bgrad \psi \cdot \bvarphi - \mathcal{D}_\e(\bvarphi).
\end{equation*}
The regularized problem then consists in finding a pair $(\bu_\e, p_\e)$ which is a saddle point of $\E_\e$. Equivalently, after differentiation of the regularized dissipation, it consists in letting the seepage law be
\begin{equation}
    \label{eq:regularized-law}
    -\bgrad p_\varepsilon + \bff = \blambda_{\delta, \e}[\bu_\e],
\end{equation}
where the regularized drag operator $\blambda_{\delta, \e}$ is given by
\begin{equation*}
  \blambda_{\delta, \e}[\bu_\e] =  \d\, G_\e' * \left[ \left( \Psi(\norm{\bu_\e}^2) \mathbbm{1}_{\Omega_\delta[\bu_\e]} + \hat\Psi(\norm{\bu_\e}^2) \mathbbm{1}_{\hat \Omega_\delta[\bu_\e]} \right) \right].
\end{equation*}
This drag operator contains the information of the threshold flux $\bar u_\delta$ via the indicator functions and the resulting \emph{regularized problem} is the following:
\begin{pb}[regularized]
  \label{pb:regul}
  Find $\bu_\e\:\Omega\to\R^d$ and $p_\e\:\Omega\to\R$ such that
  \begin{equation*}
    \begin{cases}
      \dive\bu_\e = q & \text{in $\Omega$},\\
      -\bgrad p_\e + \bff = \blambda_{\delta, \e}[\bu_\e] & \text{in $\Omega$},\\
      \bu_\e \cdot \bn = u_0 & \text{on $\Sigma_{\mt{v}}$},\\
      p_\e = p_0 & \text{on $\Sigma_{\mt{p}}$}.
    \end{cases}
  \end{equation*}
\end{pb}

Problem \ref{pb:regul} is indeed the problem we validate numerically in Section \ref{sec:numerical-validation} as a tool to generate partitions of porous media into fast and slow subdomains. Thanks to~\cite[Corollary~5.3 and Theorem~5.4]{FP23}, it is well posed and converges in a weak sense to the adaptive problem as $\e\to0$. Thus, the slow and fast subdomains $\Omega_\delta[\bu_\e]$ and $\hat \Omega_\delta[\bu_\e]$ it generates differ from those stemming from the adaptive model, namely, $\Omega_\delta[\bu]$ and $\hat \Omega_\delta[\bu]$, but are expected, though not yet proven, to be not too different for a reasonably small choice of $\e$, such as $\varepsilon = 0.1$, which is the value we take for the numerical experiments below.

\subsection{Remark: exponential seepage laws}
\label{sec:remark-expon-laws}

As an intriguing observation, we show here an alternative way of reducing the number of inertial coefficients in the untruncated law \eqref{eq:general-seepage-law-rey}, which yields laws different from those derived in Section \ref{sec:trunc-seepage-laws} and, to our knowledge, not found in the current literature, although we do not explore them further in the present paper.

\subsubsection{Derivation}
The idea is to reduce the number of coefficients to a single one by requiring that there be a coefficient $c \: \Omega \to \R$ satisfying
\begin{equation}
  \label{eq:single-coeff}
  i!\, c_i = c^i, \qquad \text{$i\in\N$}.
\end{equation}
Moreover, we give ourselves a truncating parameter $m\in\N_0\cup\{\infty\}$.

From \eqref{eq:inertial-coeff}, for all $\bx\in\Omega$ and $a\geq0$, we yield
\begin{align*}
  1 + \eta(\bx,a) &= \sum_{i=0}^m \frac{1}{i!}\left( \left(c\sqrt{\k}/\mu\right)(\bx)\, a \right)^i\\
  &= G_m\left( \left(c\sqrt{\k}/\mu\right)(\bx)\, a \right) \exp\left( \left(c\sqrt{\k}/\mu\right)(\bx)\, a \right),
\end{align*}
where $G_m\colon \R_+ \to (0, 1]$ is defined, for all $b\in\R_+$, as
\begin{equation}
  \label{eq:gm}
  G_m(b) = 
  \begin{cases}
    \dfrac{\Gamma(m+1,b)}{\Gamma(m+1,0)} & \text{if $m\in\N_0$},\\
    1 & \text{if $m=\infty$},
  \end{cases}
\end{equation}
where $\Gamma$ is the upper incomplete gamma function~\cite{CZ01}. Thus, \eqref{eq:general-seepage-law-rey} gives
\begin{equation}
  \label{eq:seepage-law-exp}
  -\bgrad p + \bff = G_m\left( c\, \mt{Re}[\bu] \right) \exp\left( c\, \mt{Re}[\bu] \right) \d\, \bu \quad \text{in $\Omega$}.
\end{equation}
We may call this law \emph{full exponential law} if $m=\infty$ and \emph{truncated exponential law} if $m\in\N_0$. This is again a generalized Forchheimer law since, from it, we retrieve the Darcy--Forchheimer law if $m=1$ and $c = c_{\mathrm{F}}$, and indeed Darcy's law if $m=0$. However, when $m>1$, we yield laws different from those presented in Section \ref{sec:trunc-seepage-laws}.

Mathematically, this formulation makes sense for some nonpositive and noninteger values of $m$ since the gamma function is well defined on some of those values. Further note that $G_m\uparrow G_\infty$ pointwise as $m\to\infty$, so that indeed the truncated exponential law approaches the full law as the truncating parameter $m$ grows. Note, however, that well-posedness for such exponential laws is not guaranteed by the arguments in~\cite{FP23} since operator boundedness does not hold in this case.

\subsubsection{Local error}

Following the reasoning of Section \ref{sec:truncated-error} for exponential laws, using \eqref{eq:seepage-law-exp} with coefficients $c$ and $\hat c$ and truncating parameters $m$ and $\hat m$, the local error $d_{c, m; \hat c, \hat m}[\bu] \colon \Omega \to \R_+$ may be defined as 
\begin{equation*}
  d_{c, m; \hat c, \hat m}[\bu] = \frac{\abs{ G_m(c\, \mt{Re}[\bu])\exp(c\, \mt{Re}[\bu]) - G_{\hat m}(\hat c\,\mt{Re}[\bu])\exp(\hat c\,\mt{Re}[\bu])}}{G_{\hat m}(\hat c\,\mt{Re}[\bu])\exp(\hat c\,\mt{Re}[\bu])}.
\end{equation*}

\begin{ex}[equal coefficients and full exponential]
  \label{ex:full-exp}
  The particular case $c = \hat c$ and $\hat m > m$ simplifies into 
  \begin{equation*}
    d_{m; c, \hat m}[\bu] \coloneq d_{c, m; c, \hat m}[\bu] = 1 - \frac{G_m(c\,\mt{Re}[\bu])}{G_{\hat m}(c\,\mt{Re}[\bu])} < 1,
  \end{equation*}
  which, when $\hat m = \infty$, further cuts down to
  \begin{equation}
    \label{eq:error-full-exp}
    d_{m; c}[\bu] \coloneq d_{c, m; c, \infty}[\bu] = 1 - G_m(c\, \mt{Re}[\bu]).
  \end{equation}
\end{ex}

The considerations in Section \ref{sec:flux-thresh-subd} still hold in the exponential setting given in Example \ref{ex:full-exp}. Indeed, one can find $\Delta_{m, ;c, \hat m} \colon \Omega\times \R_+\to [0, 1)$ so that, for all $\bx\in\Omega$, there holds
\begin{equation*}
  d_{m; c, \hat m}[\bu](\bx) = \Delta_{m; c, \hat m}(\bx, \norm{\bu(\bx)}),
\end{equation*} 
where $\Delta_{m; c, \hat m}(\bx, \cdot)$ is an increasing bijection. Using the same definition as in \eqref{eq:thresh}, when the reference law is the full exponential, from \eqref{eq:error-full-exp}, one gets
\begin{equation*}
  w_\delta = G_m^{-1}(1 - \delta) \frac{\mu}{c \sqrt{\k}},
\end{equation*}
which yields
\begin{equation*}
  \bar u_\delta = G_m^{-1}(1 - \delta) \inf_\Omega \frac{\mu}{c \sqrt{\k}}.
\end{equation*}
The inverse $G_m^{-1}$ can, in general, be computed numerically; but, if, for instance $m=0$, that is, the drag operator $\blambda$ is Darcy's, then $G_0^{-1} = -\ln$, so that
\begin{equation*}
  \bar u_\delta = -\ln(1 - \delta) \inf_\Omega \frac{\mu}{c\sqrt{\k}}.
\end{equation*}

%%% Local Variables:
%%% mode: latex
%%% TeX-master: "master"
%%% End:

%% file: validation.tex
\section{Numerical validation}
\label{sec:numerical-validation}

In this section, we present four test cases to validate the performance and quality of our approach, that is, the results of the numerical resolution of Problem \ref{pb:regul} using Raviart--Thomas finite elements and an explicit fixed-point algorithm, as described in detail in~\cite{FP23}. Note that, in this section, as opposed to the rest of the paper, we refer to the model in Problem \ref{pb:regul} as the \emph{adaptive} model, rather than \emph{regularized}; this is to keep stressing the fact that there always is an underlying adaptability between the nonlinear and linear laws in the numerically solved model.

In Section \ref{subsec:case1}, two examples are proposed with two or more highly conductive channels and a background heterogeneity, as inspired by~\cite{WVCF22}. The examples in Sections \ref{subsec:case2} and \ref{subsec:case3} are based on the SPE10 benchmark study~\cite{Christie2001}, where we consider, respectively, a single layer and multiple layers.

\subsection{Highly permeable channels} \label{subsec:case1}
% (case2 git)

This test case is based on the examples reported in \cite{WVCF22} for a Darcy--Forchheimer system, and adapted here for our purposes. We consider in particular two configurations depending on the number of highly conductive channels that are present in the problem. In the first case, reported in Section \ref{subsubsub:two_channels}, two vertical channels are present, while in Section \ref{subsubsub:network} a network of vertical and horizontal channels is considered.

\begin{figure}[!ht]
    \centering
    \includegraphics[scale=0.1]{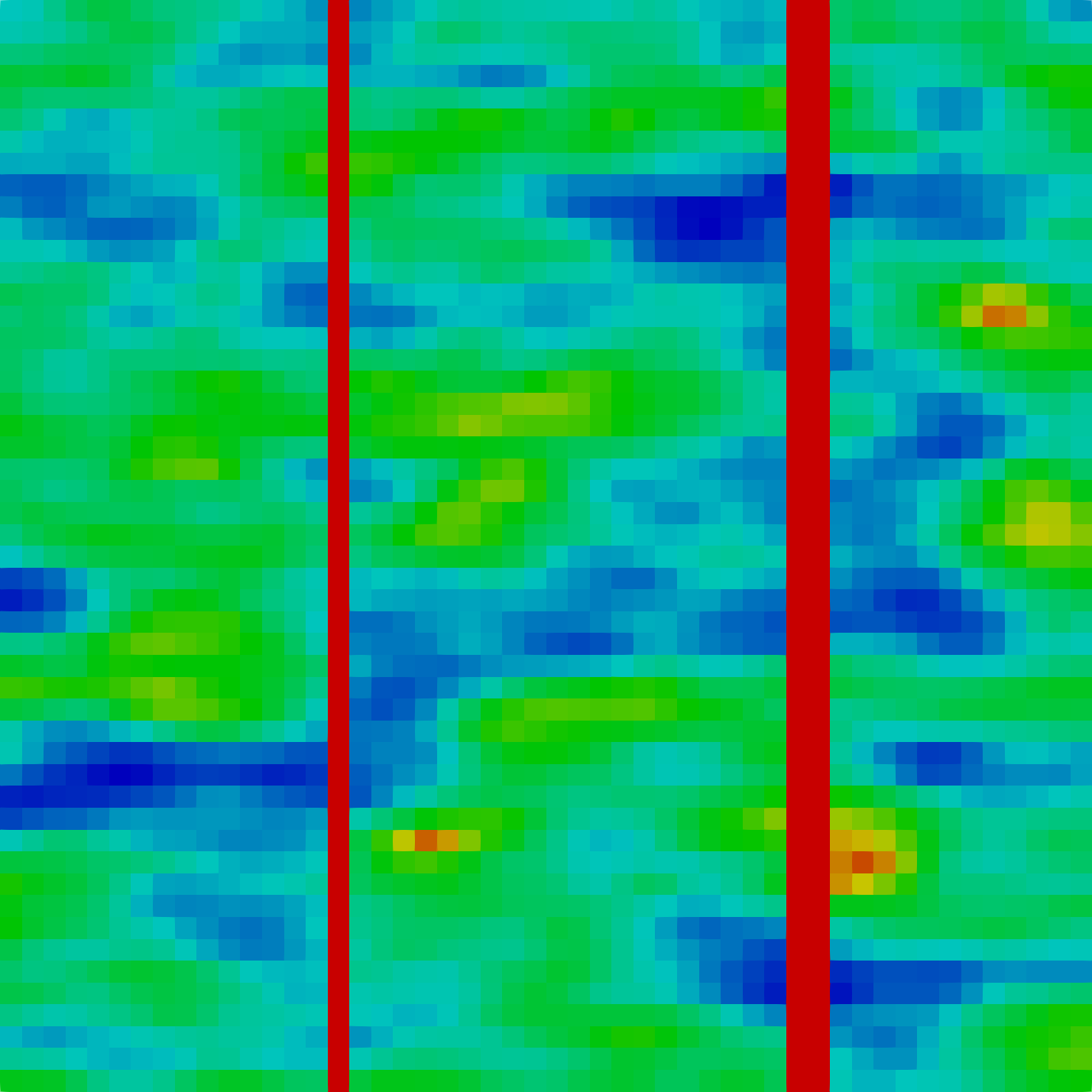}%
    \hspace*{0.025\textwidth}%
    \includegraphics[scale=0.15]{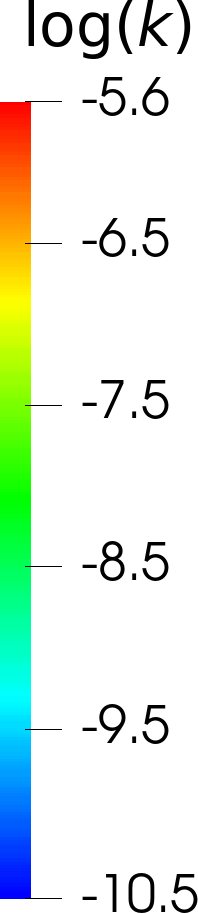}%
    \hspace*{0.06\textwidth}%
    \includegraphics[scale=0.1]{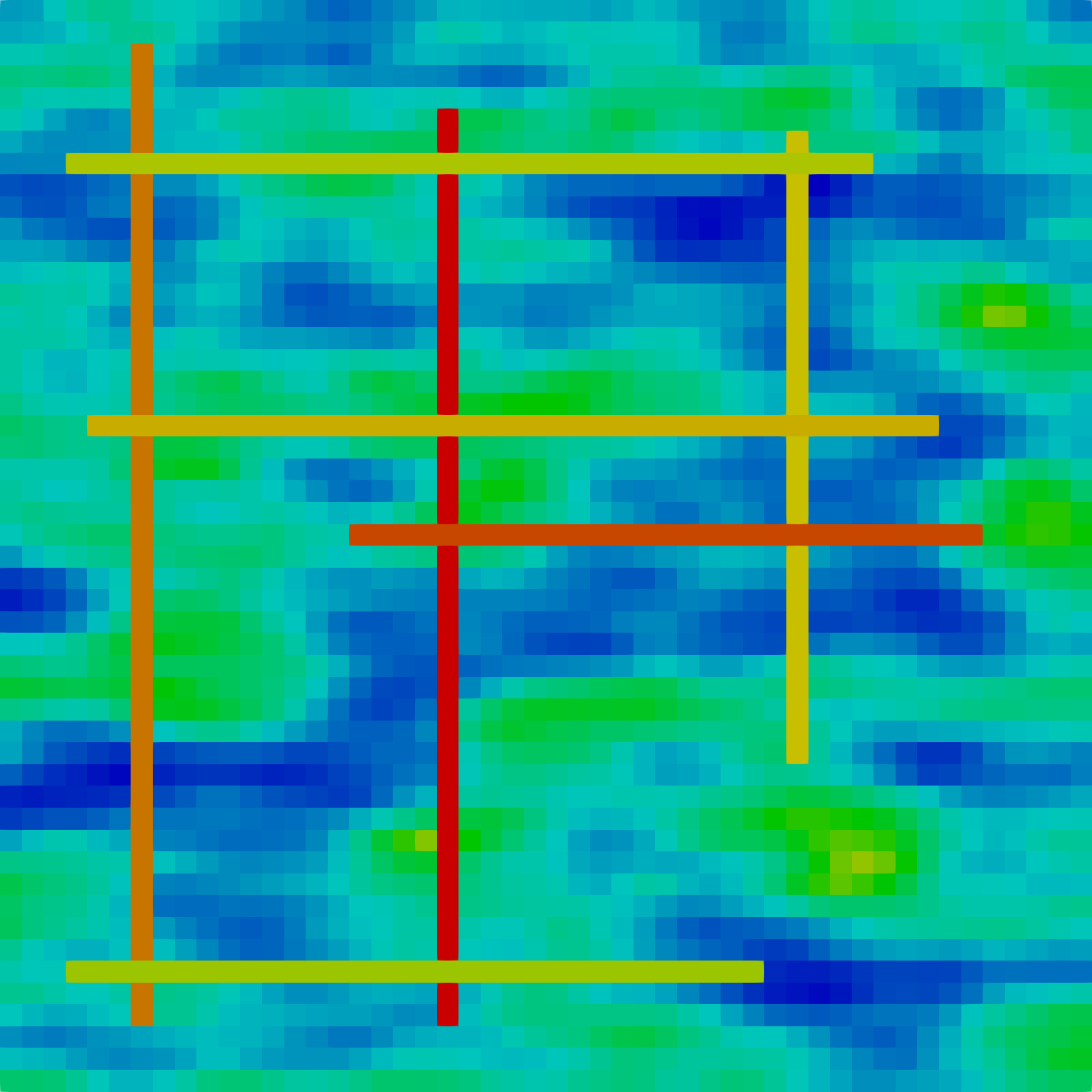}%
    \hspace*{0.025\textwidth}%
    \includegraphics[scale=0.15]{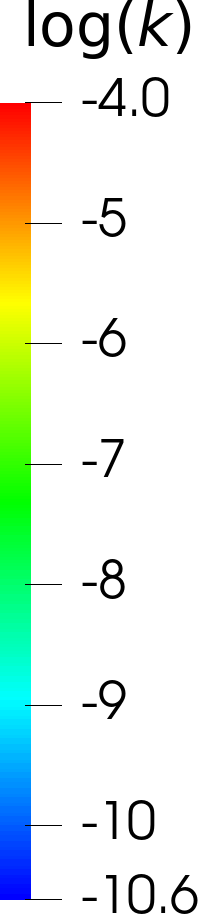}%
    \caption{On the left, the log-permeability for the example of Section \ref{subsubsub:two_channels}; on the right, for Section \ref{subsubsub:network}.}
    \label{fig:case1}
\end{figure}

For both tests, we consider the following physical data for the fluid: dynamic viscosity $\mu = 0.001 \sib{\pascal\second}$ and density $\rho = 998.0 \sib{\kilo\gram\per\cubic\meter}$. For the porous medium, 
%we take the porosity $\phi = 94.0 \sib{\%}$ everywhere and 
the permeability is reported in Figure \ref{fig:case1} for each configuration. To pass from permeability to porosity and vice-versa when needed, we use the Kozeny--Carman law:
\begin{equation*}
  \k = \k_{\mt{ref}}\frac{(1-\phi_{\mt{ref}})^2}{\phi_{\mt{ref}}^3} \frac{\phi^3}{(1-\phi)^2},
\end{equation*}
where $\k_{\mt{ref}} = 1.0152441851\cdot 10^{-9} \sib{\square\meter}$ and $\phi_{\mt{ref}} = 35 \sib{\percent}$ (cf.~\cite{WVCF22}).
%$\bar{u} = 0.13297 \sib{\kilo\gram\per\square\meter\per\second}$

Unless stated otherwise, we fix the error tolerance $\delta = 0.1 \sib{-}$, yielding the critical Forchheimer number $\mt{Fo}_\delta \simeq 0.11 \sib{-}$ (cf. Section \ref{sec:comp-seep-laws}); moreover, as chosen in~\cite{WVCF22}, we let $m=1$, corresponding to the classical Darcy--Forchheimer model, with coefficient $c_{\mt{F}} = 0.55\sib{-}$ (cf. Section \ref{sec:trunc-seepage-laws}). 

For each case, we consider two scenarios for different values of the top inflow magnitude: \textit{Scenario a}, with $u_{\mt{in}} = 7 \cdot 10^{-3} \sib{\kilo\gram\per\square\meter\per\second}$; \textit{Scenario b}, with $u_{\mt{in}} = 7\cdot 10^{-2} \sib{\kilo\gram\per\square\meter\per\second}$. On the left and right parts of the boundary, we set zero flux, while on the bottom part, we impose the hydrostatic pressure $p_{\mt{atm}} + \rho g h$, where $p_{\mt{atm}} = 1.01325 \cdot 10^5 \sib{\pascal}$, $h=10 \sib{\meter}$ and $g=9.81\sib{\meter\per\square\second}$. 

\subsubsection{Two vertical channels} \label{subsubsub:two_channels}

We consider here the case where only two highly conductive channels are present in the domain, spanning the domain from top to bottom. Their permeability is roughly two orders of magnitude higher than that of the surrounding porous medium, and they are aligned with the flow direction given by the top boundary condition (cf. left in Figure \ref{fig:case1}). Depending on the flux intensity, \textit{Scenario a} or \textit{Scenario b}, we expect to observe a Darcy-Forchheimer, i.e., nonlinear, flow develop inside these channels and a Darcy, i.e., linear, flow outside.

\begin{figure}[!ht]
    \centering
    \includegraphics[scale=0.075]{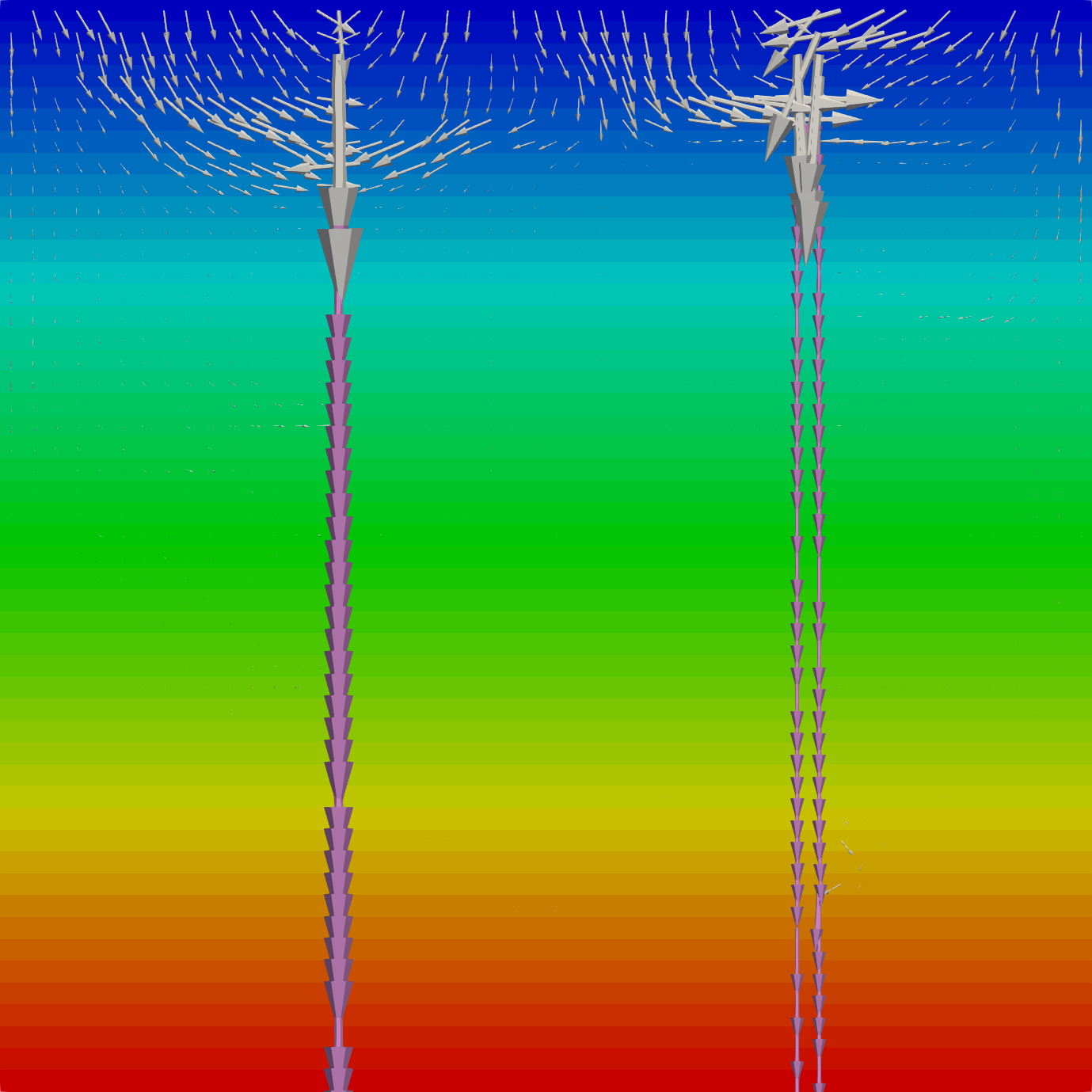}%
    \hspace*{0.0125\textwidth}%
    \includegraphics[scale=0.1]{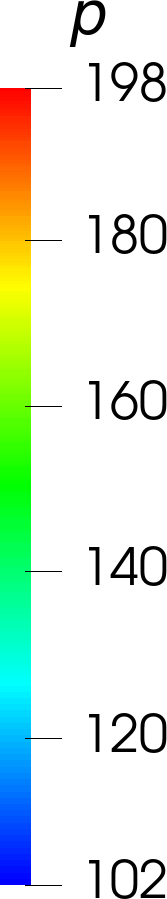}%
    \hspace*{0.02\textwidth}%
    \includegraphics[scale=0.075]{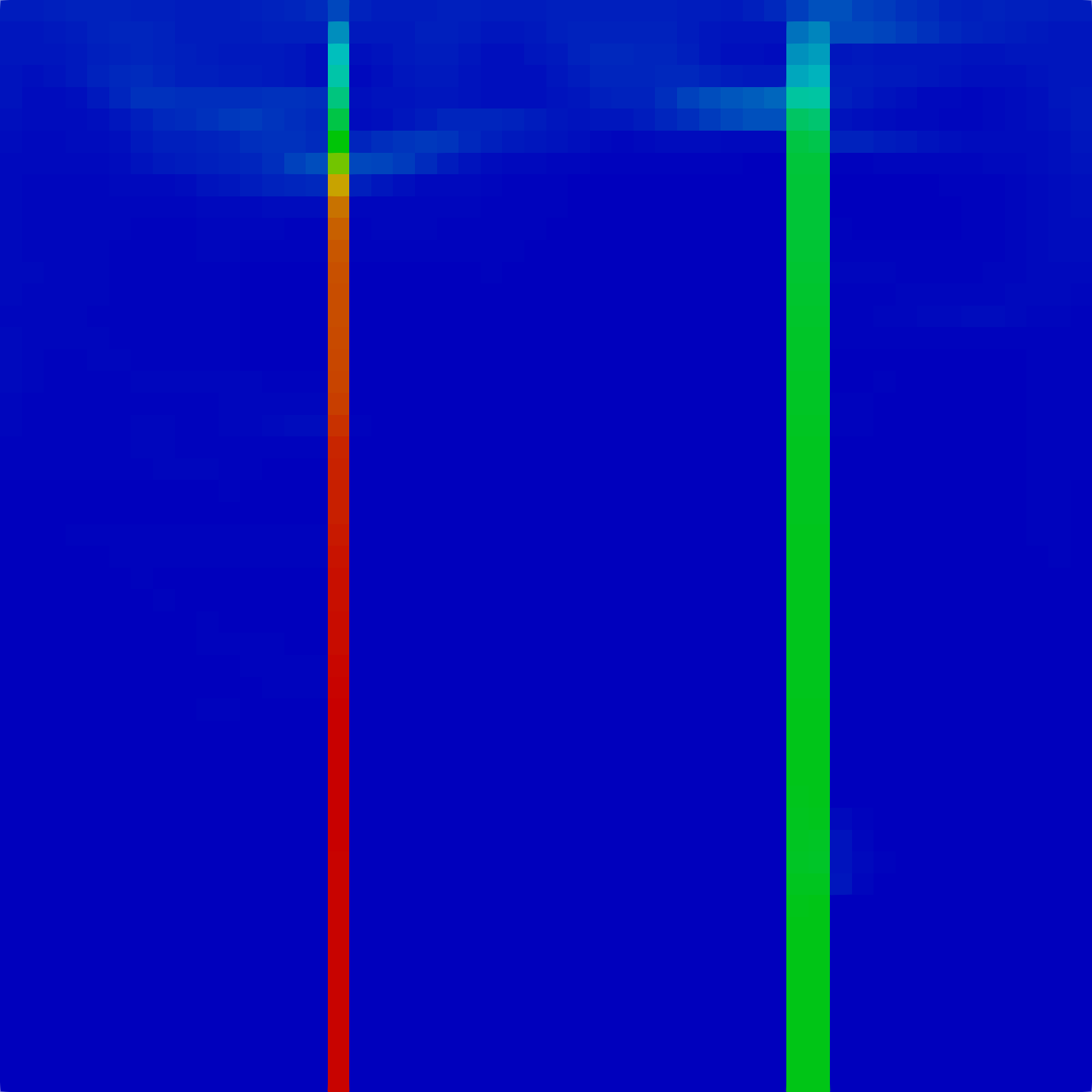}%
    \hspace*{0.0125\textwidth}%
    \includegraphics[scale=0.1]{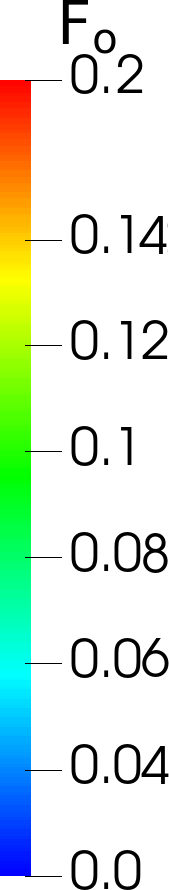}%
    \hspace*{0.02\textwidth}%
    \includegraphics[scale=0.075]{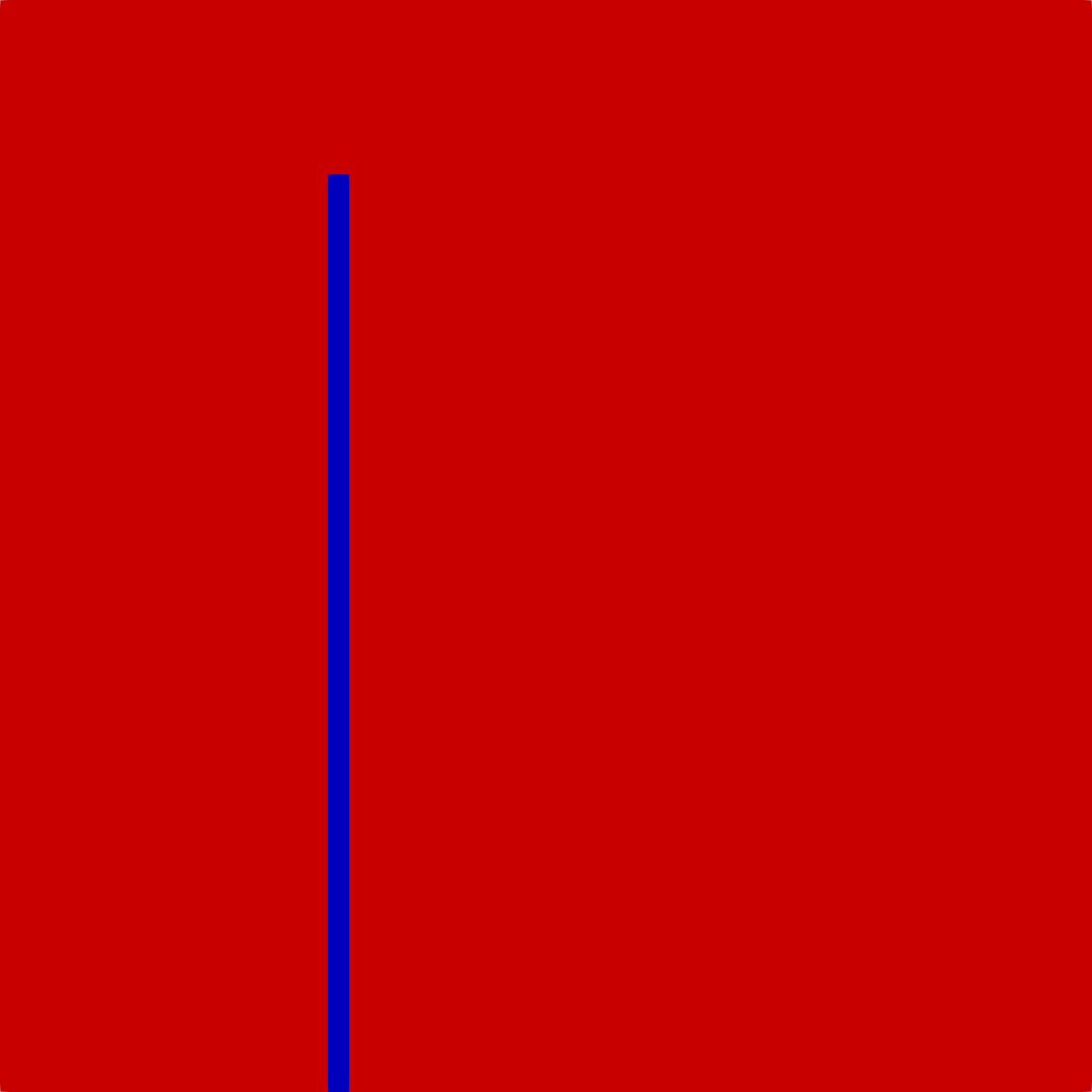}%
    \hspace*{0.0125\textwidth}%
    \includegraphics[scale=0.1]{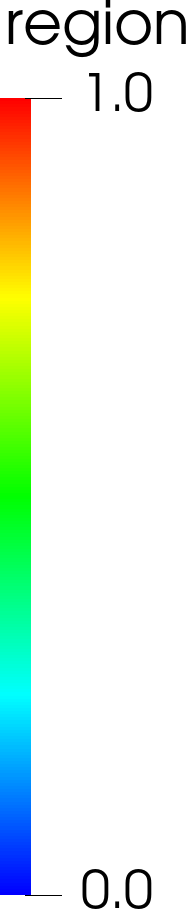}\\[0.2cm]
    \includegraphics[scale=0.075]{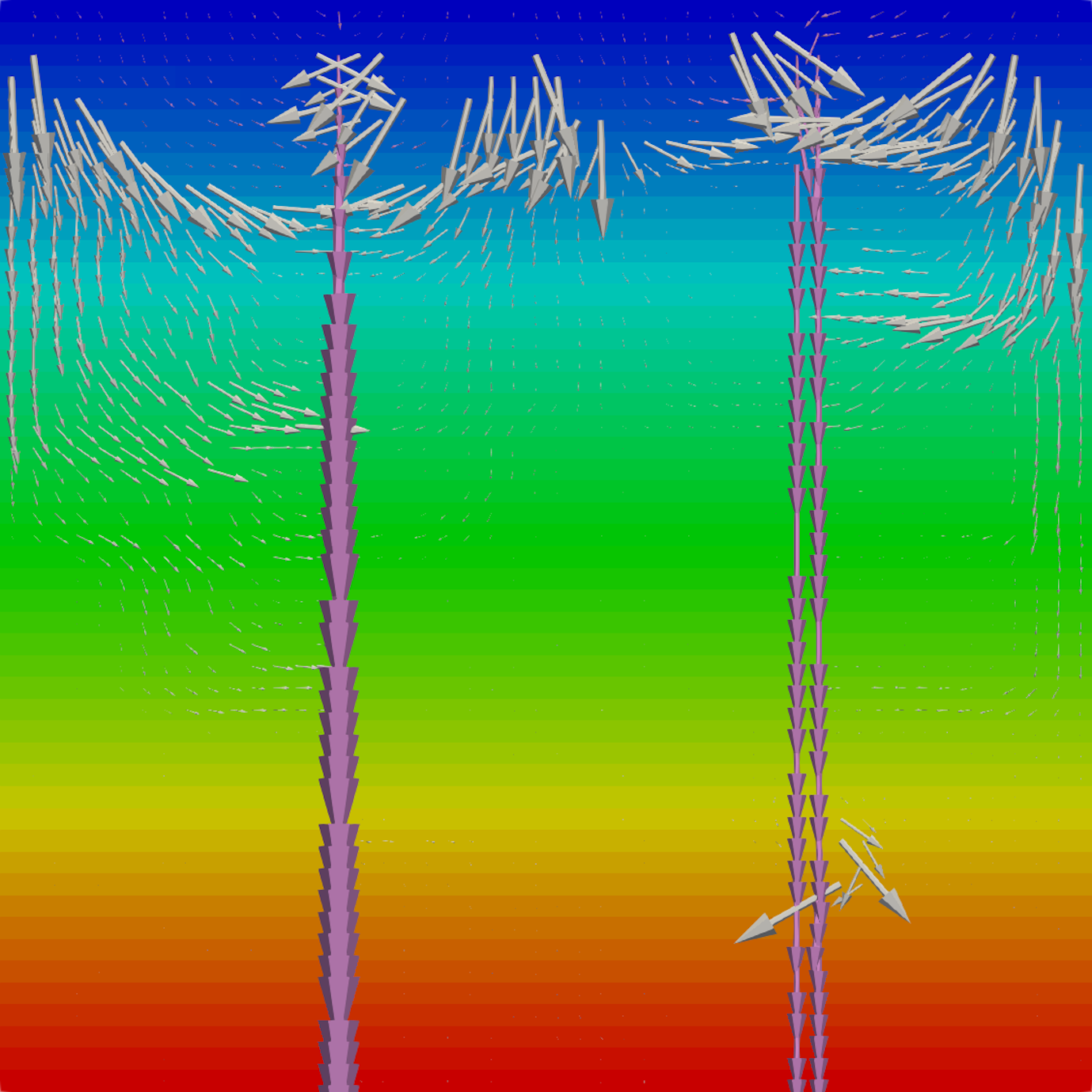}%
    \hspace*{0.0125\textwidth}%
    \includegraphics[scale=0.1]{figs/case1_adapt_pressure_legend.png}%
    \hspace*{0.02\textwidth}%
    \includegraphics[scale=0.075]{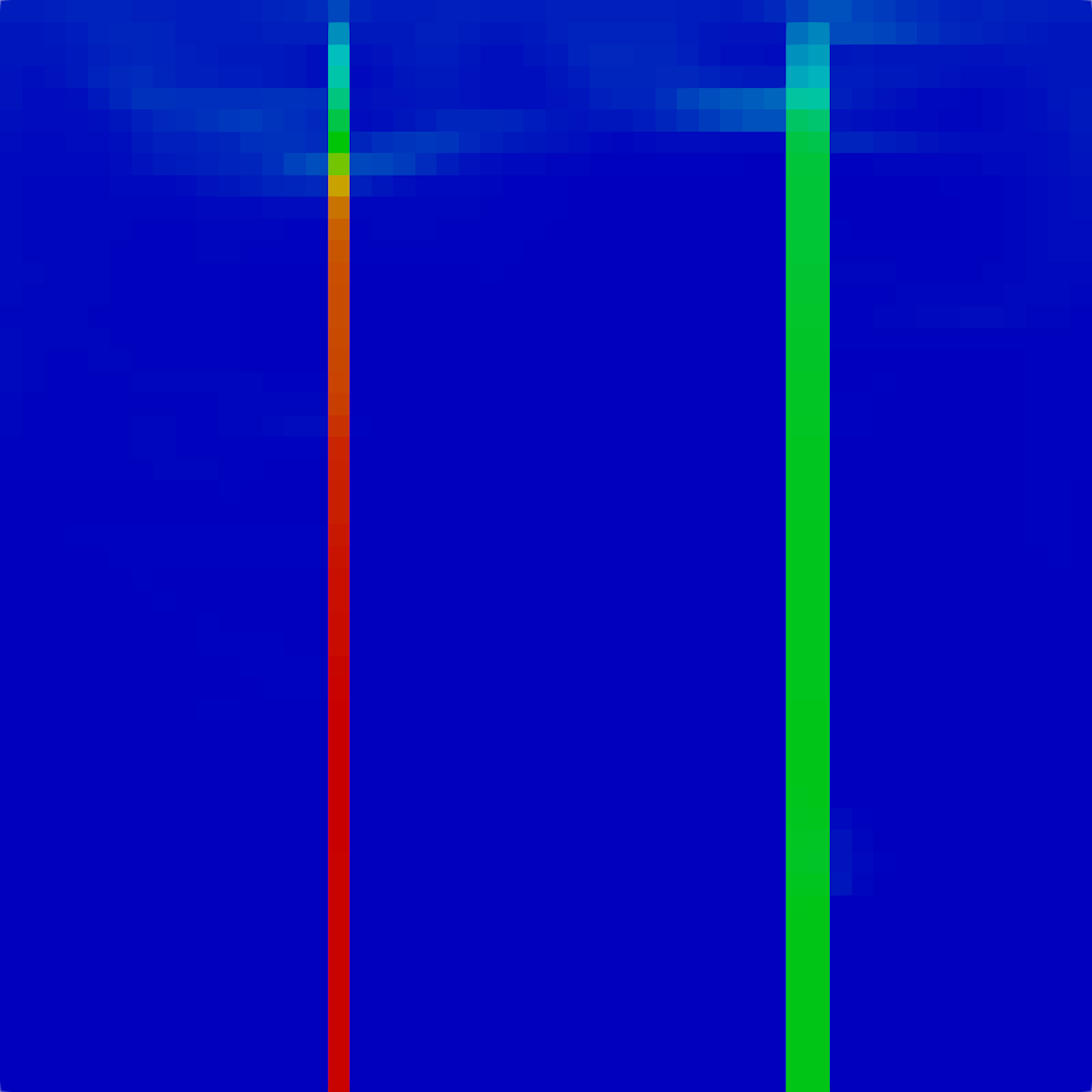}%
    \hspace*{0.0125\textwidth}%
    \includegraphics[scale=0.1]{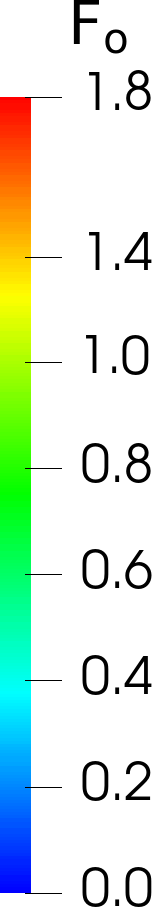}%
    \hspace*{0.02\textwidth}%
    \includegraphics[scale=0.075]{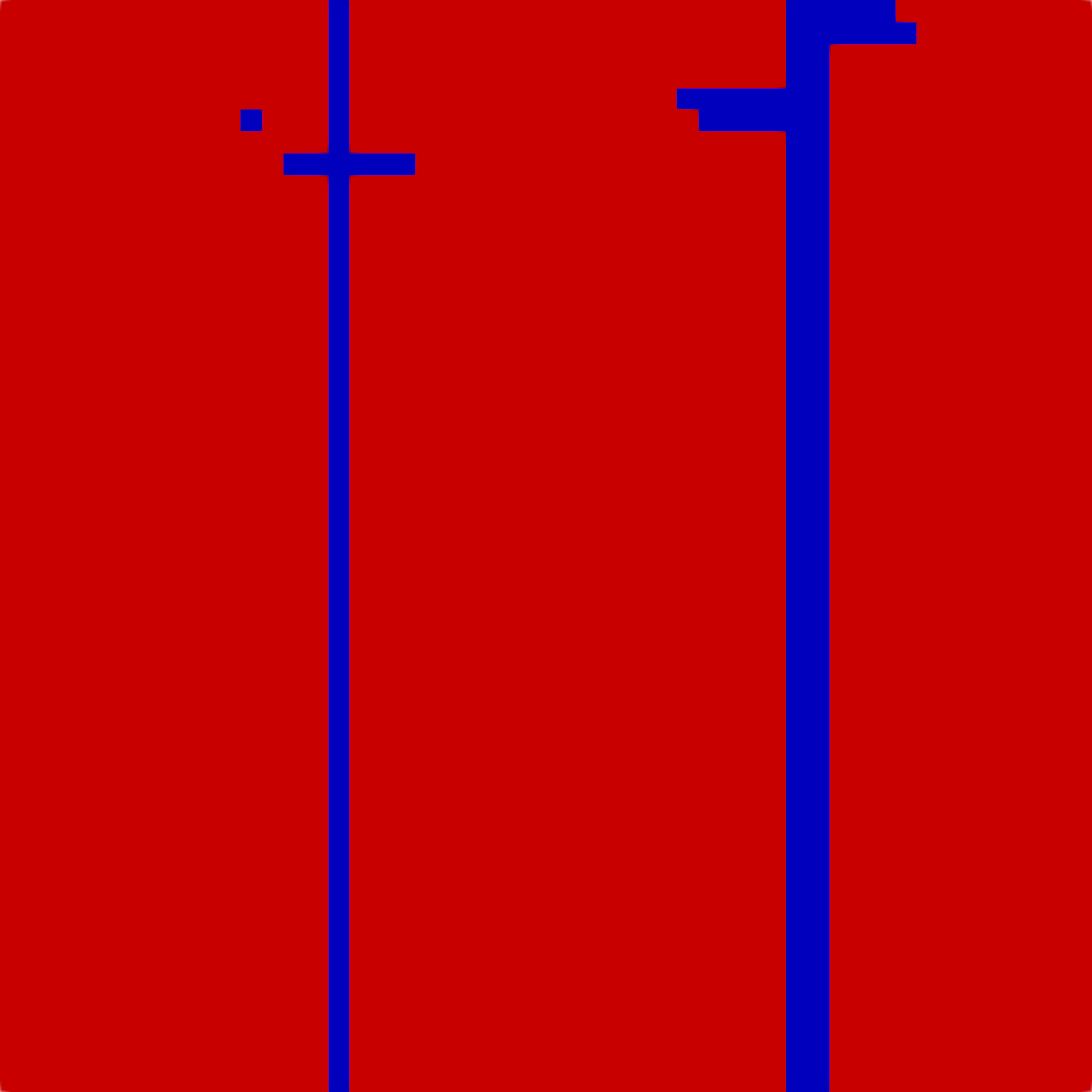}%
    \hspace*{0.0125\textwidth}%
    \includegraphics[scale=0.1]{figs/spe10_region_10_legend.png}
    \caption{Solutions obtained to the adaptive model for the example in Section \ref{subsubsub:two_channels}, in the top row for \textit{Scenario a} and in the bottom row for \textit{Scenario b}. The pressure is multiplied by $1000$, and the velocity arrows inside the channels are scaled by $1/6$ and $1/21$ in comparison with those outside, respectively for each scenario.}
    \label{fig:case1a}
\end{figure}

In Figure \ref{fig:case1a}, we report the solution obtained for both scenarios. We notice that for \textit{Scenario b}, as expected, the flow in the two channels is higher than in \textit{Scenario a}, resulting in a higher Forchheimer number so that more cells are selected in the Darcy--Forchheimer region, or subdomain (identified in the legend of the figure with the value $1$). For \textit{Scenario b}, both channels along with some extra cells are selected to be in the Darcy--Forchheimer region. For \textit{Scenario a}, only cells in the channels are selected; more precisely, all cells are selected in the left channel, possibly because it is narrower than the right one and thus develops a higher flow rate.

Table \ref{tab:case1_couple} shows the $L^2$ errors computed relative to a reference solution, namely, the solution to the global Darcy--Forchheimer model, that is, the Darcy--Forchheimer model applied in the whole domain, without any adaptivity involved. The errors are rather small, although we see that if a global Darcy model is used, we make an error that reaches several orders of magnitude more than with the adaptive model, especially for \textit{Scenario b}, even if the cells selected in the Darcy--Forchheimer subdomain in the adaptive model are many fewer than those in the Darcy subdomain. In both scenarios, the flux $\bu$ is affected by higher errors than the pressure $p$.

%all schemes they take 5 iterations, except the Darcy. The obtained errors are
\begin{table}[!ht]
    \centering
    \begin{tabular}{cccccccc}
                                                &          &
        \multicolumn{2}{c|}{Forchheimer region} &
        \multicolumn{2}{|c|}{Darcy region}      &
        \multicolumn{2}{|c}{Whole domain}                                                                         \\
        \cline{3-8}
                                                &
                                                & Darcy    & Adaptive
                                                & Darcy    & Adaptive
                                                & Darcy    & Adaptive
        \\
        \cline{3-8}
        \multirow{ 2}{*}{\textit{Sce. a}}       &
        $err_p$                                 & 5.01e-10 & 3.64e-12 & 4.31e-10 & 8.36e-11 & 4.33e-10 & 8.28e-11 \\
                                                &
        $err_{\bu}$                                 & 6.89e-06 & 1.47e-06 & 4.95e-06 & 2.61e-06 & 6.28e-06 & 1.95e-06 \\
        \cline{3-8}
        \multirow{ 2}{*}{\textit{Sce. b}}       &
        $err_p$                                 & 3.98e-08 & 1.70e-12 & 4.35e-08 & 9.68e-12 & 4.33e-08 & 9.37e-12 \\
                                                &
        $err_{\bu}$                                 & 6.20e-05 & 3.43e-08 & 1.27e-04 & 3.68e-07 & 6.28e-05 & 4.83e-08
    \end{tabular}
    \caption{$L^2$ errors of the adaptive and global Darcy solutions relative to the global Darcy--Forchheimer solutions, broken down in the Darcy--Forchheimer, Darcy and whole domains, for both scenarios of the example in Section \ref{subsubsub:two_channels}.}
    \label{tab:case1_couple}
\end{table}

In Figure \ref{fig:case1a_smaller_error}, we decrease the error tolerance $\delta$ for both scenarios. We notice that by setting lower and lower values of $\delta$, coherently, the Darcy--Forchheimer cells become more and more numerous. First, the two channels become entirely of Darcy--Forchheimer type; then, the top boundary condition and the heterogeneities of the background permeability start to play a role in determining where the rest of the Darcy--Forchheimer subdomain develops and thus which cells should be selected. By taking $\delta$ even smaller, we expect that all the cells will be selected to be of Darcy--Forchheimer type.

\begin{figure}[!ht]
    \centering
    \includegraphics[scale=0.0913]{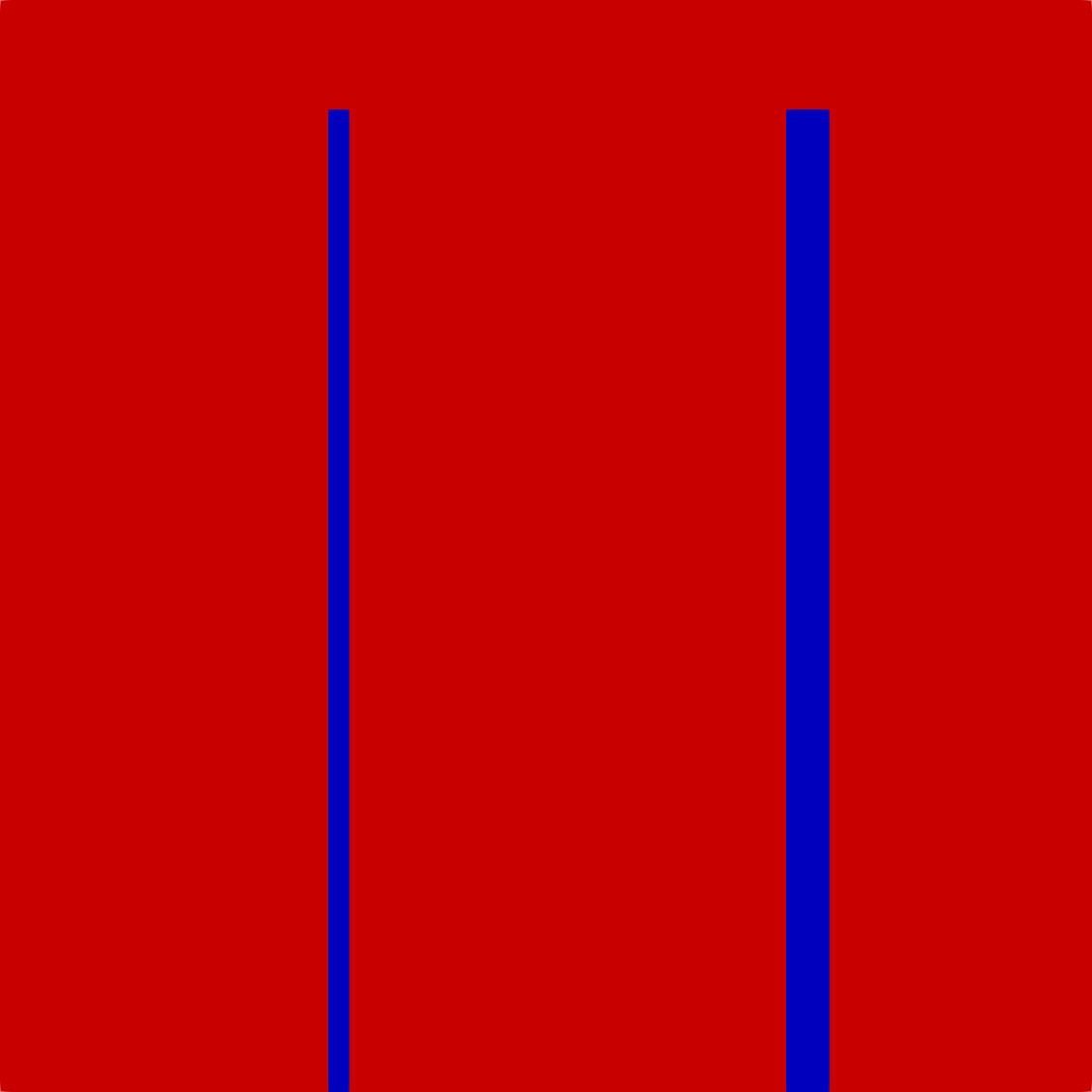}%
    \hspace*{0.02\textwidth}%
    \includegraphics[scale=0.0913]{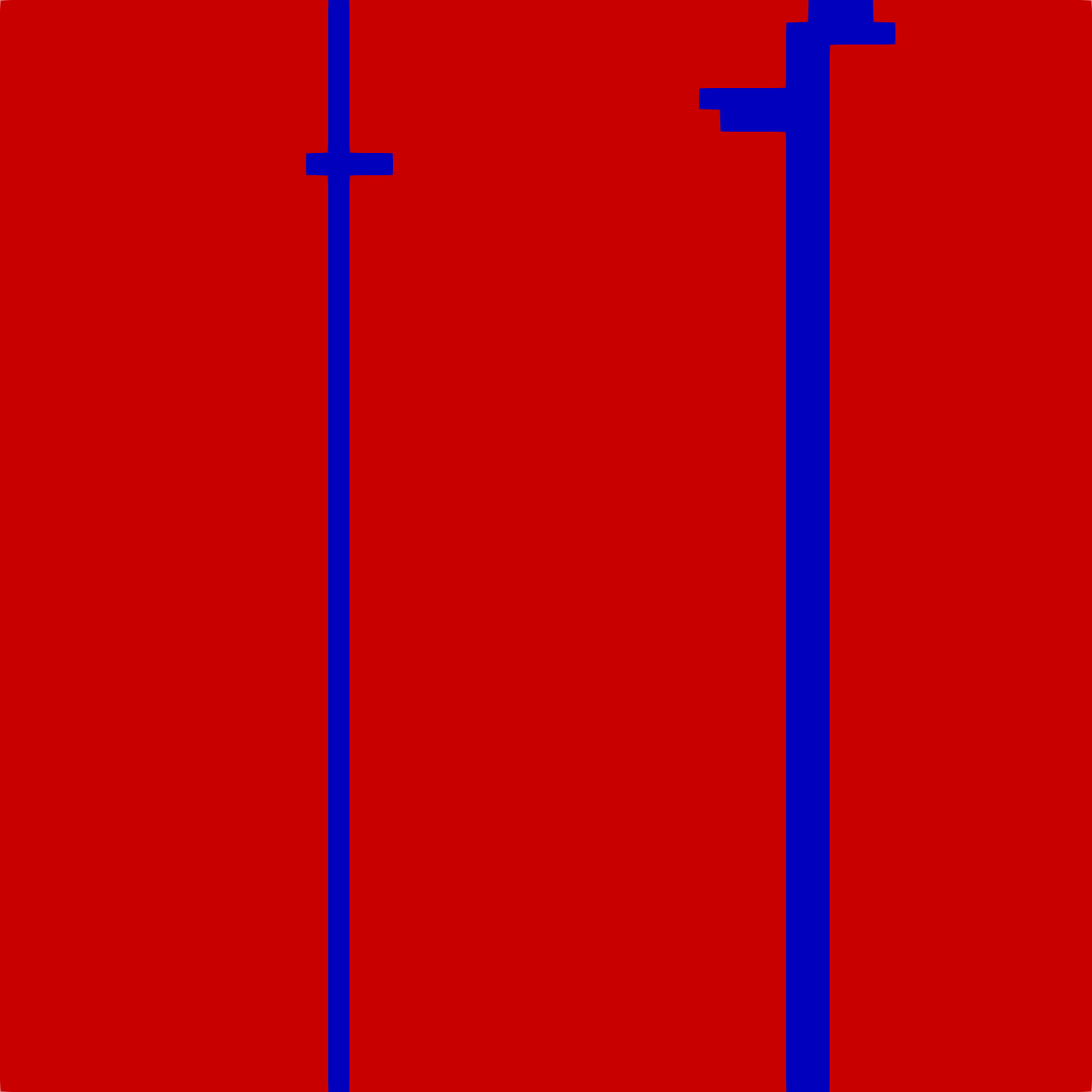}%
    \hspace*{0.02\textwidth}%
    \includegraphics[scale=0.0913]{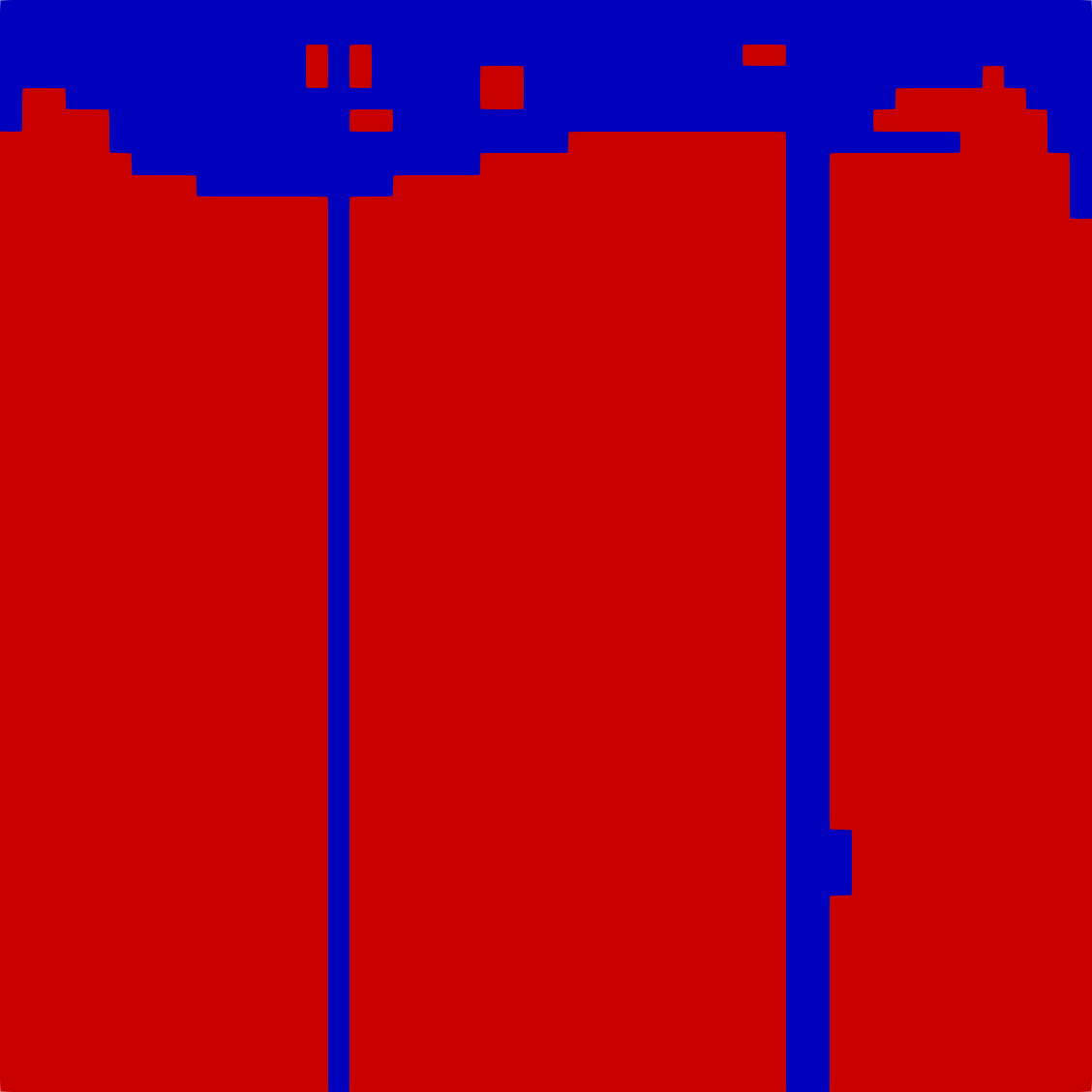}%
    \hspace*{0.0125\textwidth}%
    \includegraphics[scale=0.1]{figs/spe10_region_10_legend.png}\\[0.2cm]
    \includegraphics[scale=0.075]{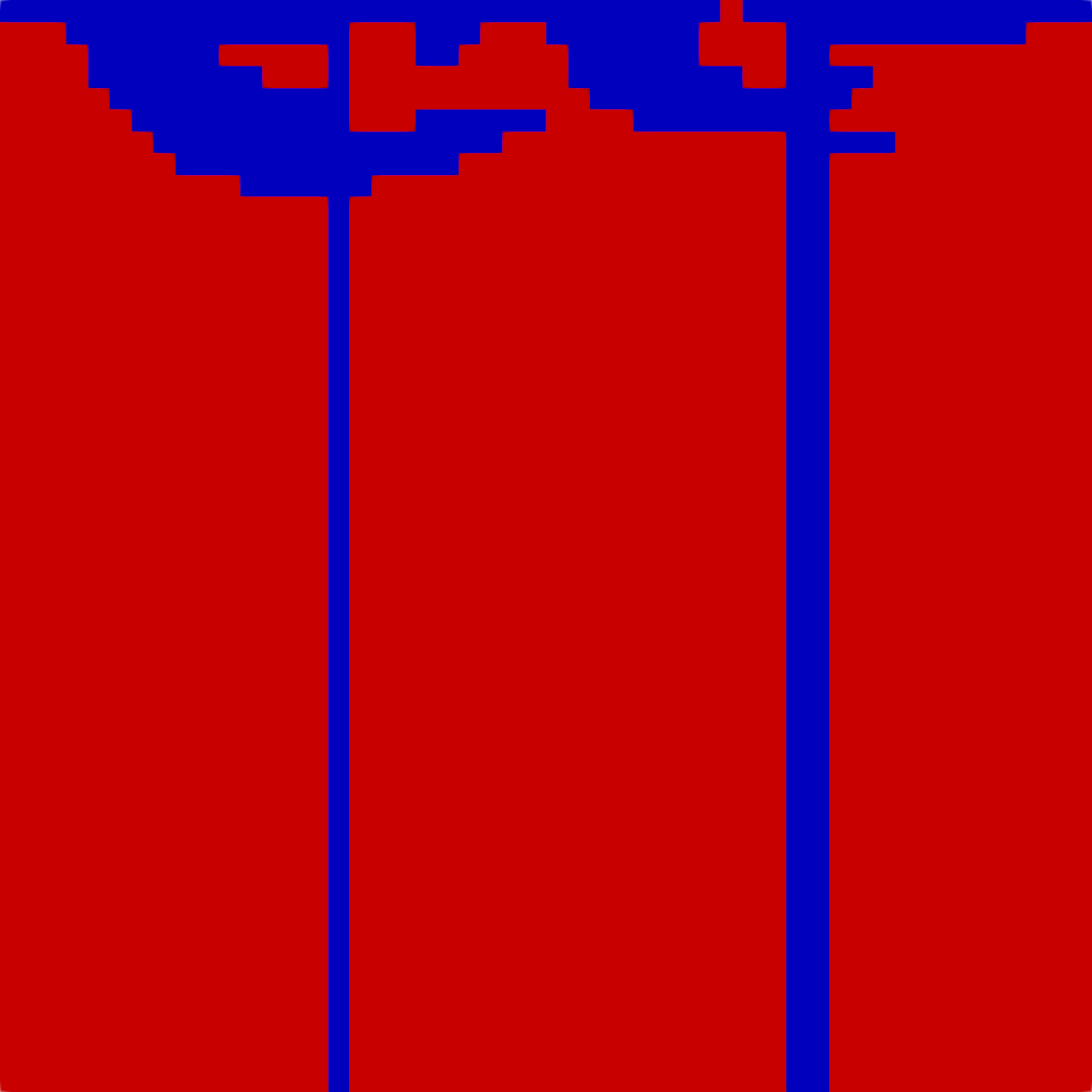}%
    \hspace*{0.02\textwidth}%
    \includegraphics[scale=0.075]{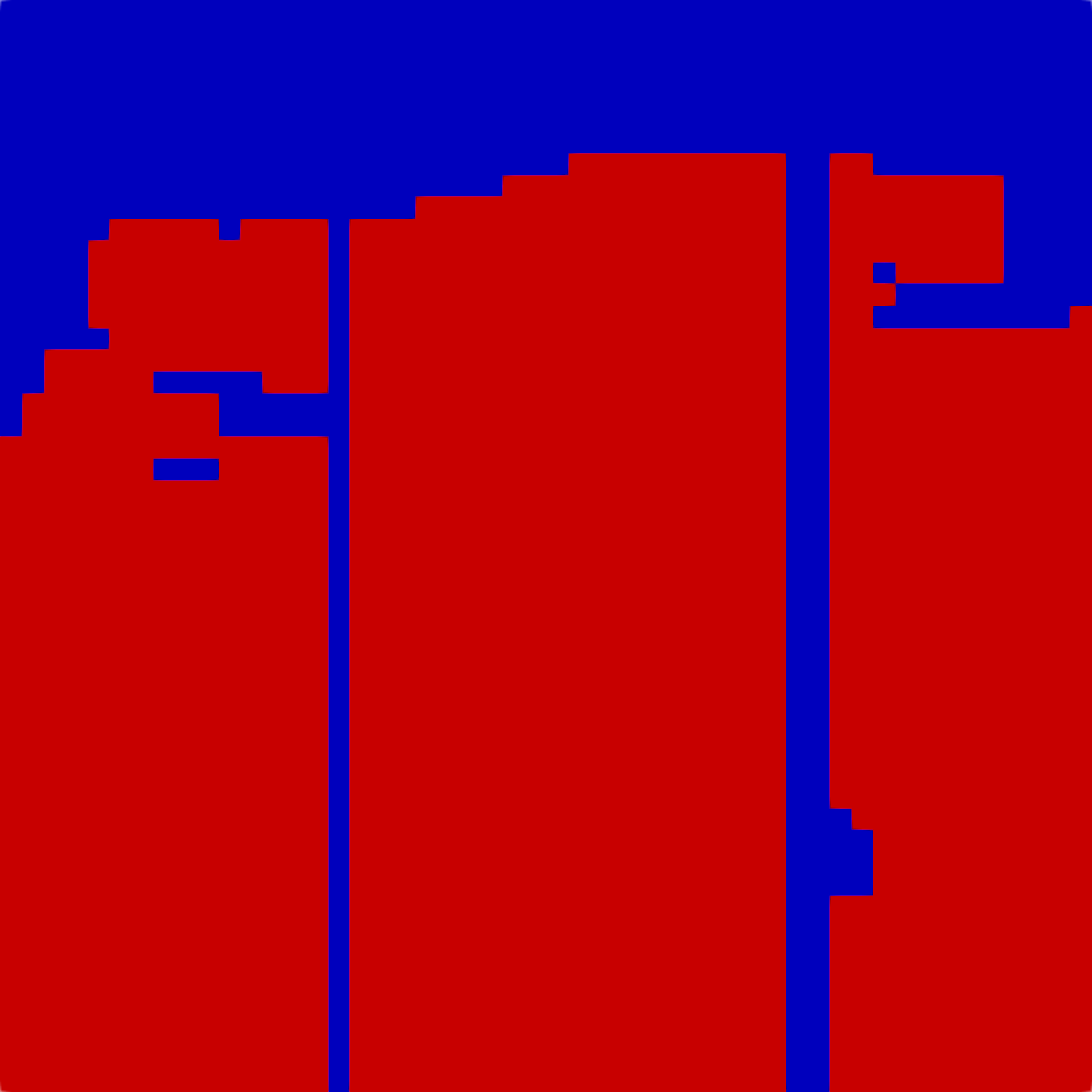}%
    \hspace*{0.02\textwidth}%
    \includegraphics[scale=0.075]{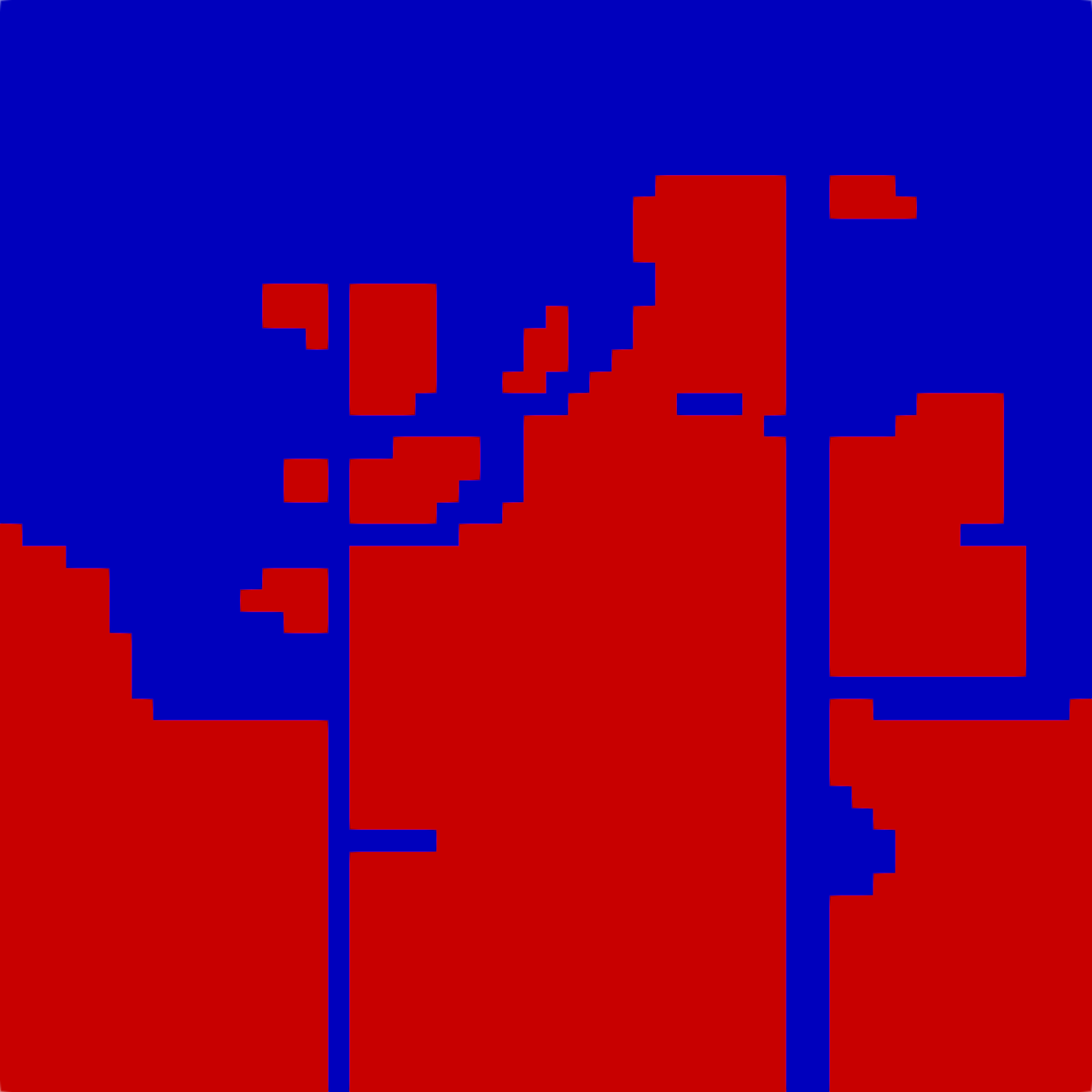}%
    \hspace*{0.0125\textwidth}%
    \includegraphics[scale=0.1]{figs/spe10_region_10_legend.png}
    \caption{Darcy (red) and Darcy--Forchheimer (blue) regions depending on the error tolerance $\delta$ decreasing from left to right according to $\delta=0.05$, $\delta=0.0125$ and $\delta=0.003125$, in the top row for \textit{Scenario a} and in the bottom row for \textit{Scenario b}, for the example in Section \ref{subsubsub:two_channels}.}
    \label{fig:case1a_smaller_error}
\end{figure}

\subsubsection{Network of channels} \label{subsubsub:network}

In this example, we consider the same background permeability as in the previous case, on which we place a network of intersecting horizontal and vertical channels with same width. All these channels are completely immersed in the domain without touching the boundary and they have permeabilities much higher than the surrounding medium. We expect, compared with the previous case, a more intricate flow distribution and thus a complex subdivision of cells into linear and nonlinear regions.

\begin{figure}[!ht]
    \centering
    \includegraphics[scale=0.075]{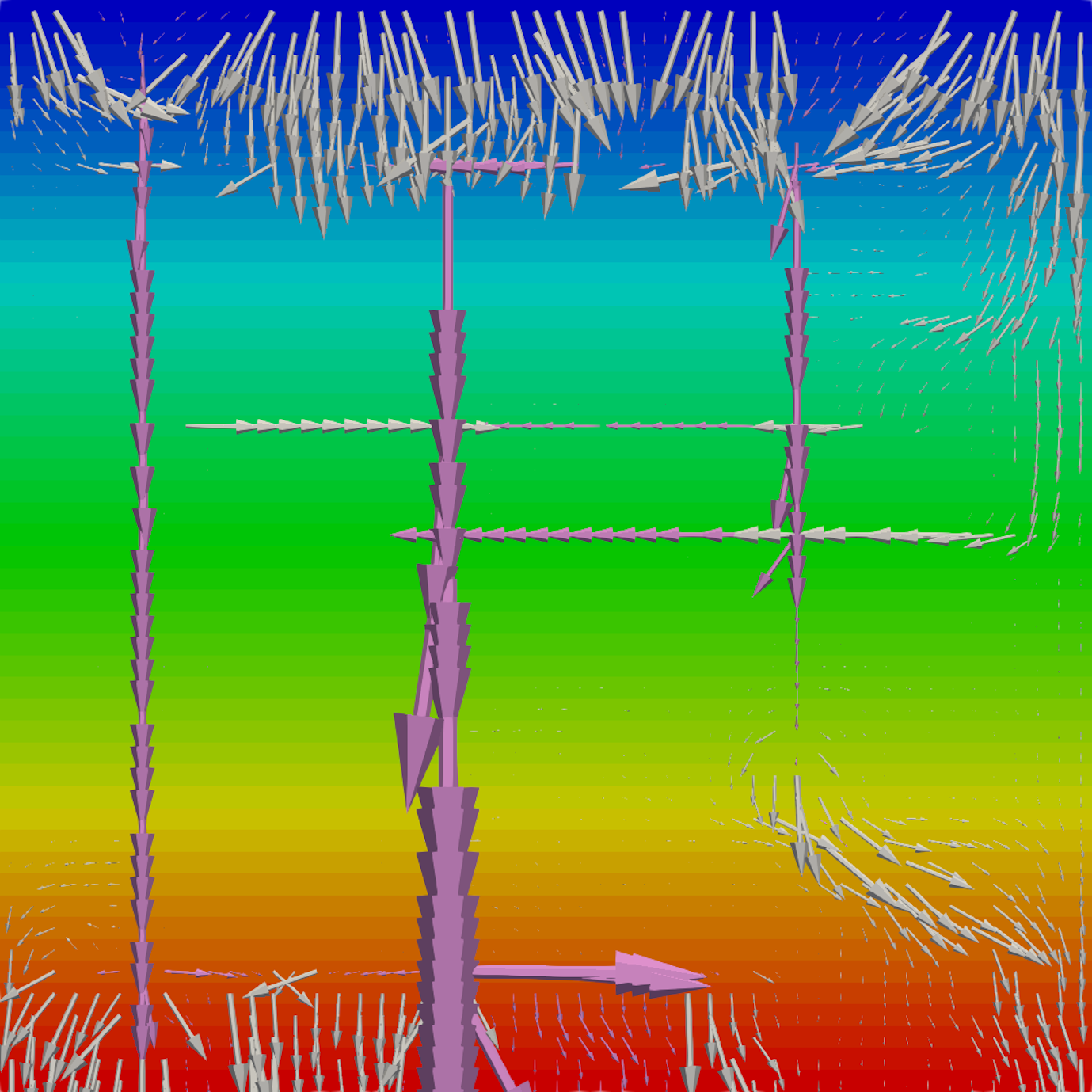}%
    \hspace*{0.0125\textwidth}%
    \includegraphics[scale=0.1]{figs/case1_adapt_pressure_legend.png}%
    \hspace*{0.02\textwidth}%
    \includegraphics[scale=0.075]{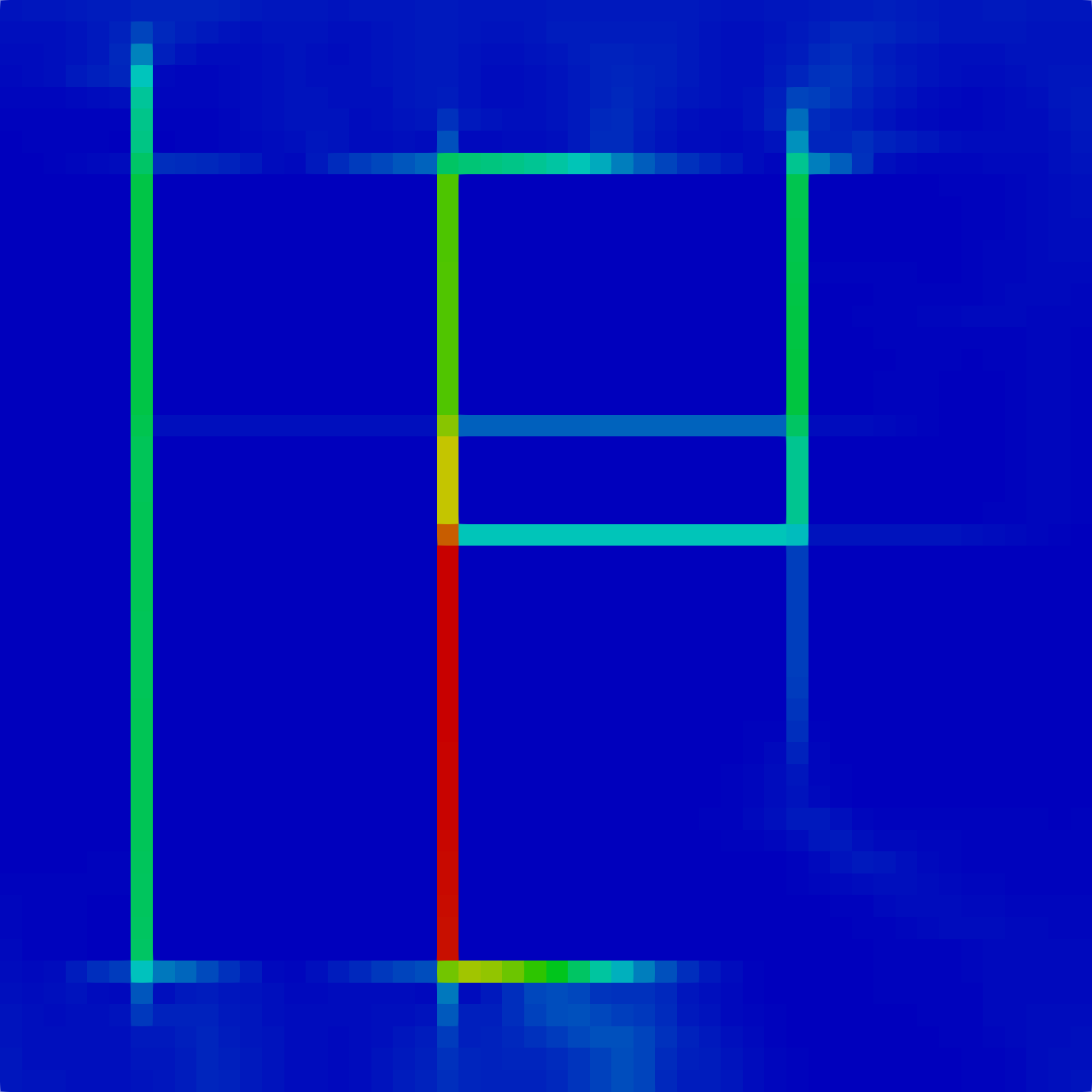}%
    \hspace*{0.0125\textwidth}%
    \includegraphics[scale=0.1]{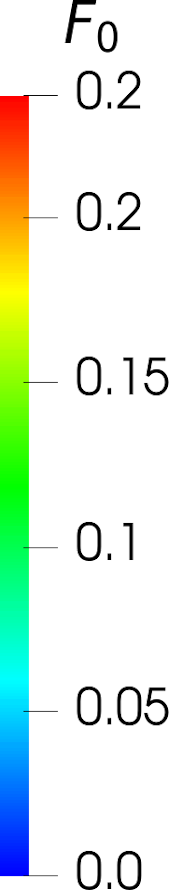}%
    \hspace*{0.02\textwidth}%
    \includegraphics[scale=0.075]{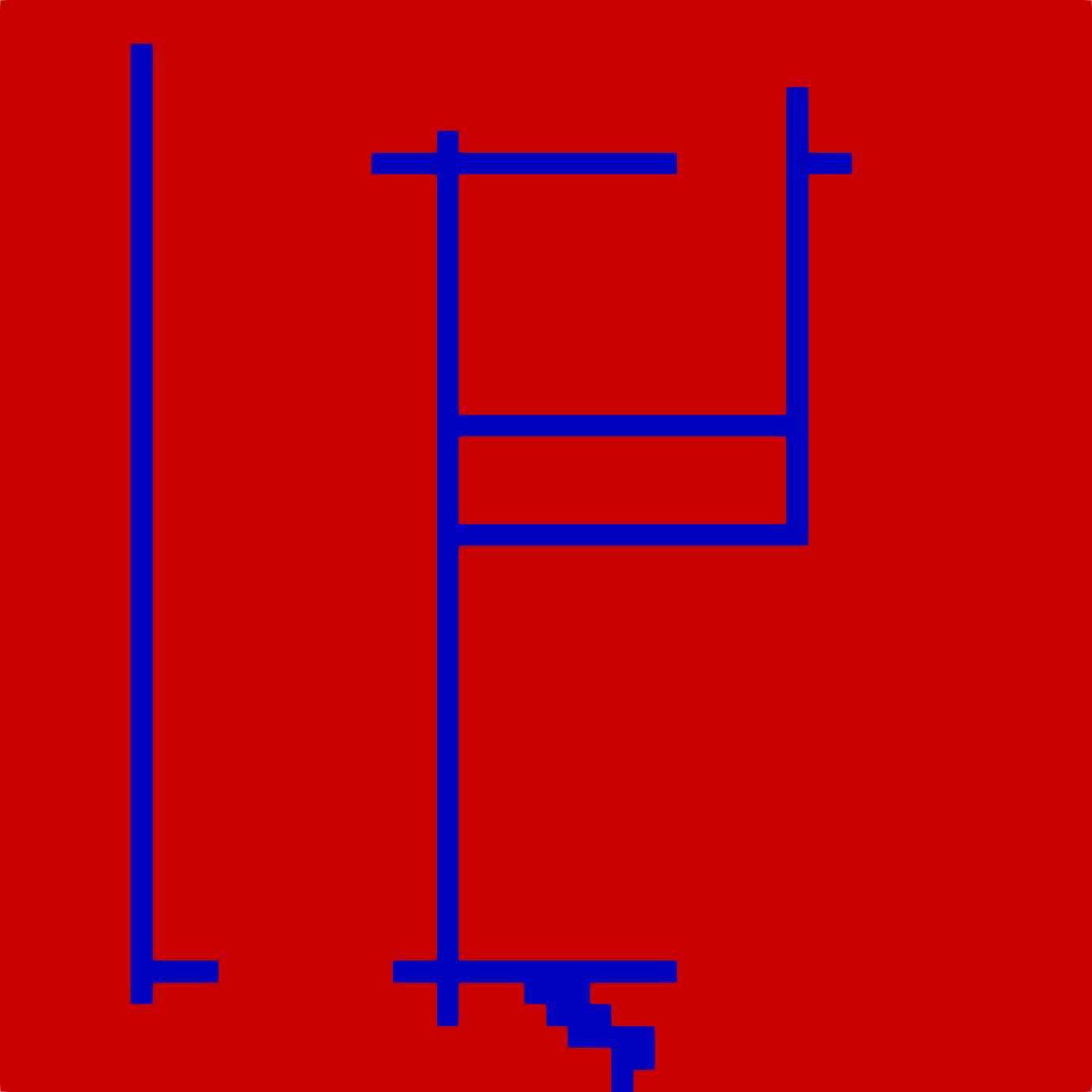}%
    \hspace*{0.0125\textwidth}%
    \includegraphics[scale=0.1]{figs/spe10_region_10_legend.png}\\[0.2cm]
    \includegraphics[scale=0.075]{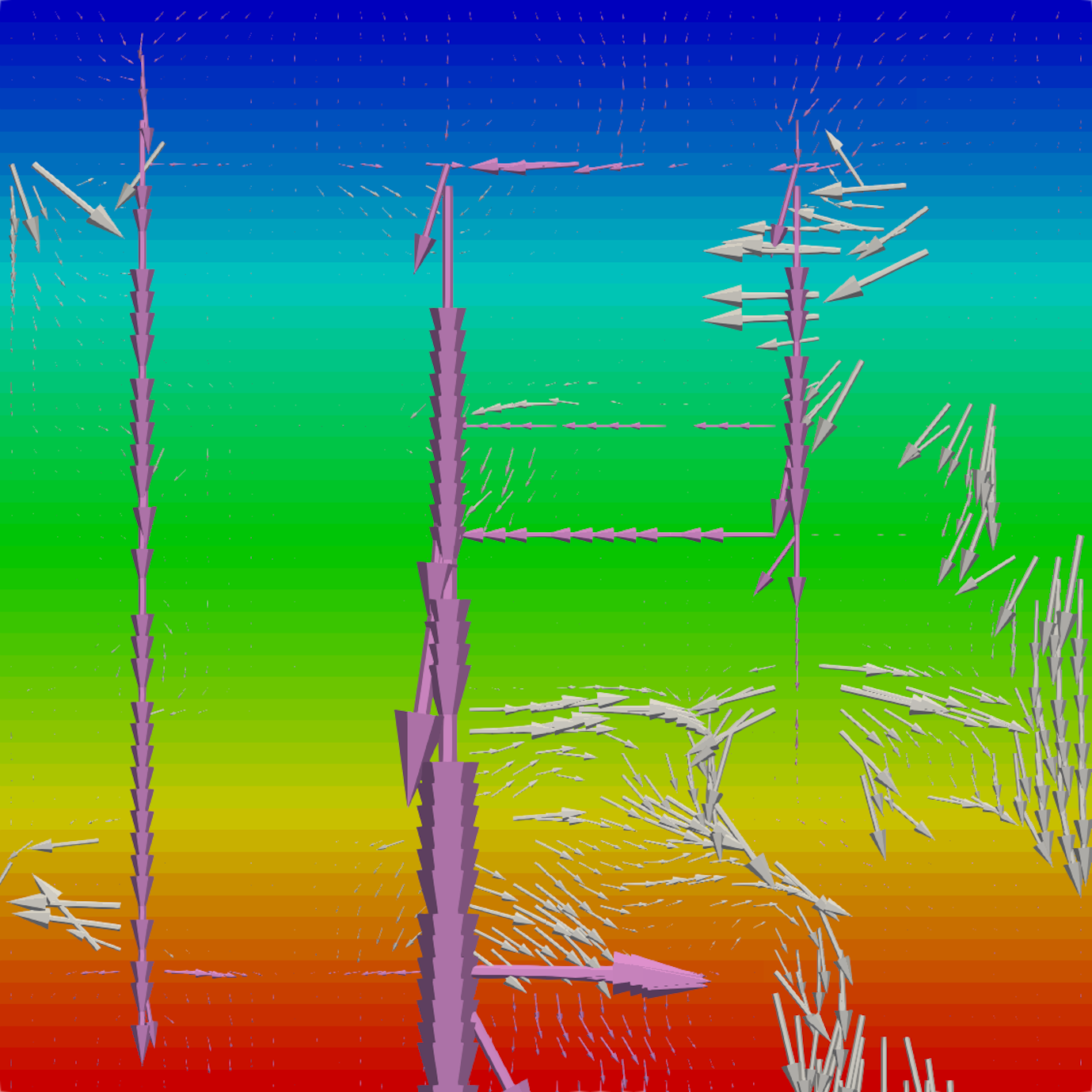}%
    \hspace*{0.0125\textwidth}%
    \includegraphics[scale=0.1]{figs/case1_adapt_pressure_legend.png}%
    \hspace*{0.02\textwidth}%
    \includegraphics[scale=0.075]{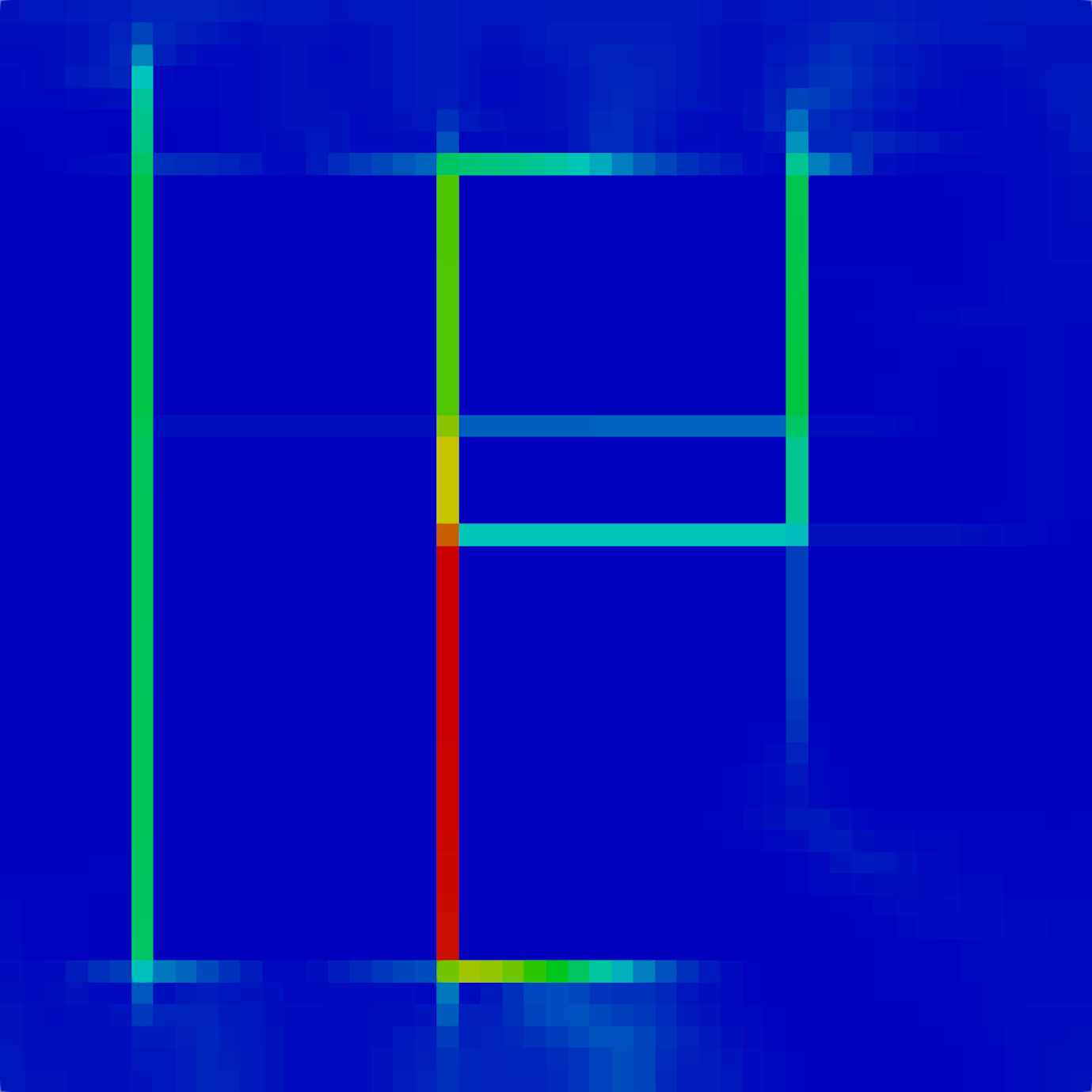}%
    \hspace*{0.0125\textwidth}%
    \includegraphics[scale=0.1]{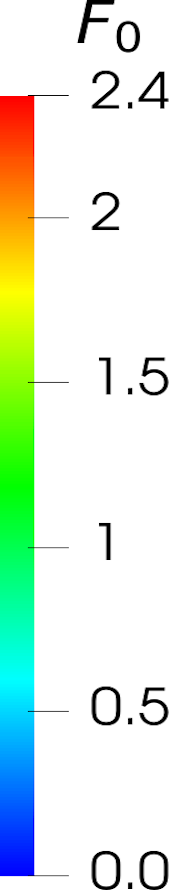}%
    \hspace*{0.02\textwidth}%
    \includegraphics[scale=0.075]{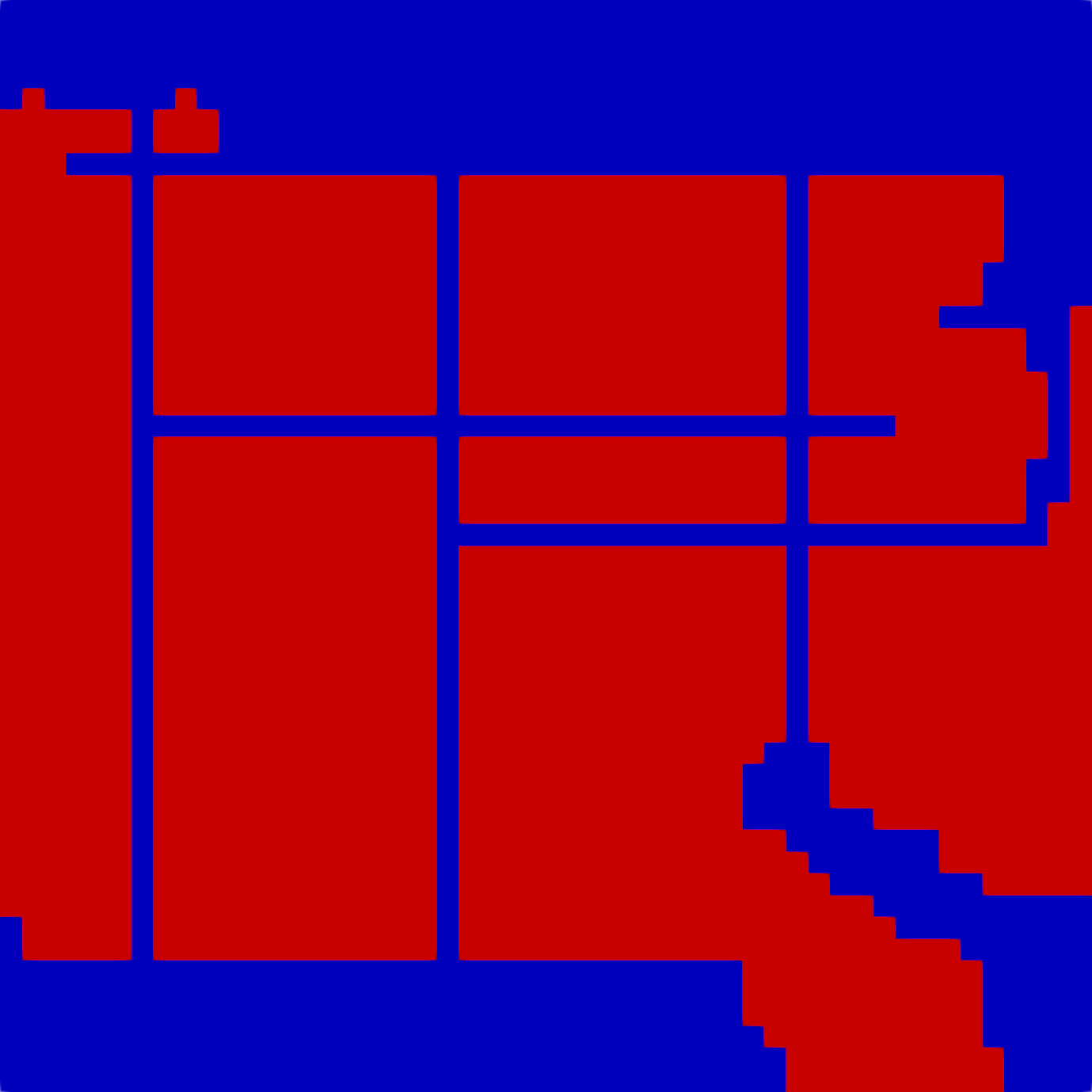}%
    \hspace*{0.0125\textwidth}%
    \includegraphics[scale=0.1]{figs/spe10_region_10_legend.png}
    \caption{Solutions obtained to the adaptive model for the example in Section \ref{subsubsub:network}, in the top row for \textit{Scenario a} and in the bottom row for \textit{Scenario b}. The pressure is multiplied by $1000$, and the velocity arrows inside the channels are scaled by $1/11$ and $1/120$ in comparison with those outside, respectively for each scenario.}
    \label{fig:case2a}
\end{figure}

In Figure \ref{fig:case2a}, we report the solutions for \textit{Scenario a} and \textit{Scenario b} by, recall, fixing the value of the error tolerance $\delta=0.1$. As expected, $\norm{\bu}$ for \textit{Scenario b} is much higher than for \textit{Scenario a}, which results in a higher Forchheimer number and consequently in more cells marked as Darcy--Forchheimer. We notice that, for \textit{Scenario a}, the Darcy--Forchheimer cells lie mostly inside the network, with only a few more cells leaking out of it and connecting with the bottom boundary. In contrast, for \textit{Scenario b}, two rather large regions, about the top and bottom of the domain, are marked as Darcy--Forchheimer and connect with the network.

Table \ref{tab:case1_network} shows the relative errors made by the global Darcy model and the adaptive model with respect to the global Darcy--Forchheimer model, considered thus to be the reference. As for the two-channel case, we notice that the adaptive scheme is more accurate, by several order of magnitude, than the global Darcy model. In particular, the detrimental impact of the error made in the cells marked as Darcy--Forchheimer brings down the error in the whole domain, even if these cells are many fewer than the Darcy cells. The adaptive approach gives thus rather accurate results.

\begin{table}[!ht]
    \centering
    \begin{tabular}{cccccccc}
                                                &          &
        \multicolumn{2}{c|}{Forchheimer region} &
        \multicolumn{2}{|c|}{Darcy region}      &
        \multicolumn{2}{|c}{Whole domain}                                                                         \\
        \cline{3-8}
                                                &
                                                & Darcy    & Adaptive
                                                & Darcy    & Adaptive
                                                & Darcy    & Adaptive
        \\
        \cline{3-8}
        \multirow{ 2}{*}{\textit{Sce. a}}       &
        $err_p$                                 & 6.08e-10 & 5.56e-13 & 6.03e-10 & 6.67e-13 & 6.03e-10 & 6.60e-13 \\
                                                &
        $err_{\bu}$                                 & 1.43e-04 & 8.02e-08 & 4.76e-04 & 5.27e-07 & 1.56e-04 & 1.07e-07 \\
        \cline{3-8}
        \multirow{ 2}{*}{\textit{Sce. b}}       &
        $err_p$                                 & 6.53e-08 & 1.85e-13 & 5.75e-08 & 1.24e-13 & 6.02e-08 & 1.47e-13 \\
                                                &
        $err_{\bu}$                                 & 1.55e-03 & 7.54e-10 & 9.02e-04 & 3.59e-08 & 1.55e-03 & 9.09e-10
    \end{tabular}
    \caption{$L^2$ errors of the adaptive and global Darcy solutions relative to the global Darcy--Forchheimer solutions, broken down in the Darcy--Forchheimer, Darcy and whole domains, for both scenarios of the example in Section \ref{subsubsub:network}.}
    \label{tab:case1_network}
\end{table}

In Figure \ref{fig:case1a_smaller_error}, we show the region split by decreasing the error tolerance $\delta$. As expected, the cells marked as Darcy--Forchheimer grow in number by filling up the zones near the top and bottom parts of the domain and then the areas in between the highly conductive channels.
\begin{figure}[!ht]
    \centering
    \includegraphics[scale=0.0913]{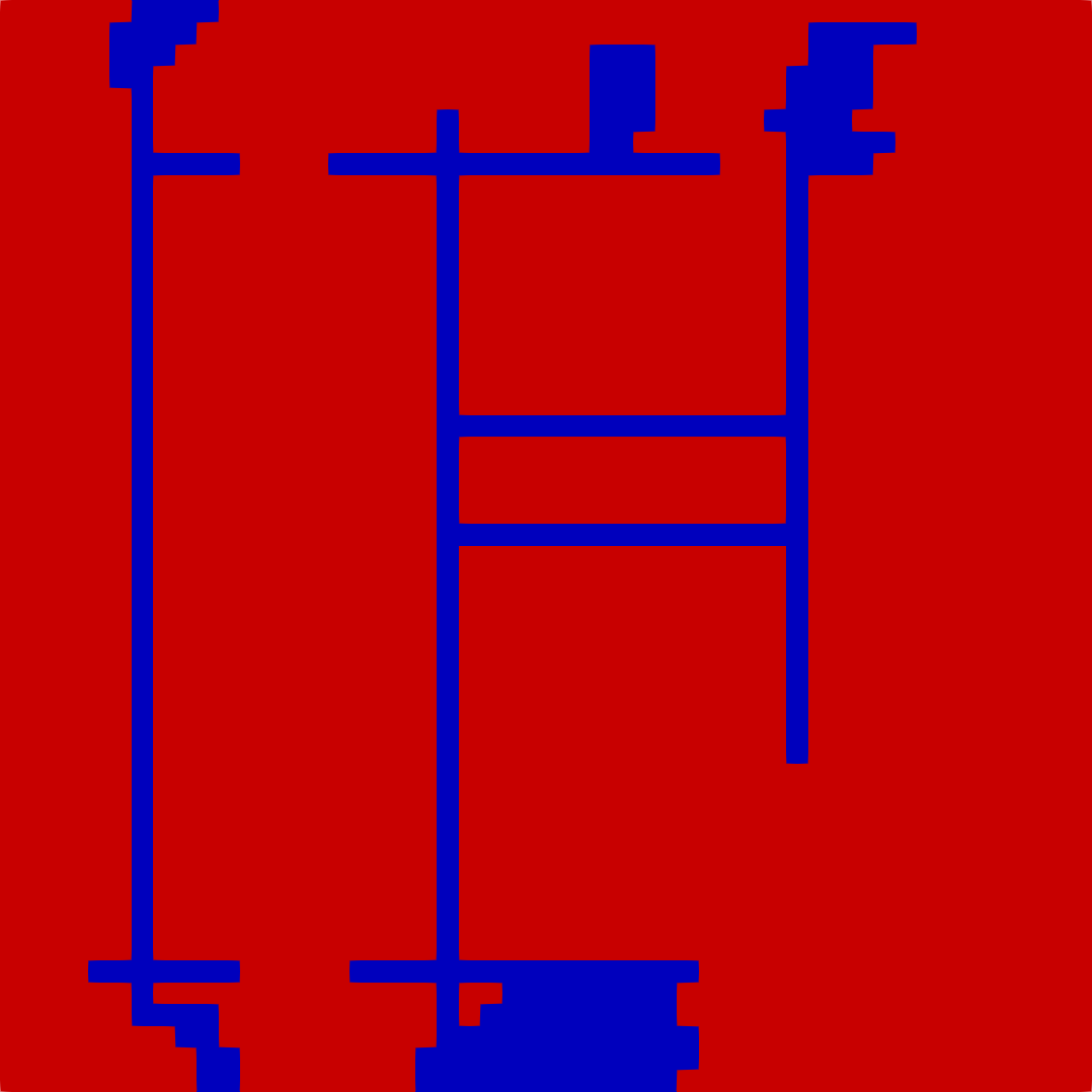}%
    \hspace*{0.02\textwidth}%
    \includegraphics[scale=0.0913]{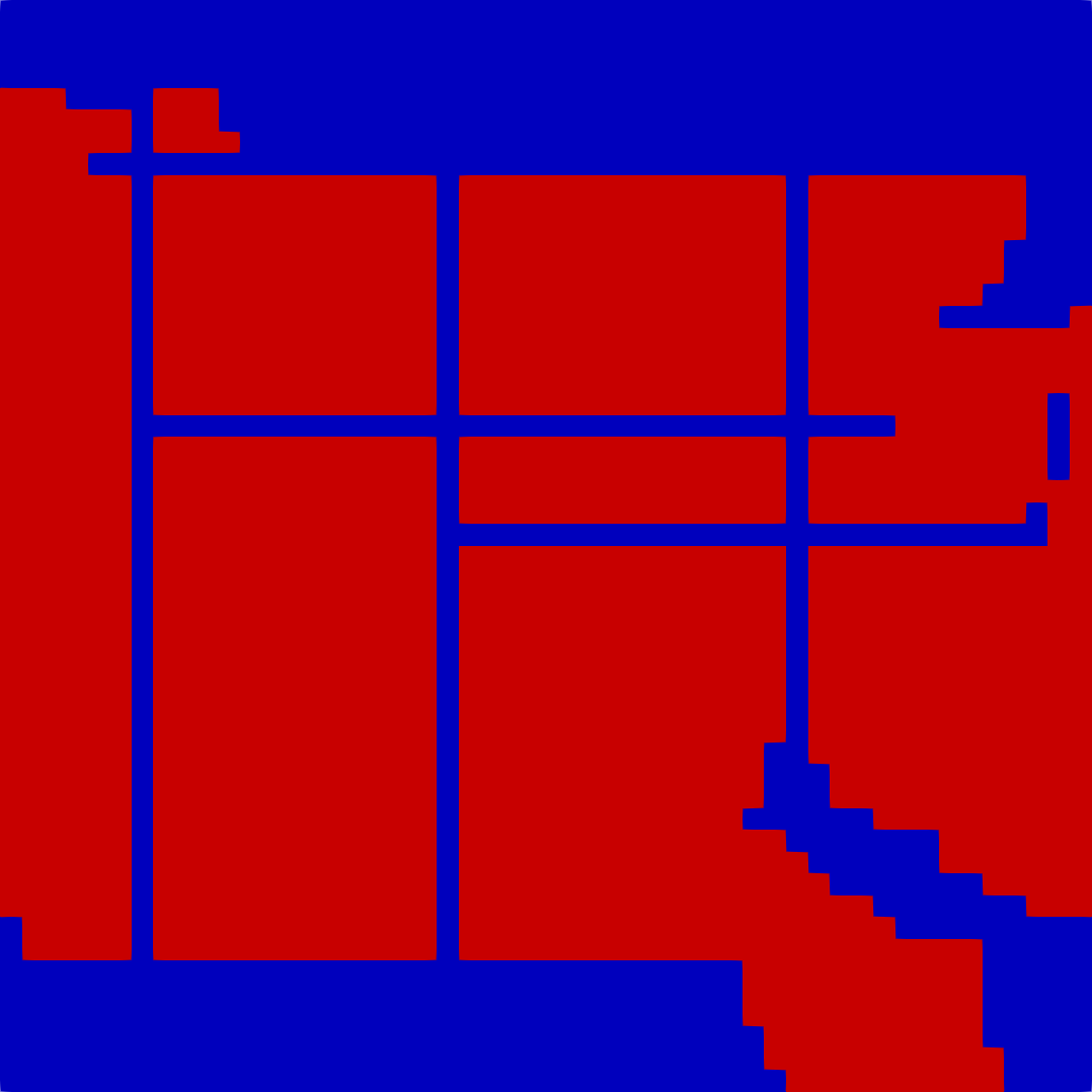}%
    \hspace*{0.02\textwidth}%
    \includegraphics[scale=0.0913]{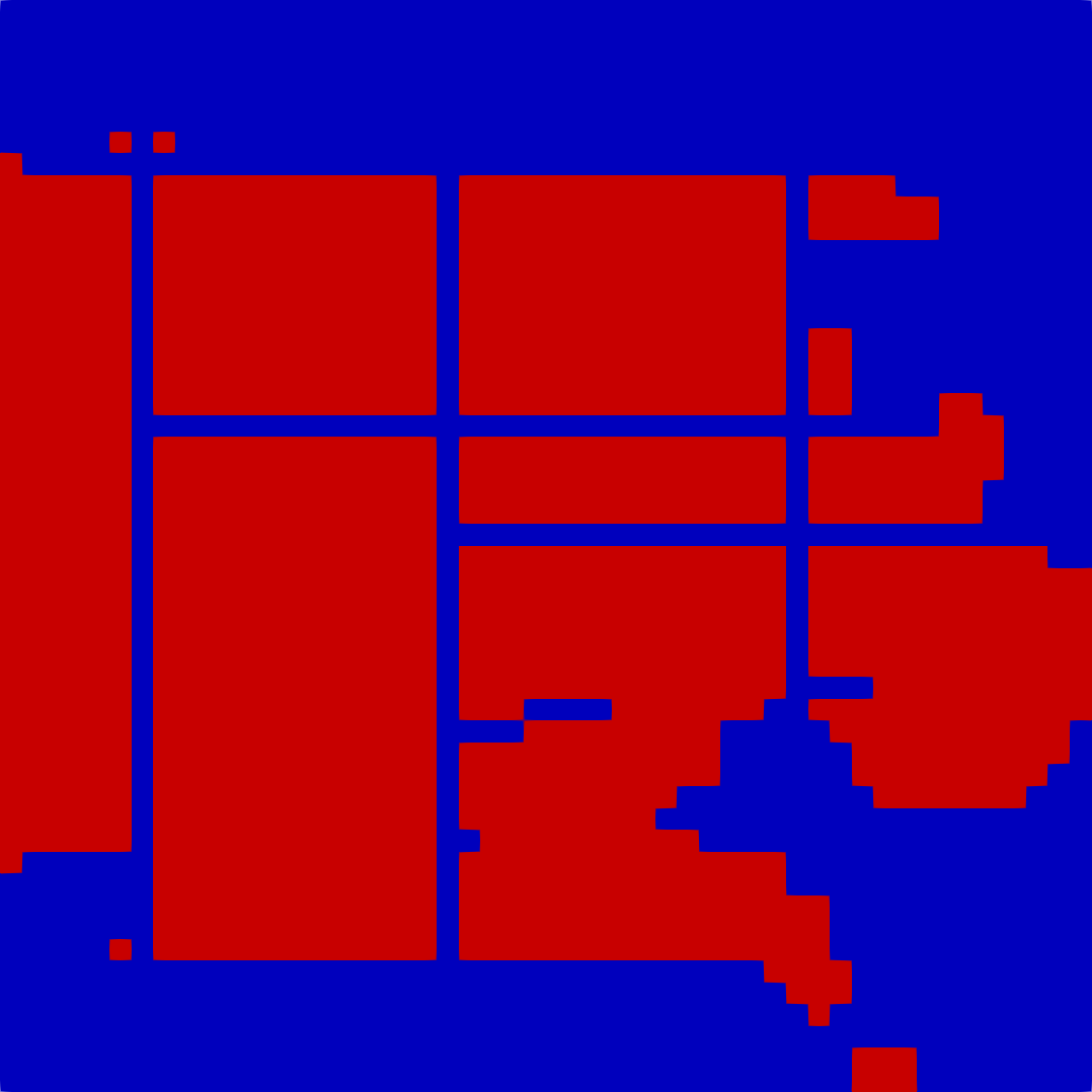}%
    \hspace*{0.0125\textwidth}%
    \includegraphics[scale=0.1]{figs/spe10_region_10_legend.png}\\[0.2cm]
    \includegraphics[scale=0.0913]{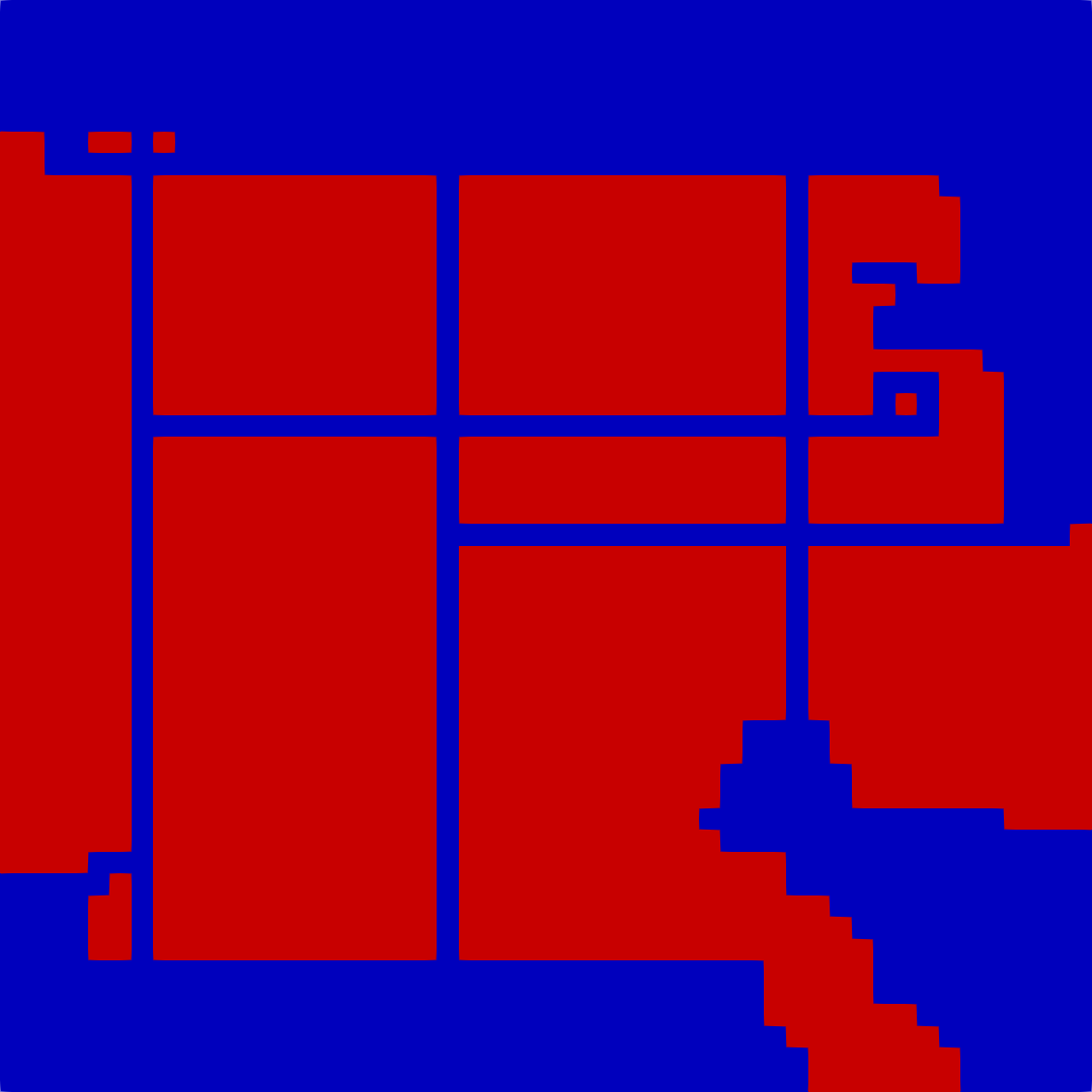}%
    \hspace*{0.02\textwidth}%
    \includegraphics[scale=0.0913]{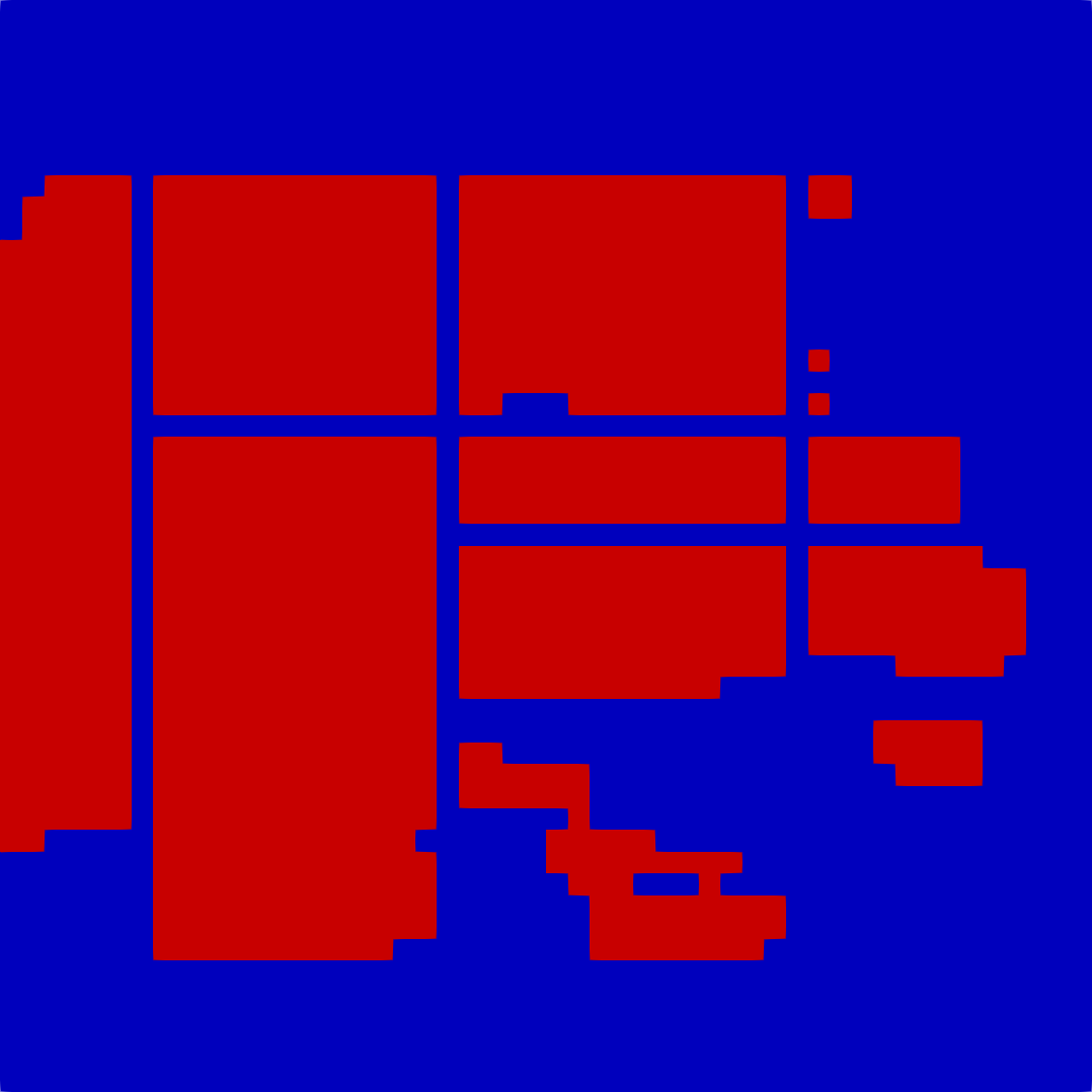}%
    \hspace*{0.02\textwidth}%
    \includegraphics[scale=0.0913]{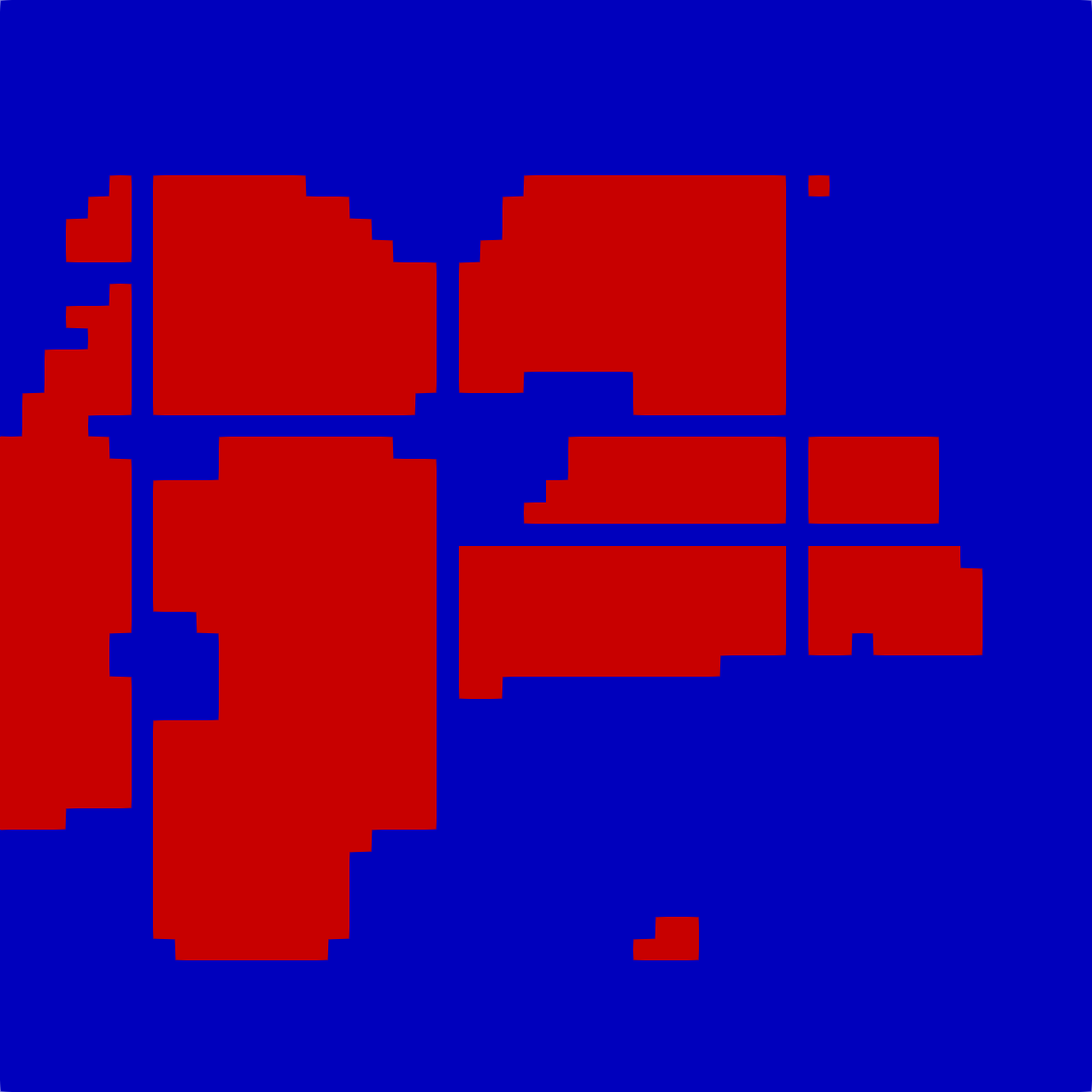}%
    \hspace*{0.0125\textwidth}%
    \includegraphics[scale=0.1]{figs/spe10_region_10_legend.png}
    \caption{Darcy (red) and Darcy--Forchheimer (blue) regions depending on the error tolerance $\delta$ decreasing from left to right according to $\delta=0.05$, $\delta=0.0125$ and $\delta=0.003125$, in the top row for \textit{Scenario a} and in the bottom row for \textit{Scenario b}, for the example in Section \ref{subsubsub:network}.}
    \label{fig:case1a_smaller_error}
\end{figure}

\subsubsection{Behavior with error tolerance}
\begin{figure}[!ht]
    \centering
    \includegraphics[width=0.75\textwidth]{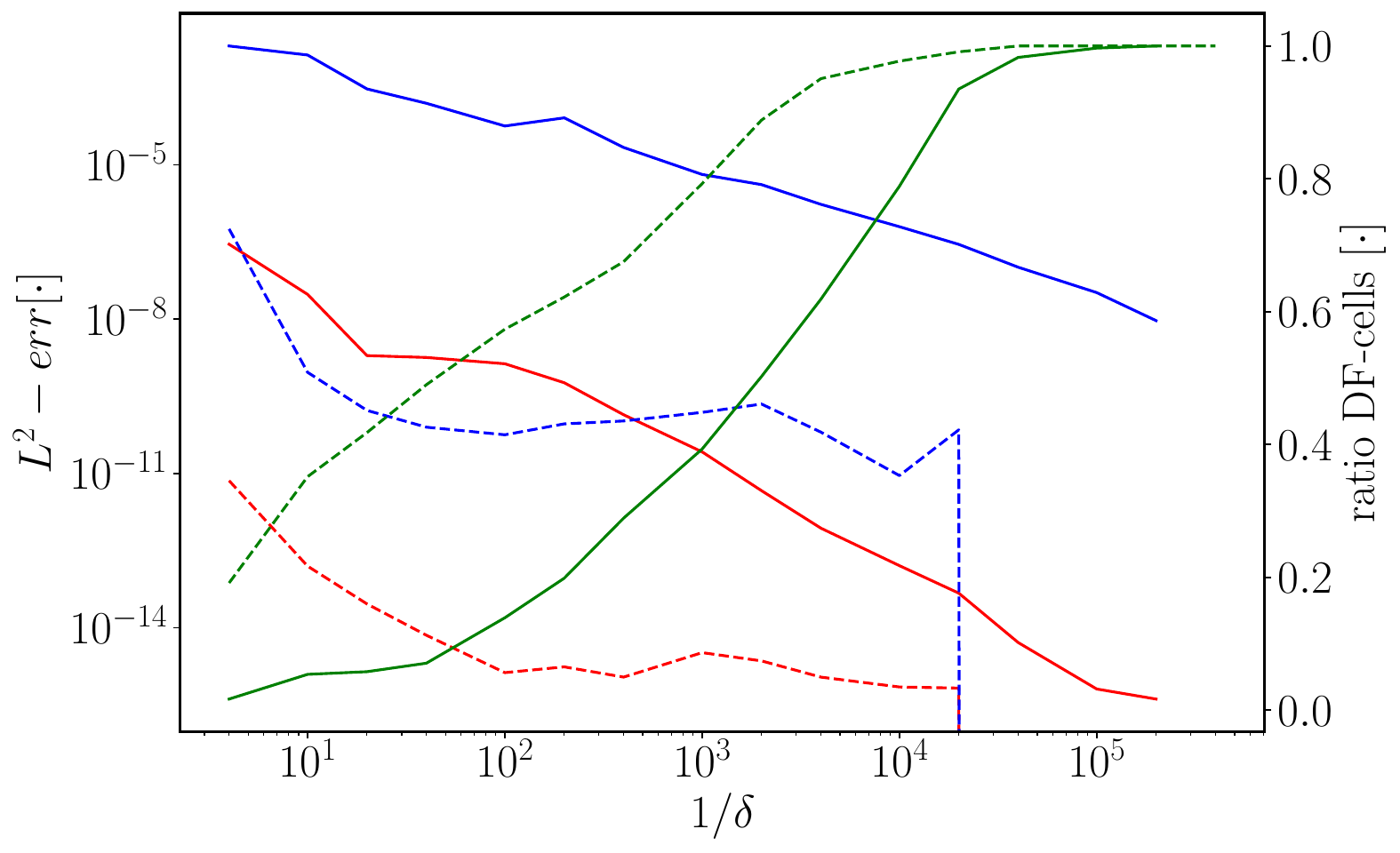}\\
    \includegraphics[width=0.75\textwidth]{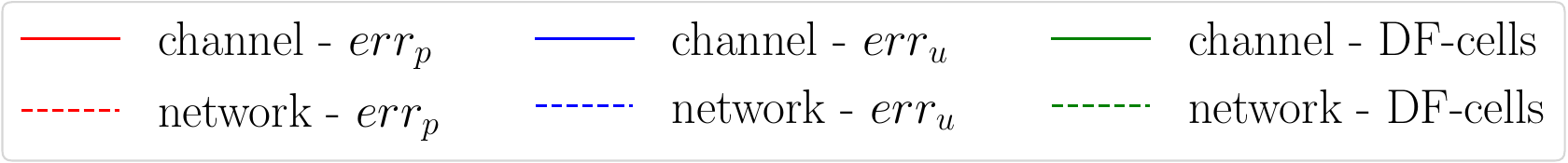}%
    \caption{Pressure (in red) and flux (in blue) $L^2$ relative errors by taking smaller values of the error tolerance $\delta$. The solid lines refer to \textit{Scenario a} and the dashed lines to \textit{Scenario b} for both the examples of Sections \ref{subsubsub:two_channels} (two-channel) and \ref{subsubsub:network} (network). In green, the ratio between the number of Darcy--Forchheimer cells and the total number of cells.}
    \label{fig:case2_error}
\end{figure}

For both the two-channel and network examples, we plot in Figure \ref{fig:case2_error} the errors obtained on the adaptive model relative the global Darcy--Forchheimer model as a function of the error tolerance $\delta$. The dashed lines for the relative errors, corresponding to \textit{Scenario b}, drop to zero starting from a certain value of $\delta$, since, from that value on, all cells in the adaptive approach are classified as Darcy--Forchheimer. This phenomenon is not seen in \textit{Scenario a} since the intensity of the top condition is one order of magnitude smaller than in the other scenario. In all examples and scenarios, the trend of the error is as expected: by decreasing $\delta$, the relative errors decrease and, simultaneously, the numbers of cells belonging to the Darcy--Forchheimer subdmain increase.

\subsection{One layer from SPE10} \label{subsec:case2}

\begin{figure}[!ht]
    \centering
    \includegraphics[scale=0.1]{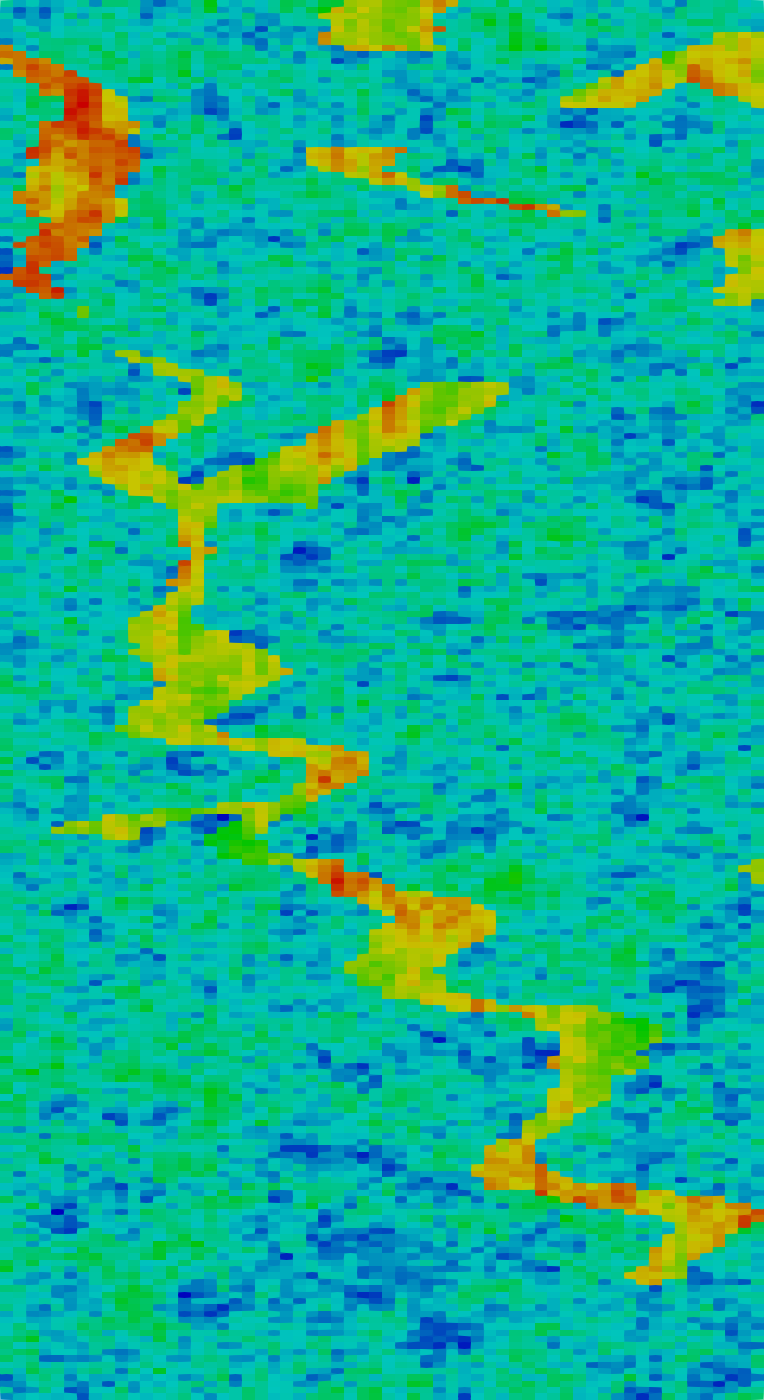}%
    \hspace*{0.025\textwidth}%
    \includegraphics[scale=0.15]{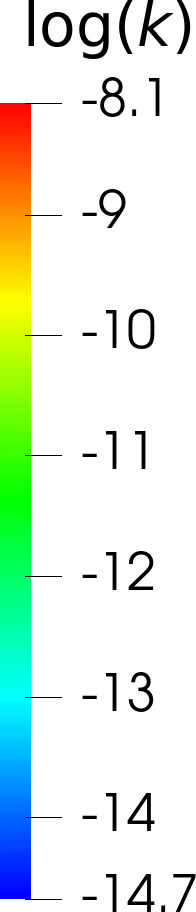}%
    \caption{Logarithm  of permeability for the example in Section \ref{subsec:case2}.}
    \label{fig:case2}
\end{figure}

In this example, we consider a single layer, namely, Layer 35, of the SPE10 benchmark study \cite{Christie2001} represented as a bidimensional domain. For the fluid, we consider the following physical parameters: the dynamic viscosity $\mu = 3\cdot 10^{-4} \sib{\pascal\cdot\second}$ and the density $\rho=1025\sib{\kilogram\per\cubic\meter}$. The permeability is given by the benchmark study and it is represented in Figure \ref{fig:case2}, where we notice some highly permeable paths surrounded by areas of low permeability.

As in the previous examples, the error tolerance is by default set to $\delta = 0.1$; the nonlinear exponent $m=1$, so that we are considering the Darcy--Forchheimer law, and the Forchheimer coefficient is set as $c_F = 0.55\sib{-}$.

We consider zero flux all around the boundary; moreover, we set five wells: one injector in the center of the domain and four identical producers at the corners. We consider different scenarios depending on the flux imposed in the injection well: in \textit{Scenario a}, we set $u_{\mt{in}} = 10 \sib{\kilogram\per\second}$; in \textit{Scenario b}, $u_{\mt{in}} = 50 \sib{\kilogram\per\second}$; in \textit{Scenario c}, $u_{\mt{in}} = 200 \sib{\kilogram\per\second}$. 

%(caso1 git)

\begin{figure}[!ht]
    \centering
    \includegraphics[scale=0.075]{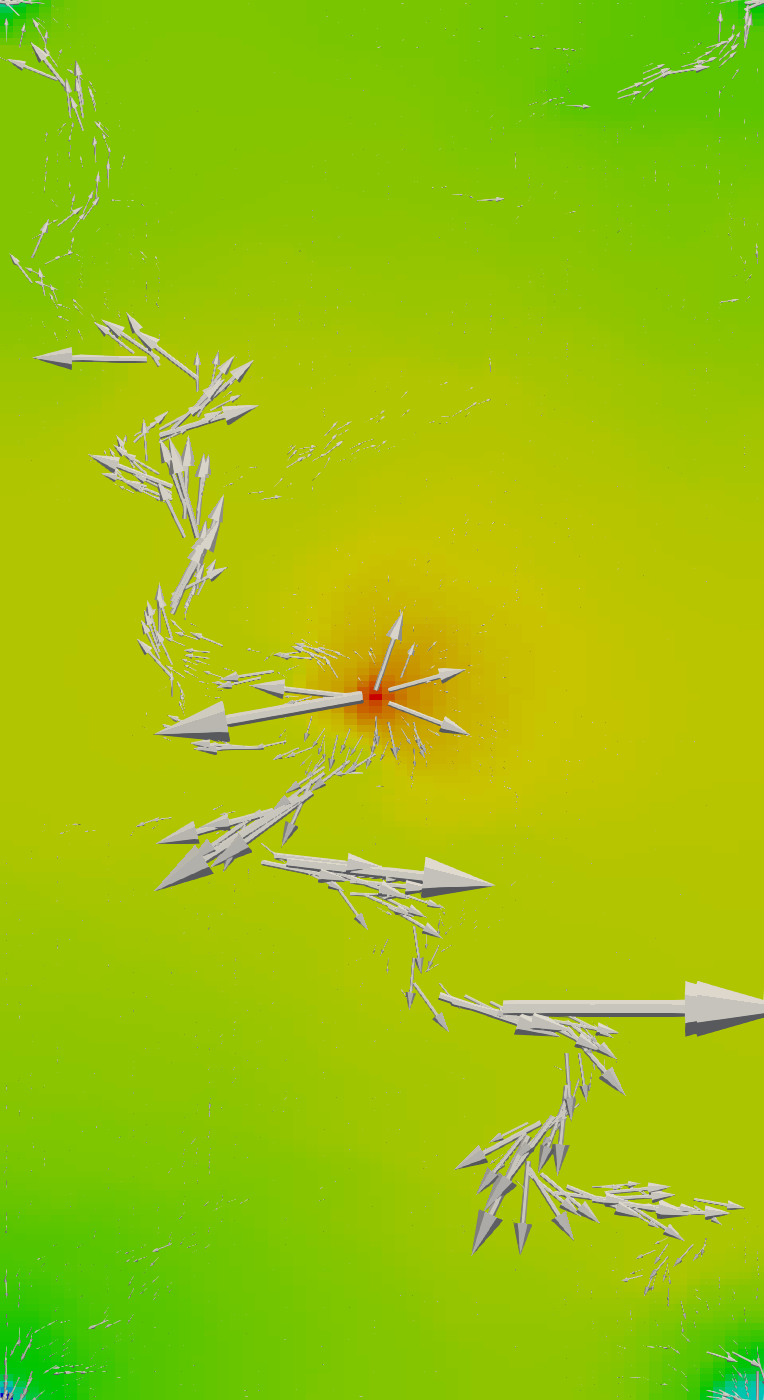}%
    \hspace*{0.0125\textwidth}%
    \includegraphics[scale=0.1]{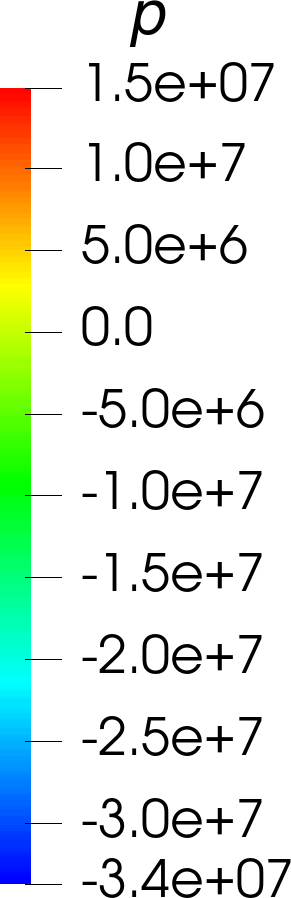}%
    \hspace*{0.02\textwidth}%
    \includegraphics[scale=0.075]{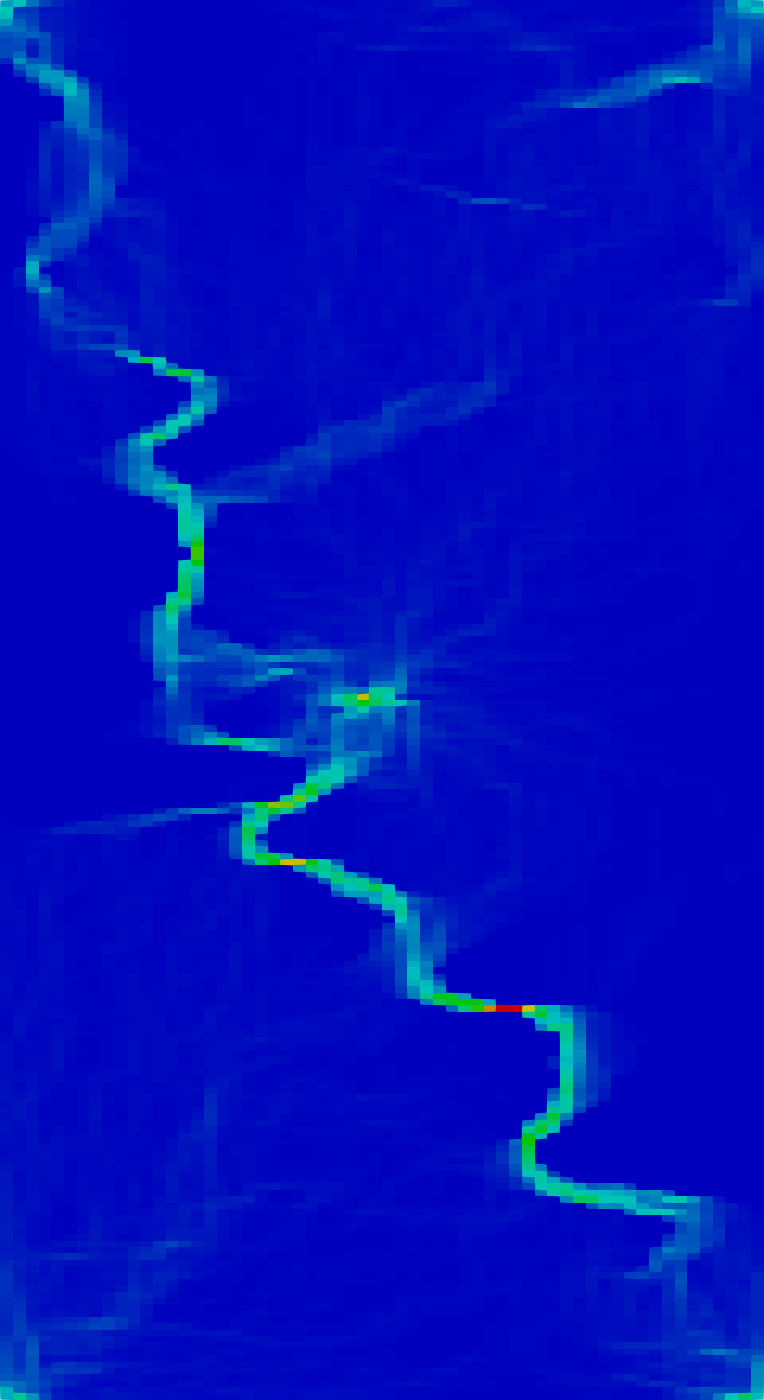}%
    \hspace*{0.0125\textwidth}%
    \includegraphics[scale=0.1]{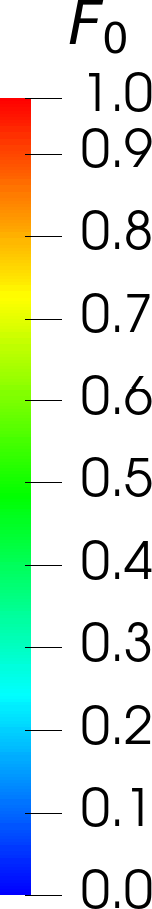}%
    \hspace*{0.02\textwidth}%
    \includegraphics[scale=0.075]{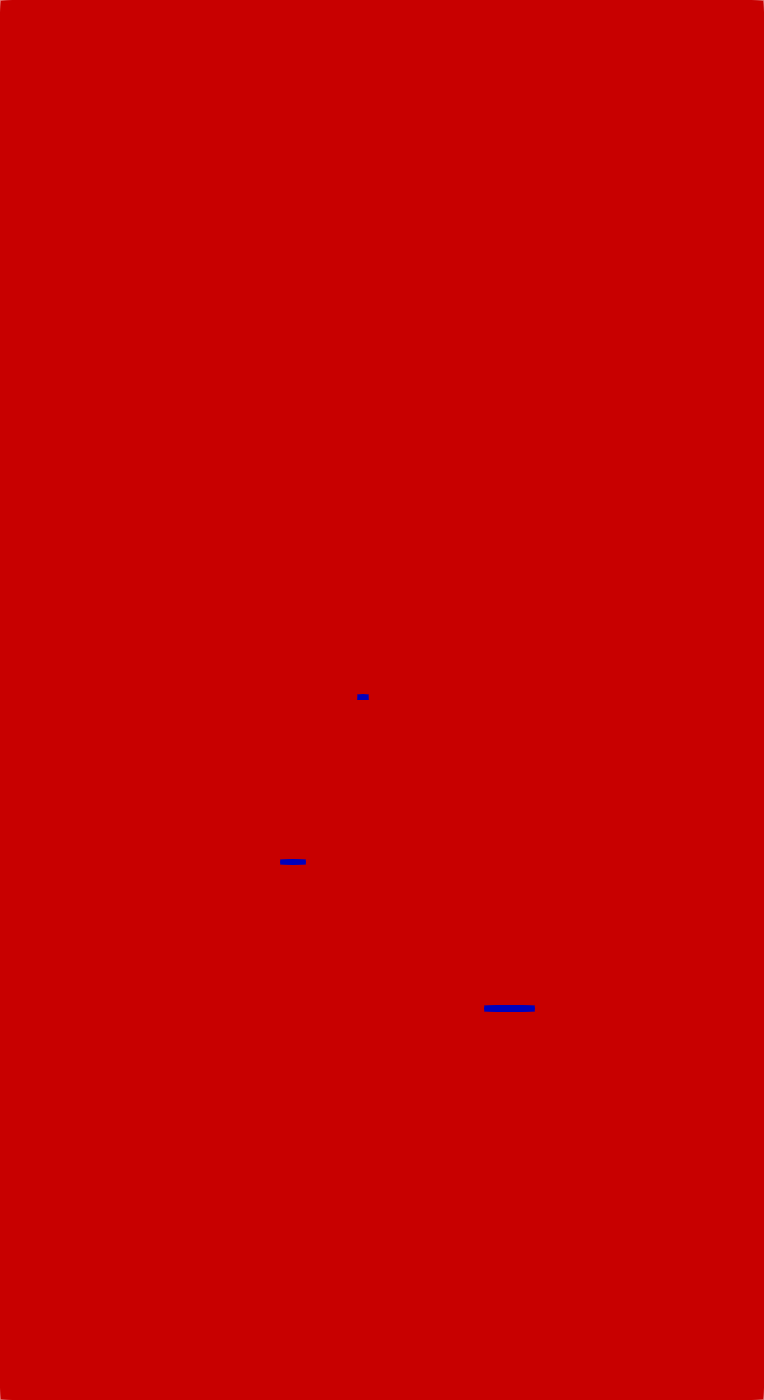}%
    \hspace*{0.0125\textwidth}%
    \includegraphics[scale=0.1]{figs/spe10_region_10_legend.png}\\[0.2cm]
    \includegraphics[scale=0.075]{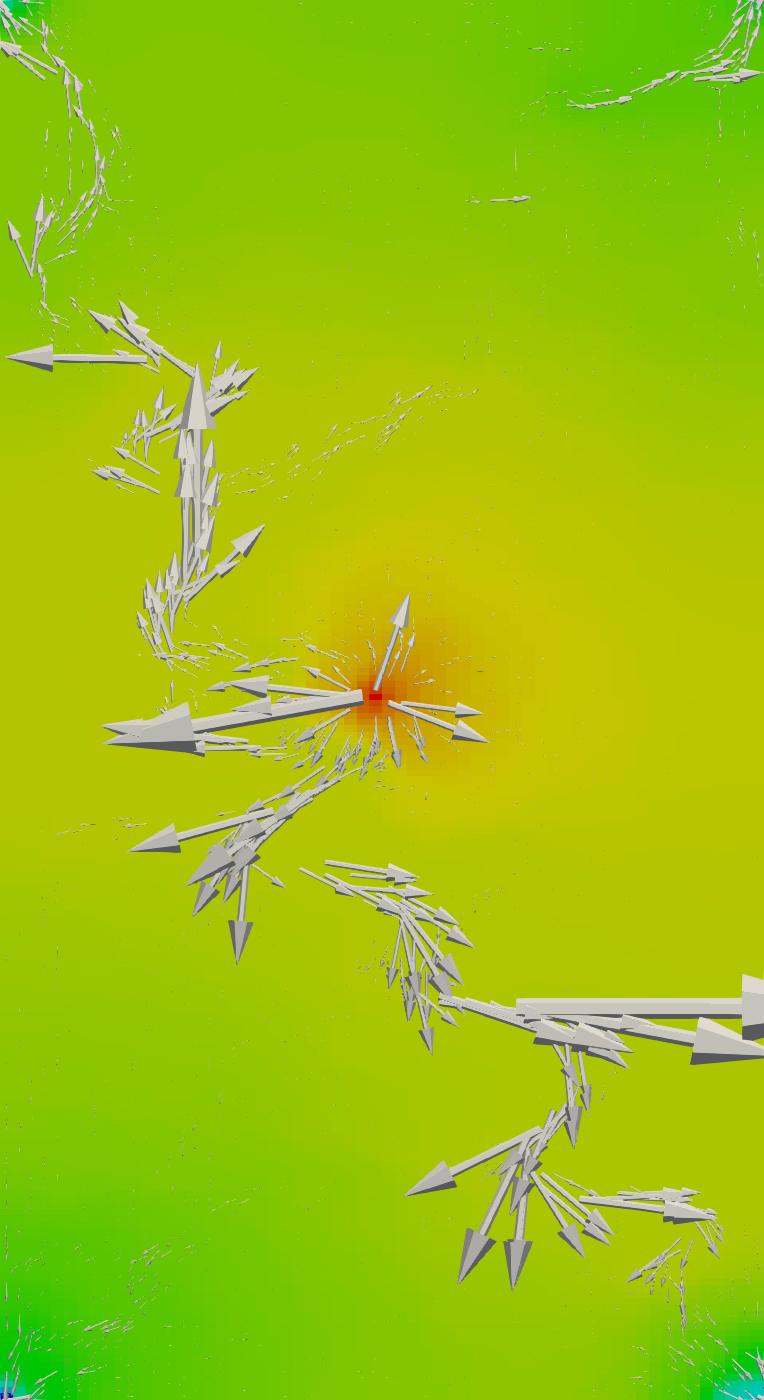}%
    \hspace*{0.0125\textwidth}%
    \includegraphics[scale=0.1]{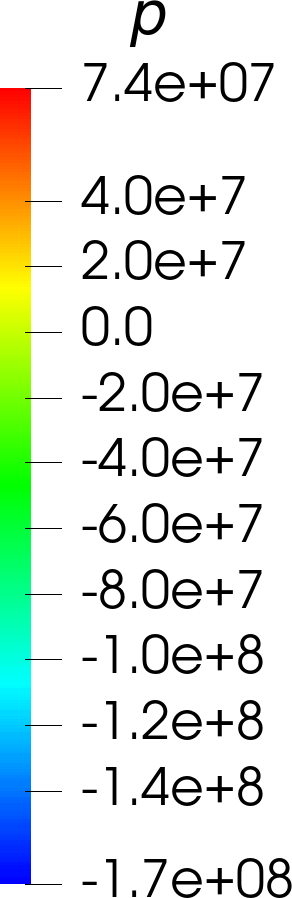}%
    \hspace*{0.02\textwidth}%
    \includegraphics[scale=0.075]{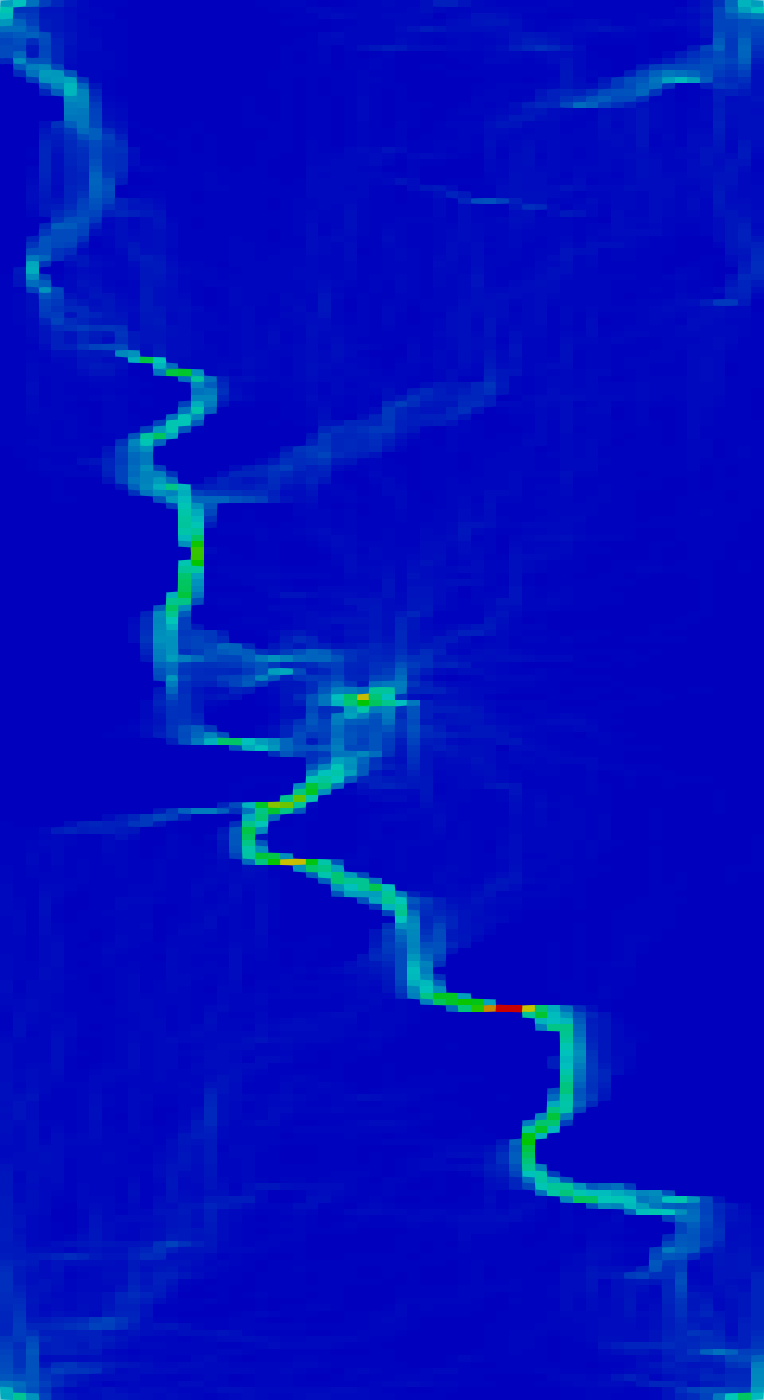}%
    \hspace*{0.0125\textwidth}%
    \includegraphics[scale=0.1]{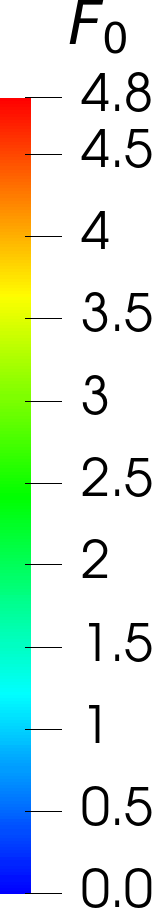}%
    \hspace*{0.02\textwidth}%
    \includegraphics[scale=0.075]{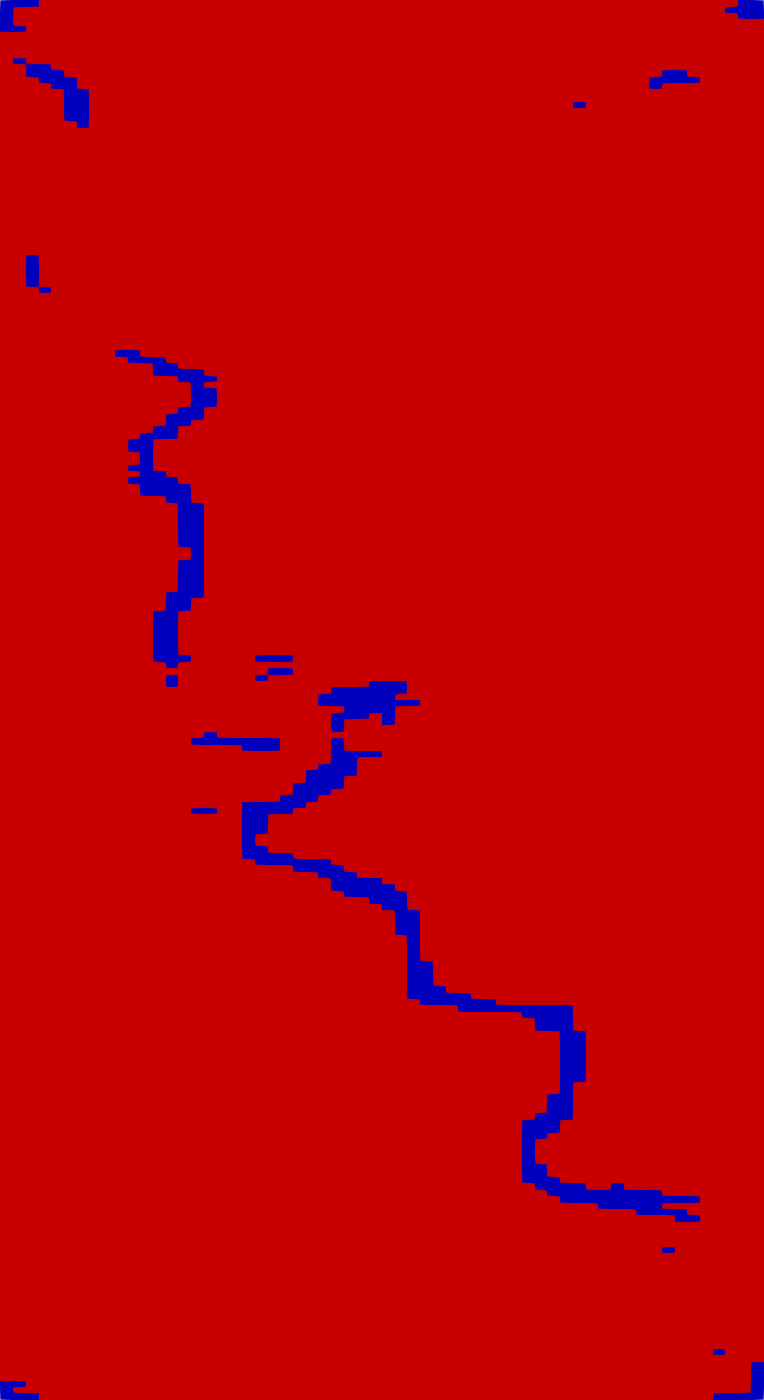}%
    \hspace*{0.0125\textwidth}%
    \includegraphics[scale=0.1]{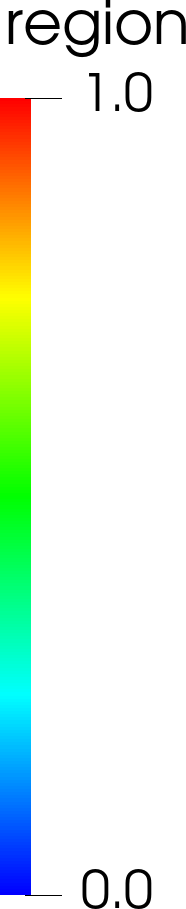}
    \\[0.2cm]
    \includegraphics[scale=0.075]{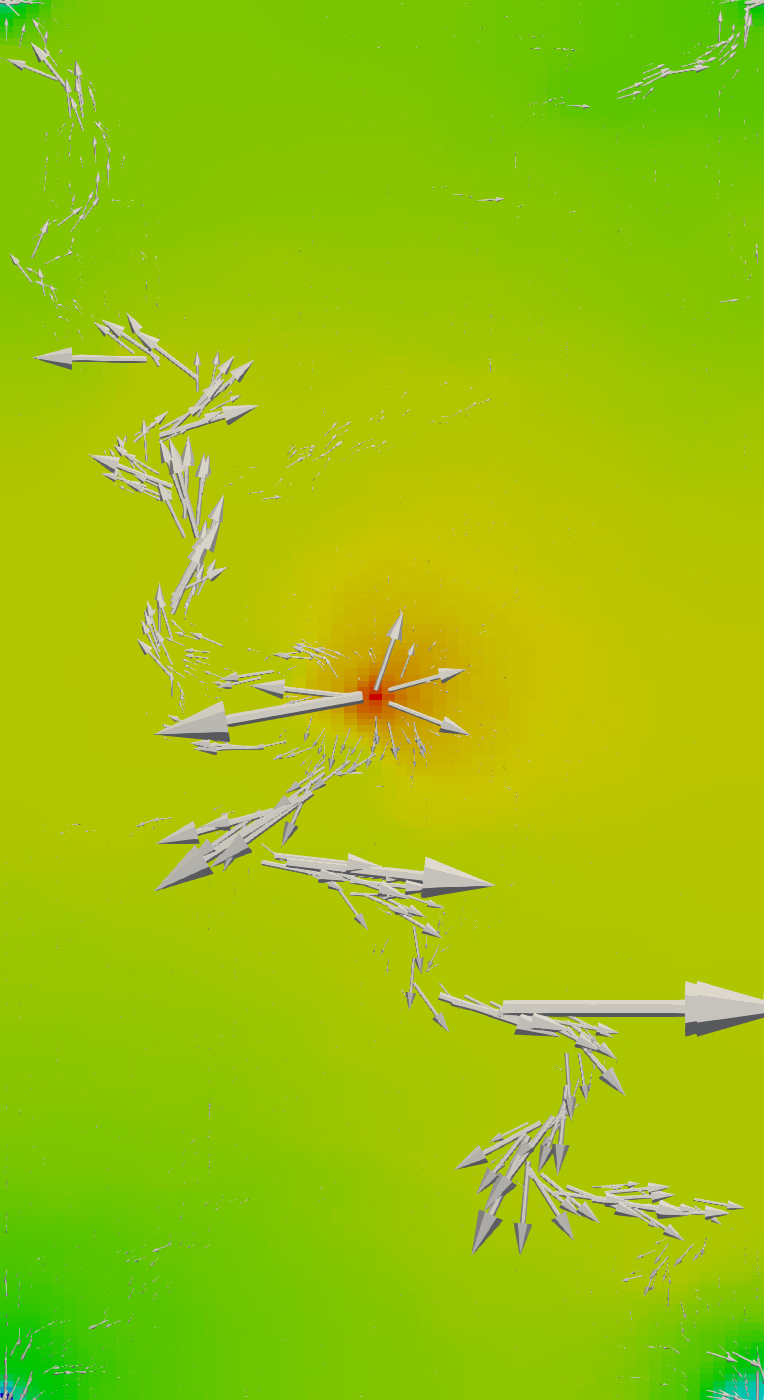}%
    \hspace*{0.0125\textwidth}%
    \includegraphics[scale=0.1]{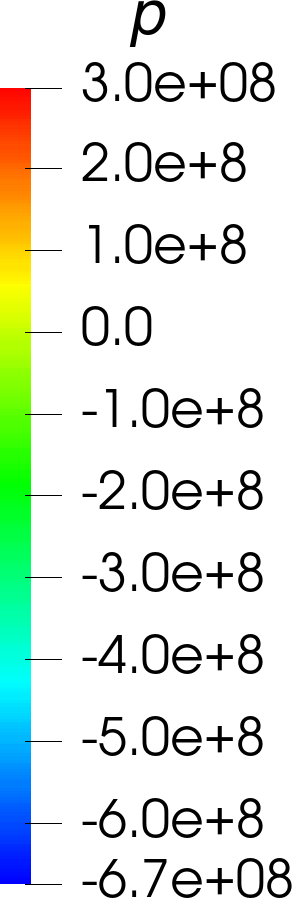}%
    \hspace*{0.02\textwidth}%
    \includegraphics[scale=0.075]{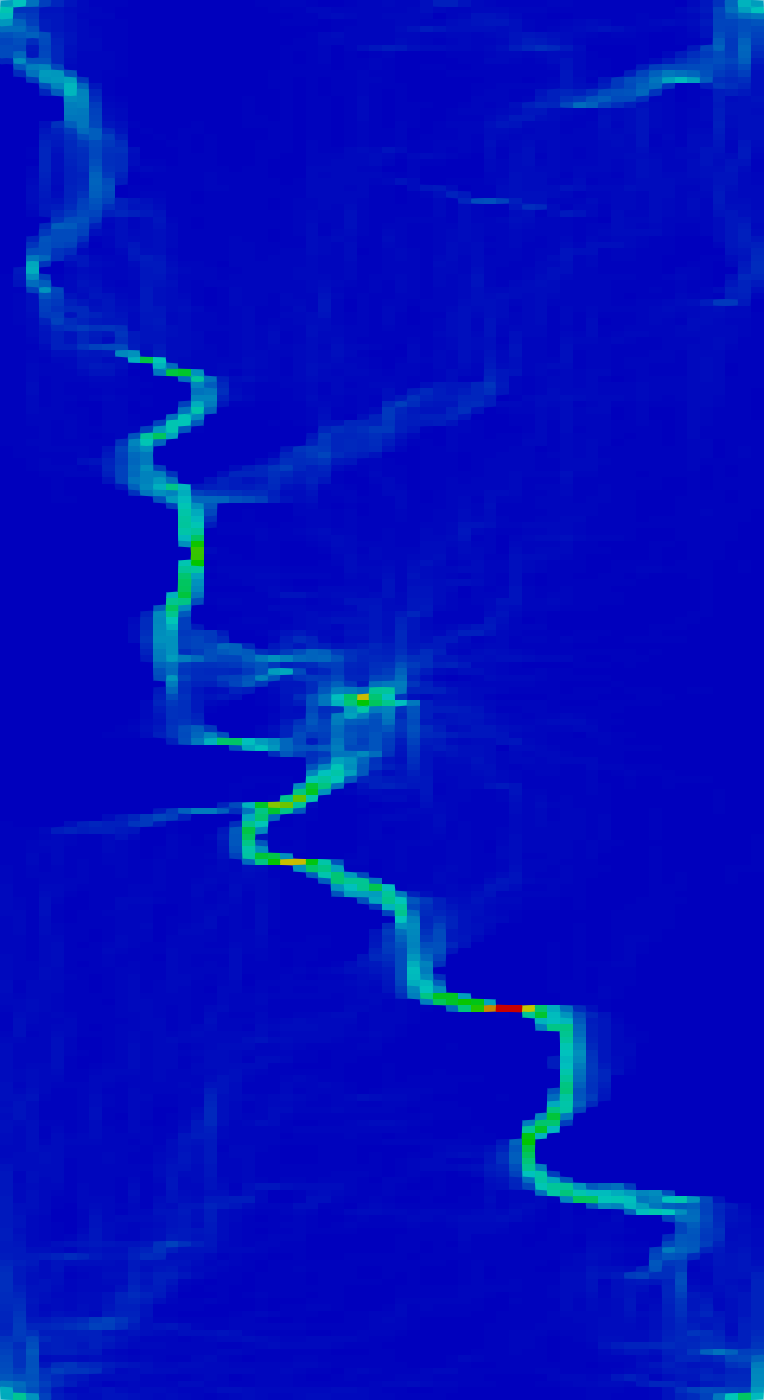}%
    \hspace*{0.0125\textwidth}%
    \includegraphics[scale=0.1]{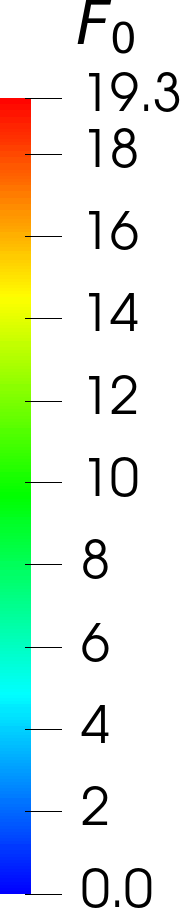}%
    \hspace*{0.02\textwidth}%
    \includegraphics[scale=0.075]{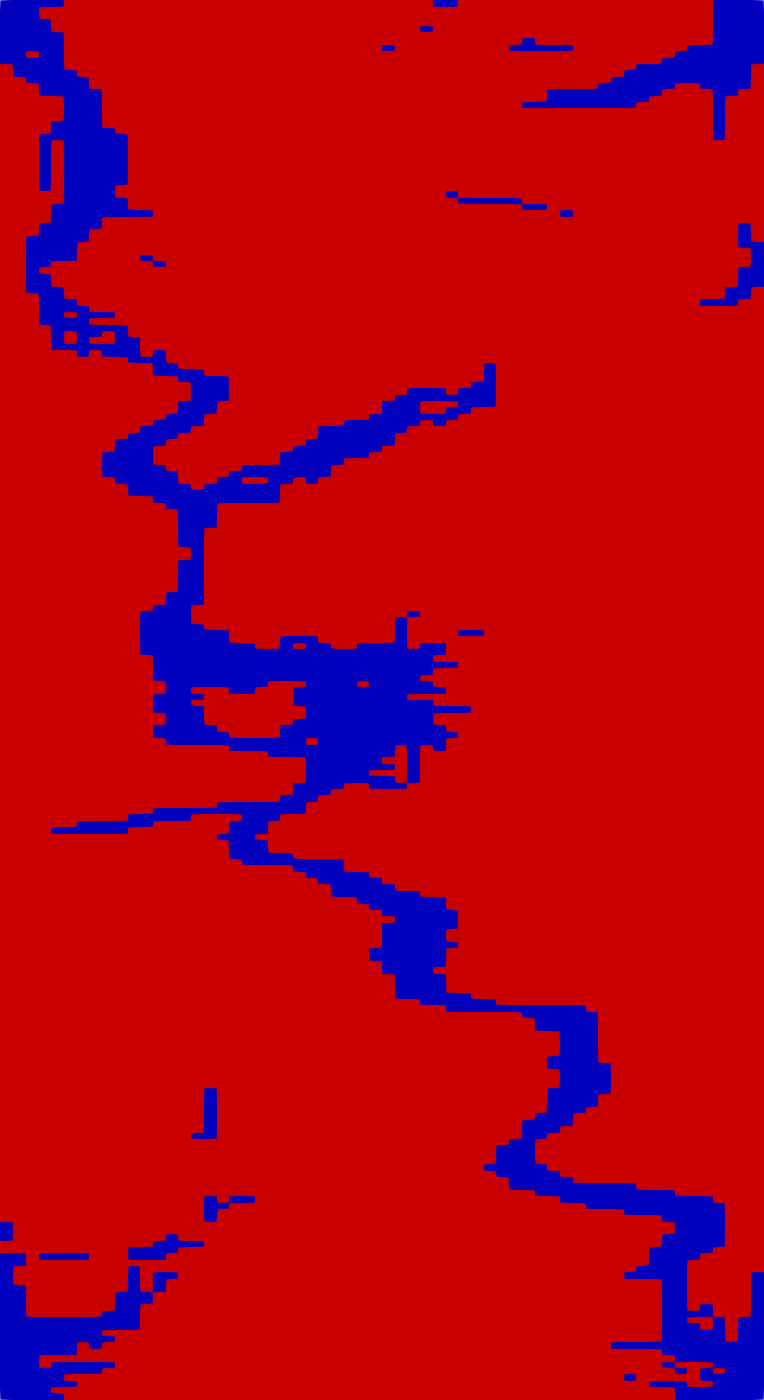}%
    \hspace*{0.0125\textwidth}%
    \includegraphics[scale=0.1]{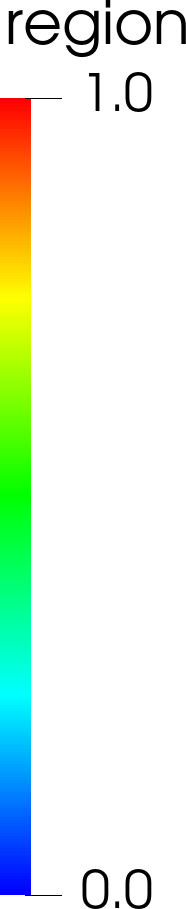}
    \caption{Solutions obtained to the adaptive model for the example in Section \ref{subsec:case2}, in the top row for \textit{Scenario a}, in the middle row for \textit{Scenario b} and in the bottom row for \textit{Scenario c}. The velocity arrows are scaled by $1/92$, $1/23$ and $1/4$, respectively for each scenario.}
    \label{fig:case1a_soo}
\end{figure}

We represent in Figure \ref{fig:case1a_soo} the solution to the adaptive model obtained for the different scenarios. We notice that by increasing the inflow at the injection well, the cells marked as Darcy--Forchheimer become more numerous and, in particular, are preferably selected along the high-permeability channels connecting some of the wells. Even if the two wells on the bottom left and top right are not in highly permeable areas, for \textit{Scenario c}, some Darcy--Forchheimer cells start to be selected around these two wells and, in particular for the second one, in the highly permeable region separating the well itself from the central injector.

\begin{figure}[!ht]
    \centering
    \includegraphics[scale=0.0913]{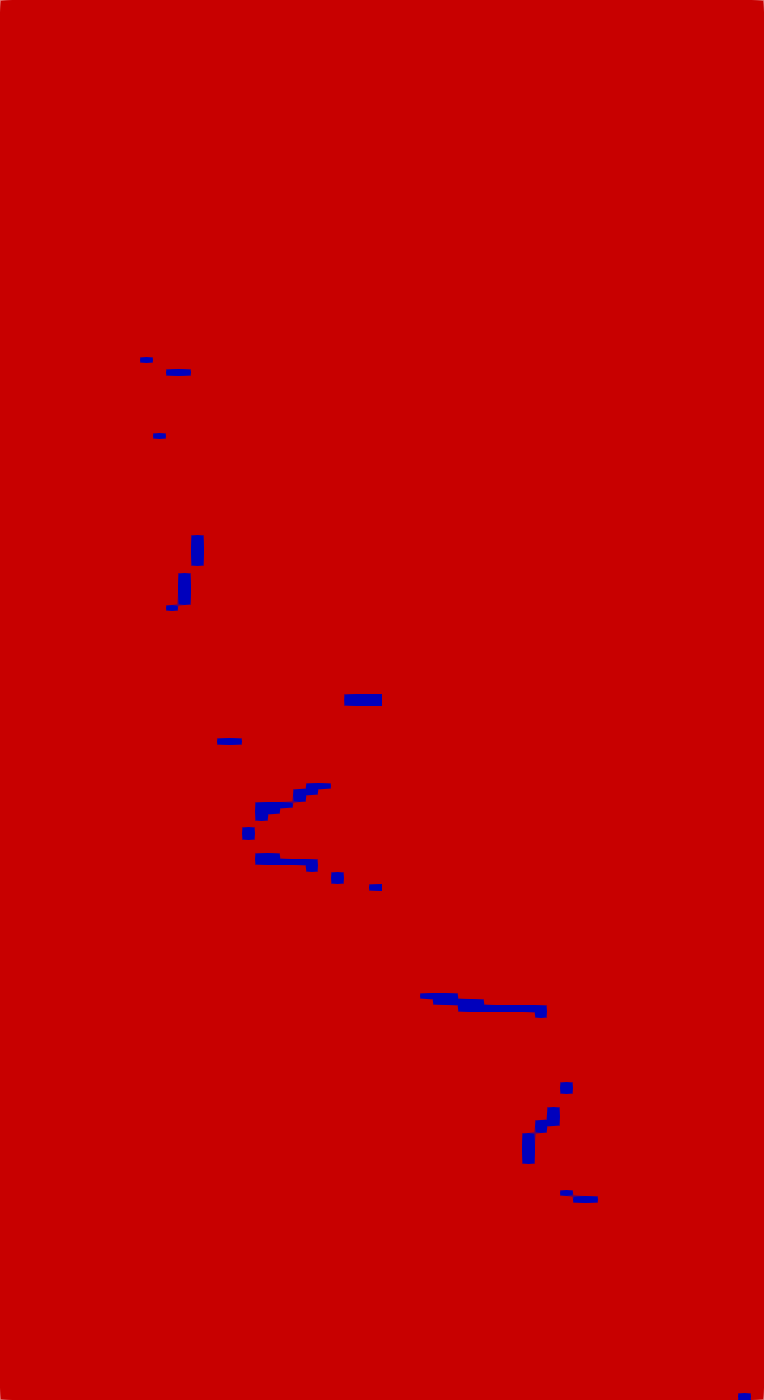}%
    \hspace*{0.02\textwidth}%
    \includegraphics[scale=0.0913]{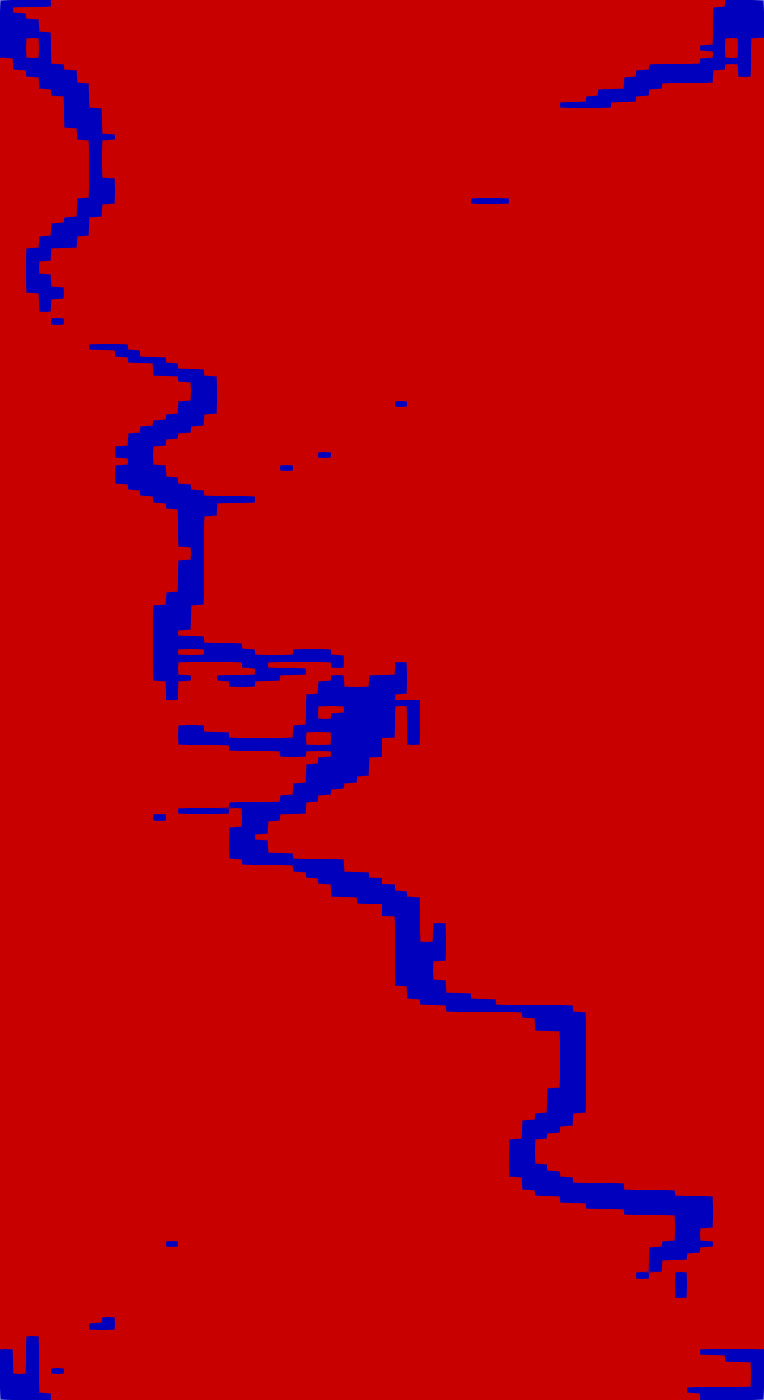}%
    \hspace*{0.02\textwidth}%
    \includegraphics[scale=0.0913]{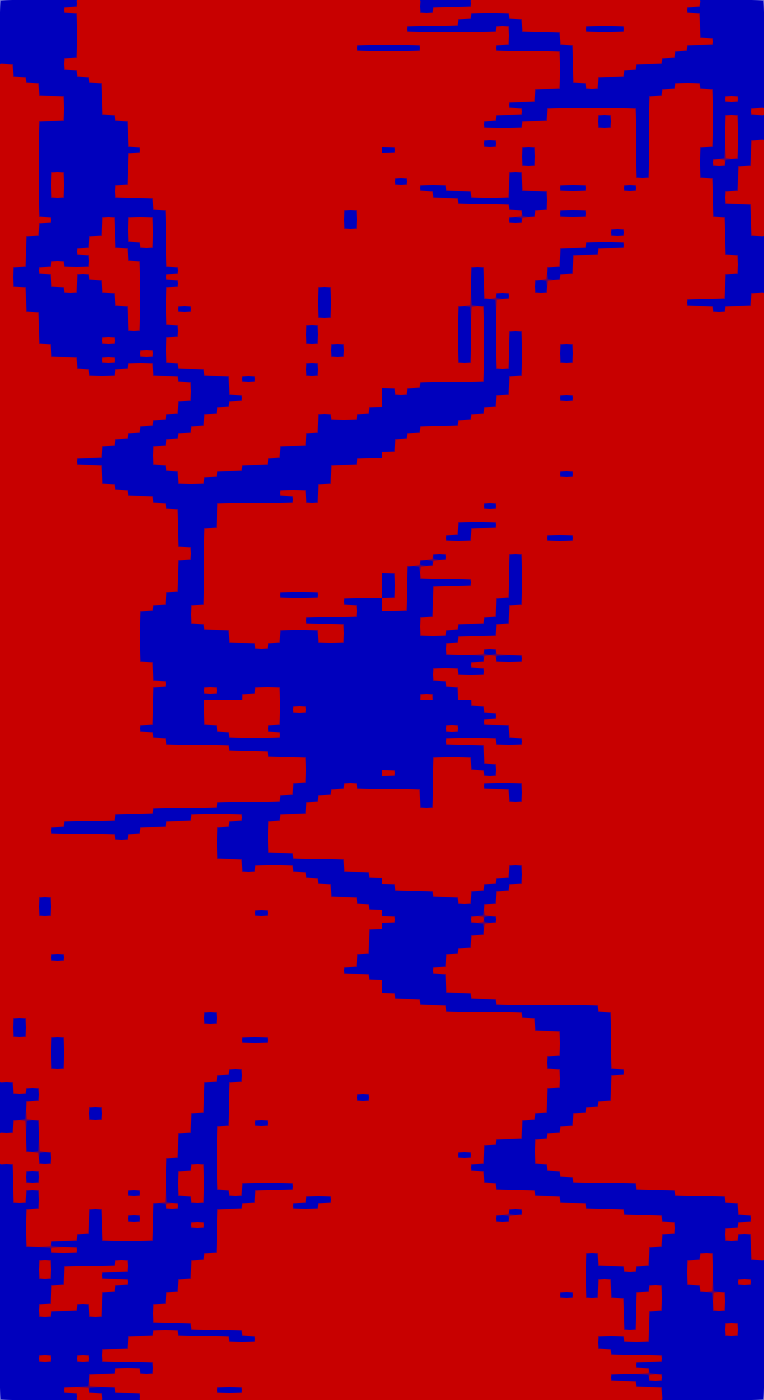}%
    \hspace*{0.0125\textwidth}%
    \includegraphics[scale=0.1]{figs/spe10_region_200_legend.png}\\[0.2cm]
    \includegraphics[scale=0.0913]{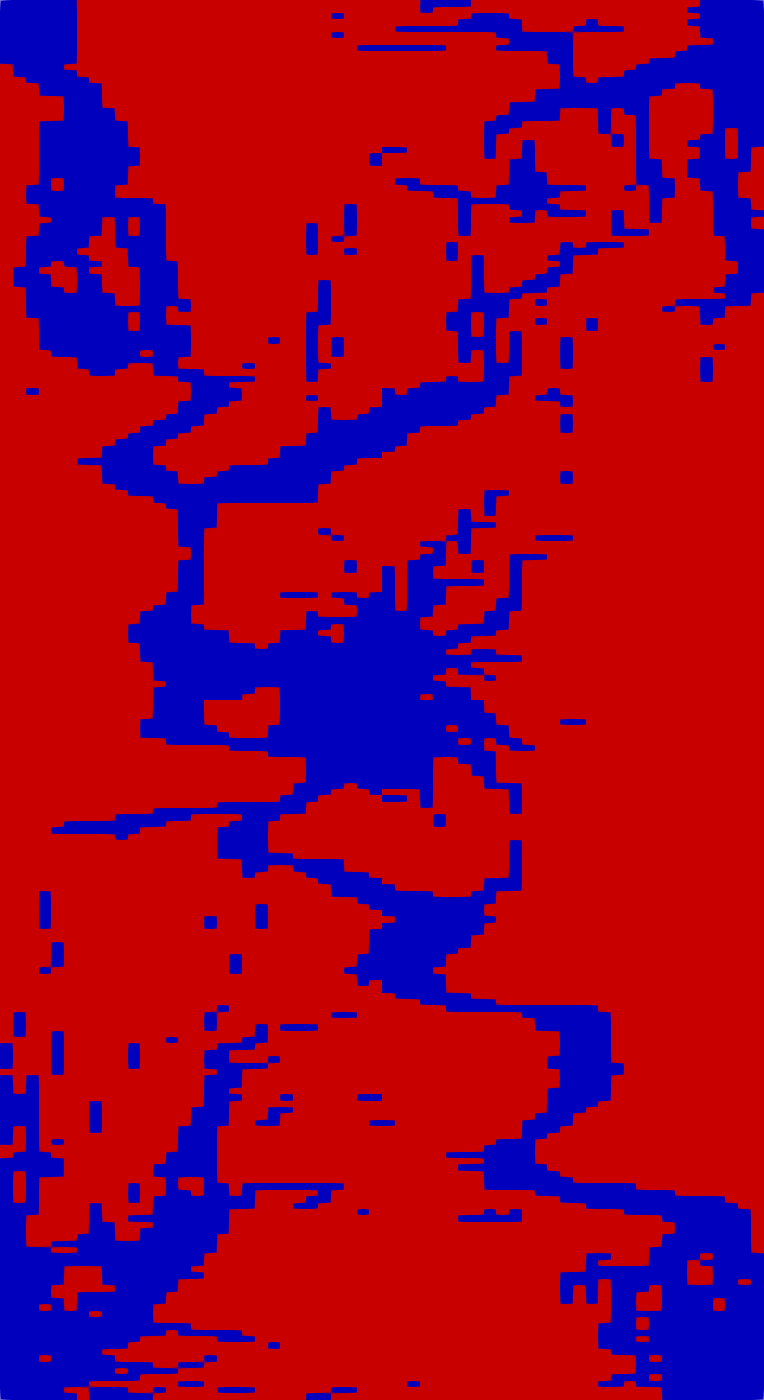}%
    \hspace*{0.02\textwidth}%
    \includegraphics[scale=0.0913]{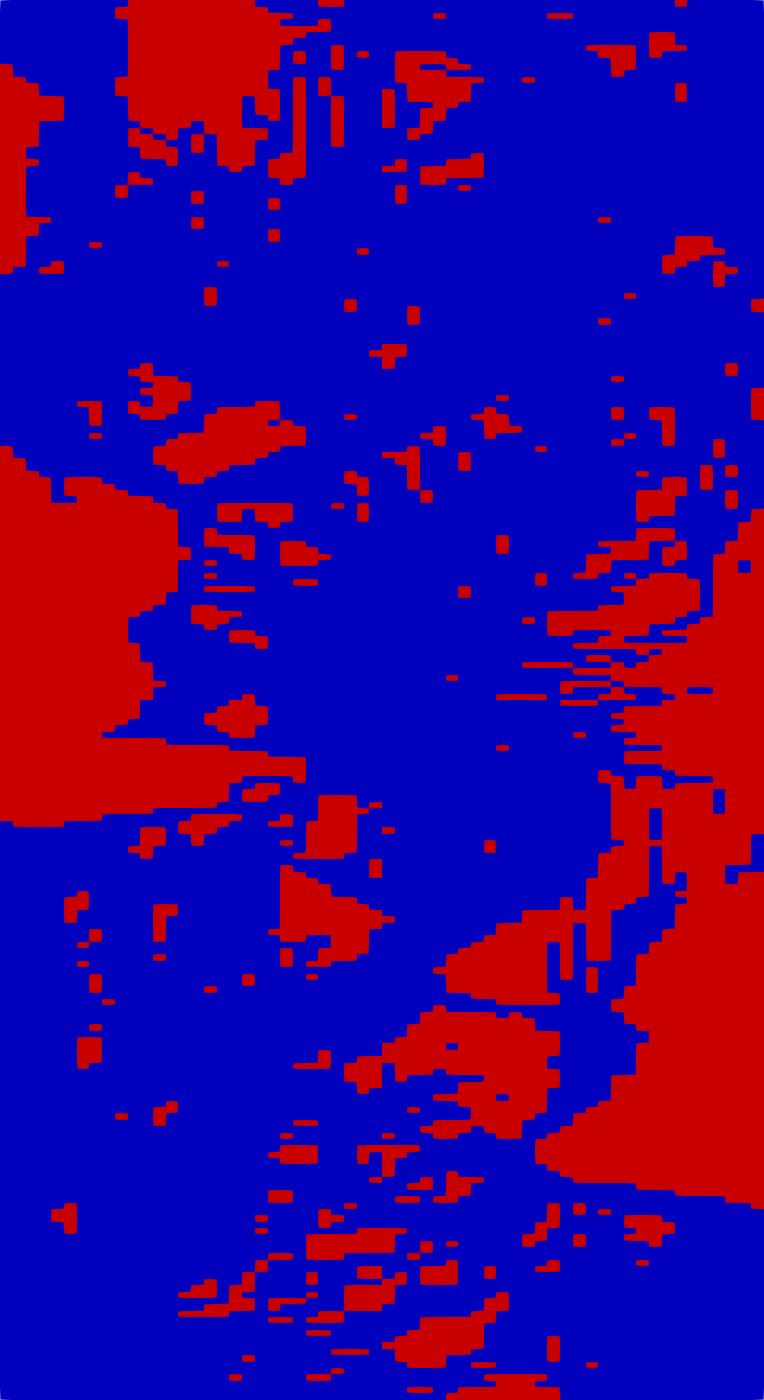}%
    \hspace*{0.02\textwidth}%
    \includegraphics[scale=0.0913]{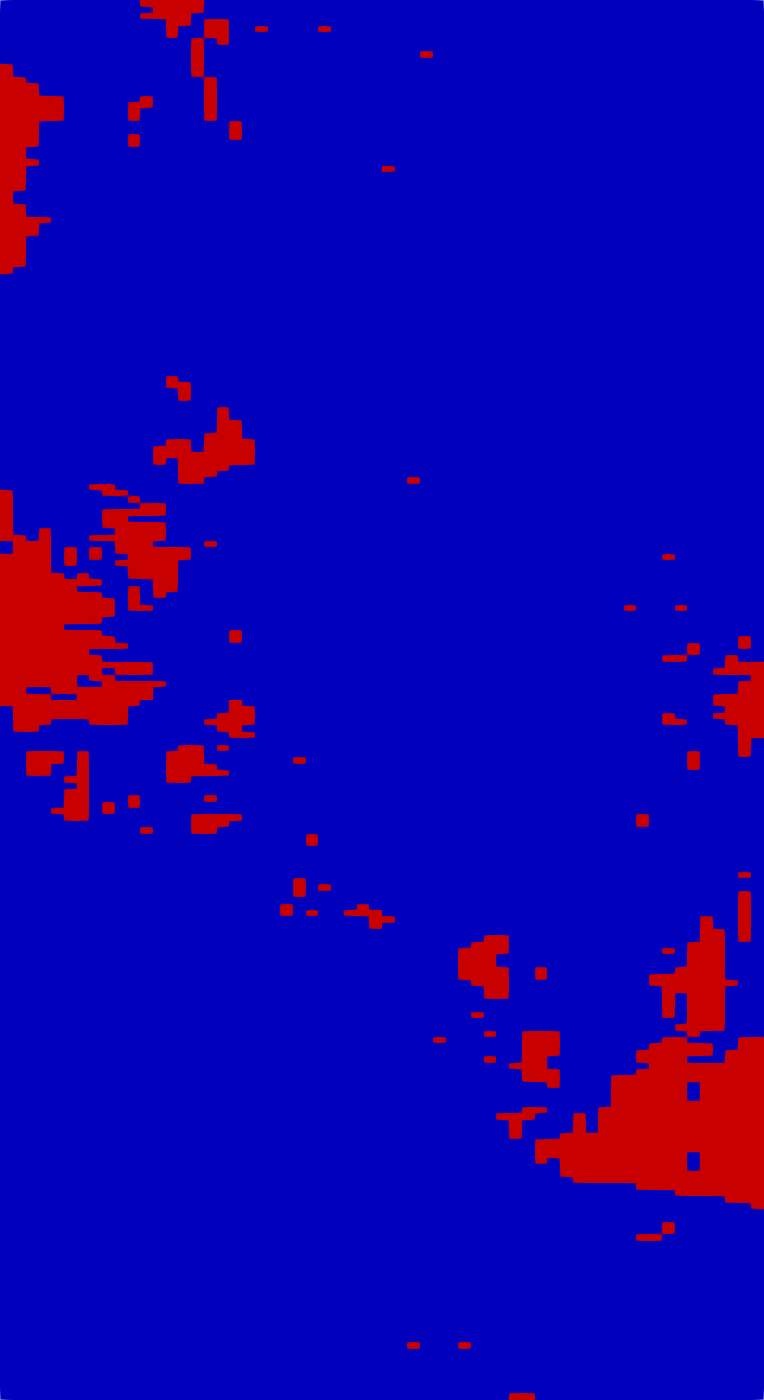}%
    \hspace*{0.0125\textwidth}%
    \includegraphics[scale=0.1]{figs/spe10_region_200_legend.png}
    \caption{Darcy (red) and Darcy--Forchheimer (blue) regions depending on the error tolerance $\delta$ decreasing from left to right according to $\delta=0.05$, $\delta=0.0125$ and $\delta=0.003125$, in the top row for \textit{Scenario a} and in the bottom row for \textit{Scenario c}, for the example in Section \ref{subsec:case2}.}
    \label{fig:case2a_smaller_error}
\end{figure}

Figure \ref{fig:case2a_smaller_error} shows the Darcy--Forchheimer cells with decreasing error tolerance $\delta$. We notice that the number of cells belonging to the Darcy--Forchheimer region increases by first following highly permeable paths from the injection well to the production wells, and then accepting also regions of lower permeability. In the last figure, the error tolerance is so small that most of the cells are categorized as Darcy--Forchheimer. We expect that by reducing $\delta$ even further, we would get all the cells into the nonlinear subdomain. 

\begin{figure}[!ht]
    \centering
    \includegraphics[width=0.75\textwidth]{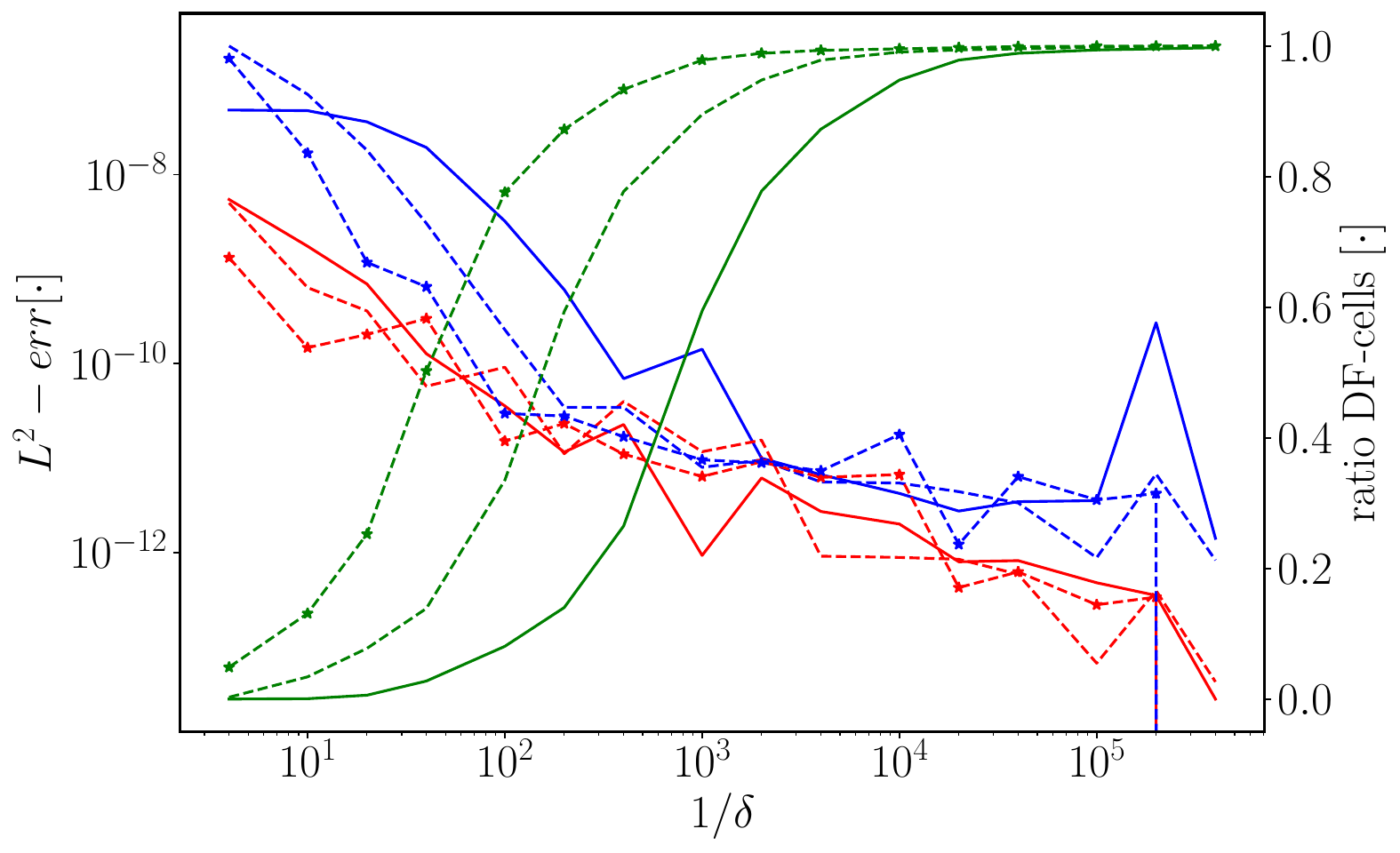}\\
    \includegraphics[width=0.75\textwidth]{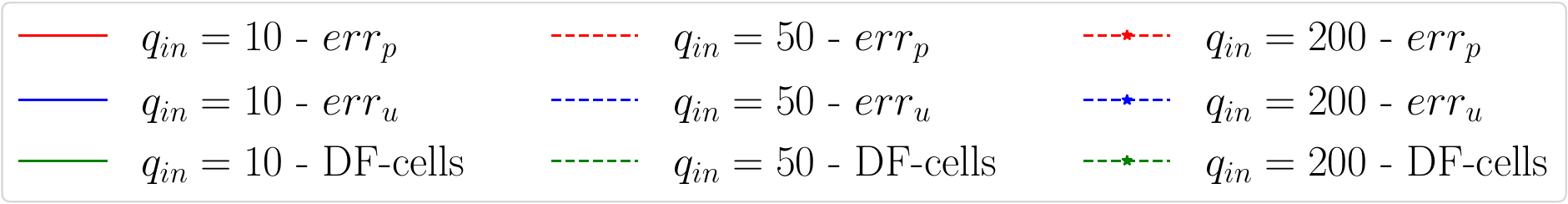}%
    \caption{Pressure (in red) and flux (in blue) $L^2$ relative errors by taking smaller values of the error tolerance $\delta$. The solid lines refer to \textit{Scenario a}, the dashed lines to \textit{Scenario b}, and the dashed-and-dotted lines to \textit{Scenario c}, for the example of Section \ref{subsec:case2}. In green, the ratio between the number of Darcy--Forchheimer cells and the total number of cells.}
    \label{fig:case1_errorii}
\end{figure}

In Figure \ref{fig:case1_errorii}, we show the $L^2$ errors in pressure and flux for the adaptive model relative to the global Darcy--Forchheimer model, as well as the number of Darcy--Forchheimer cells selected by the adaptive model, all as functions of the error tolerance $\delta$. As for the examples of Section \ref{subsec:case1}, by lowering the error tolerance, the adaptive approach tends to have a solution that is more similar to the global Darcy--Forchheimer model,  and thus yields a smaller error. In fact, when the percentage of the Darcy--Forchheimer cells reaches $1$, the adaptive model is equivalent to the global Darcy--Forchheimer one and thus has an error of exactly zero.

\subsection{Three-dimensional test case}\label{subsec:case3}
%(case4 git)

In this last test case, we consider a tridimensional portion of the SPE10 benchmark \cite{Christie2001}, composed of the first 15 layers. However, to keep the computational cost limited, we select only the first 10 cells in the $x$ direction and the first 20 in the $y$ direction. As done in the previous example, we consider fives wells: one in the centre of the domain, the injector, and four in the corners, the producers. Here, the wells are represented as monodimensional objects and a Peaceman model \cite{Peaceman1978,Peaceman1983} is used for the coupling with the surrounding porous medium. The wells are placed in the middle of each of their cells and are nonmatching with the porous-medium grid. By keeping the injection flow rate fixed, equal to $10 \sib{\kilogram\per\second}$, we change the depth of the wells and see the impact on the selection procedure of Darcy and Darcy--Forchheimer cells by the adaptive model.

\begin{figure}[!ht]
    \centering
    \includegraphics[width=0.33\textwidth]{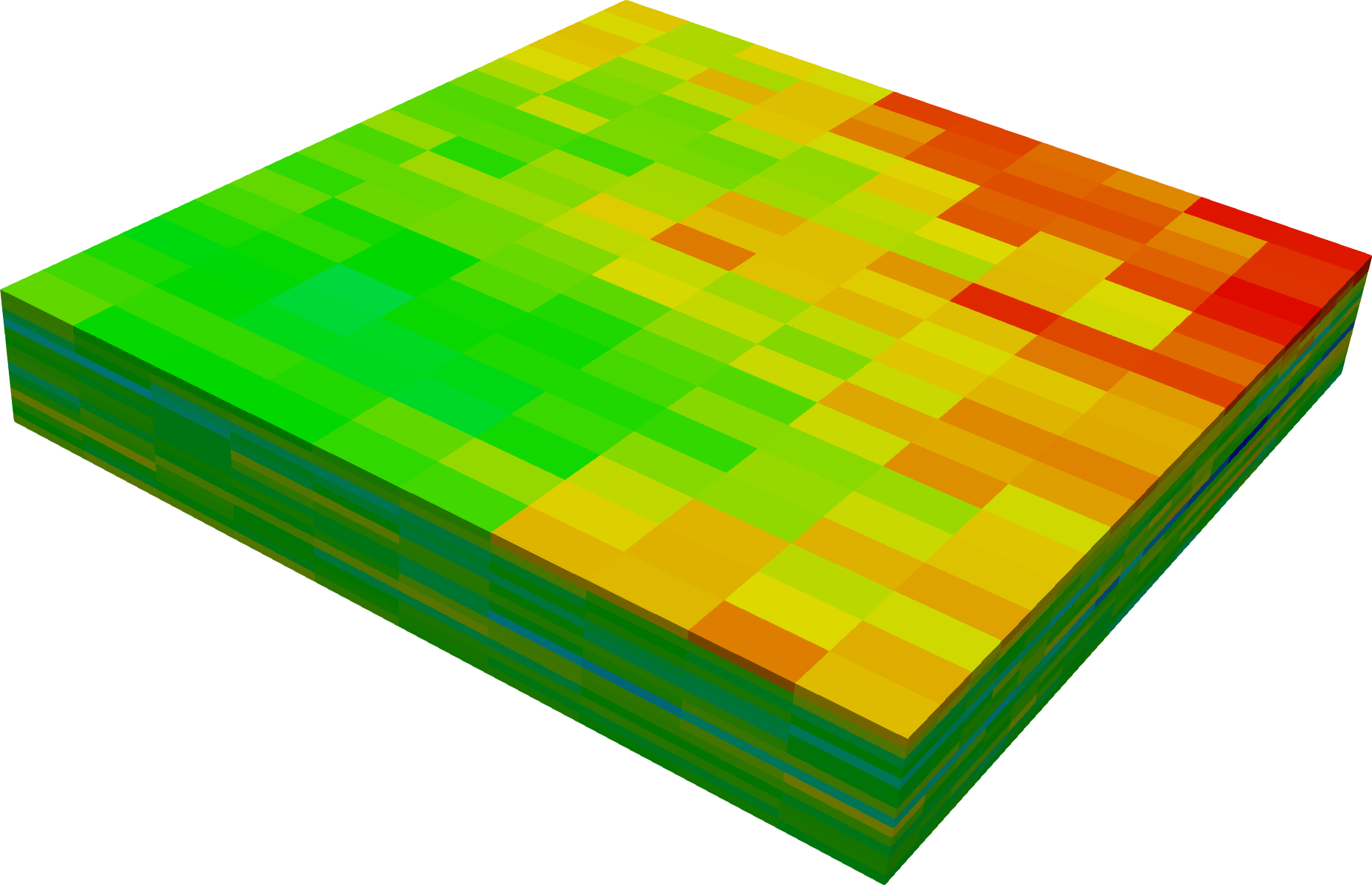}%
    \hspace{0.025\textwidth}%
    \includegraphics[width=0.05\textwidth]{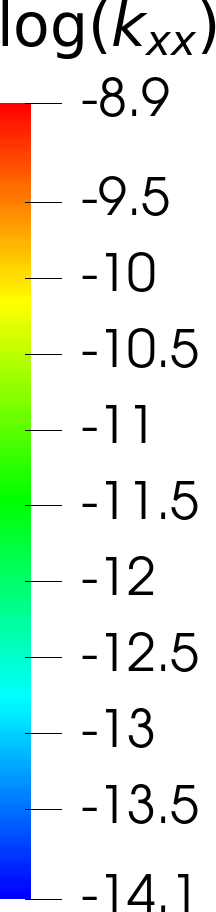}%
    \hspace{0.1\textwidth}%
    \includegraphics[width=0.33\textwidth]{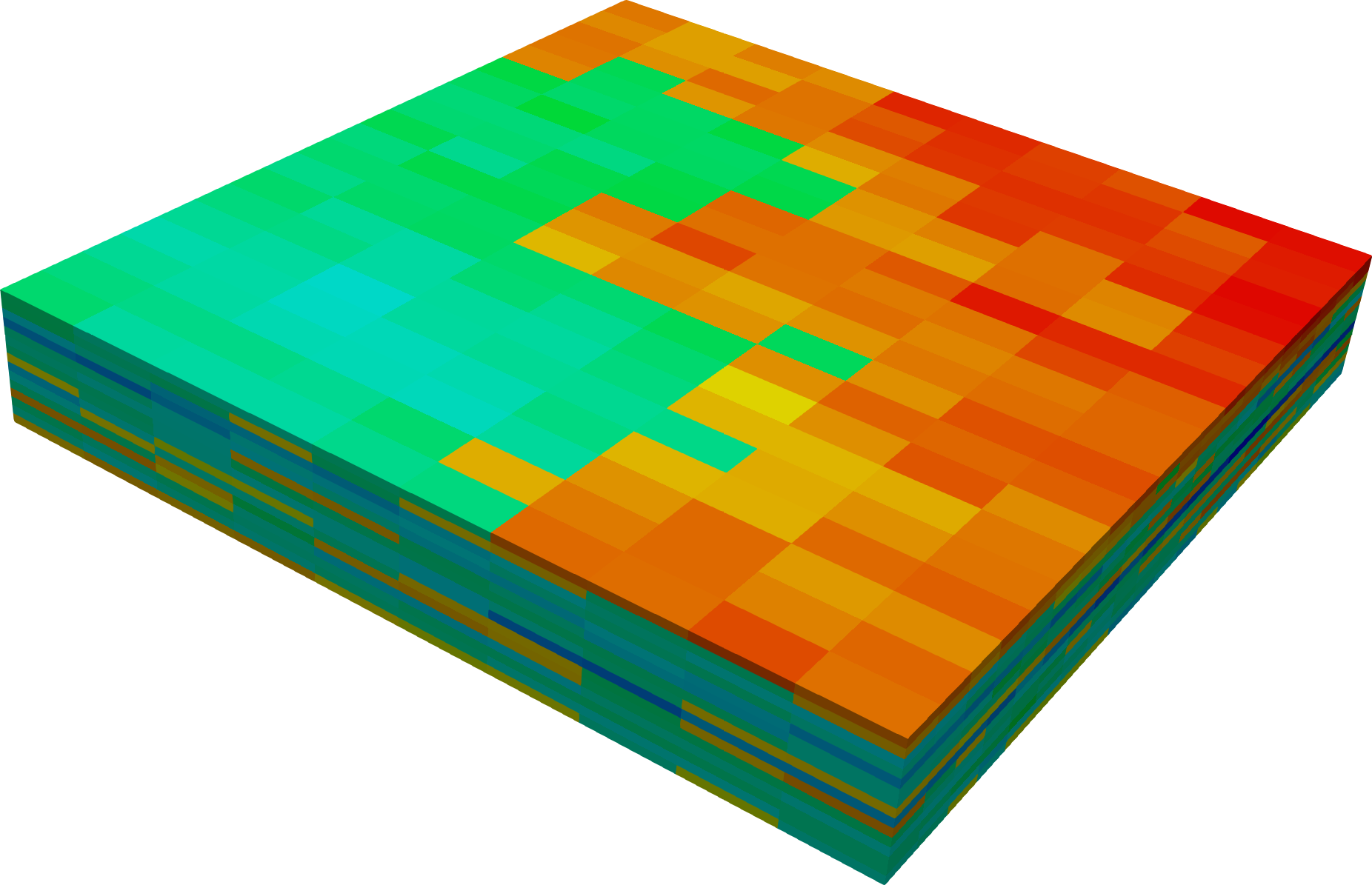}%
    \hspace{0.025\textwidth}%
    \includegraphics[width=0.05\textwidth]{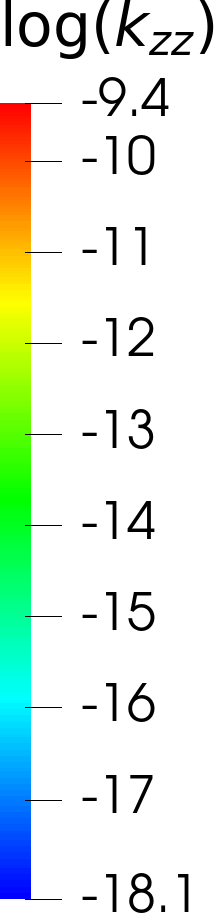}%
    \caption{Log-permeability field for the example in Section \ref{subsec:case3} in the horizontal (left) and vertical (right) directions.}
    \label{fig:case4_perm}
\end{figure}

Although the models presented in Section \ref{sec:model} only consider a scalar permeability $\k$, we wish here to take into account the anisotropy given by the SPE10 case study in the vertical direction. To do so, we alter the Darcy--Forchheimer law slightly:
\begin{equation*}
  -\bgrad p + \bff = \left( 1 + c_{\mt{F}}\, \mt{Re}[\bu] \right) \nu\K^{-1}\, \bu \quad \text{in $\Omega$},
\end{equation*}
where $\K$ is now the permeability tensor and the Reynolds number is updated into
\begin{equation*}
  \mt{Re}[\bu] = \frac{\det(\K)^{1/6}}{\mu} \norm{\bu}.
\end{equation*}
We refer the reader to~\cite{SV01,KL95} for a justification of such a law. Note, however, that this form does not fit in the context of~\cite{FP21,FP23}, where the adaptive model is derived and shown to be well posed, since, for that, the norm of the flux in the Reynolds number would need to be weighed by $\K^{-1}$. We leave the precise numerical treatment of the tensor model for future investigation and accept for now the validity, in our context, of the above three-dimensional version of the Darcy--Forchheimer law. We have in fact, for the SPE10 study, a permeability tensor which is diagonal with equal values in the $x$ and $y$ directions and a different value in the $z$ direction. A representation of the permeability field is reported in Figure \ref{fig:case4_perm}.

\begin{figure}[!ht]
    \centering
    \includegraphics[width=0.33\textwidth]{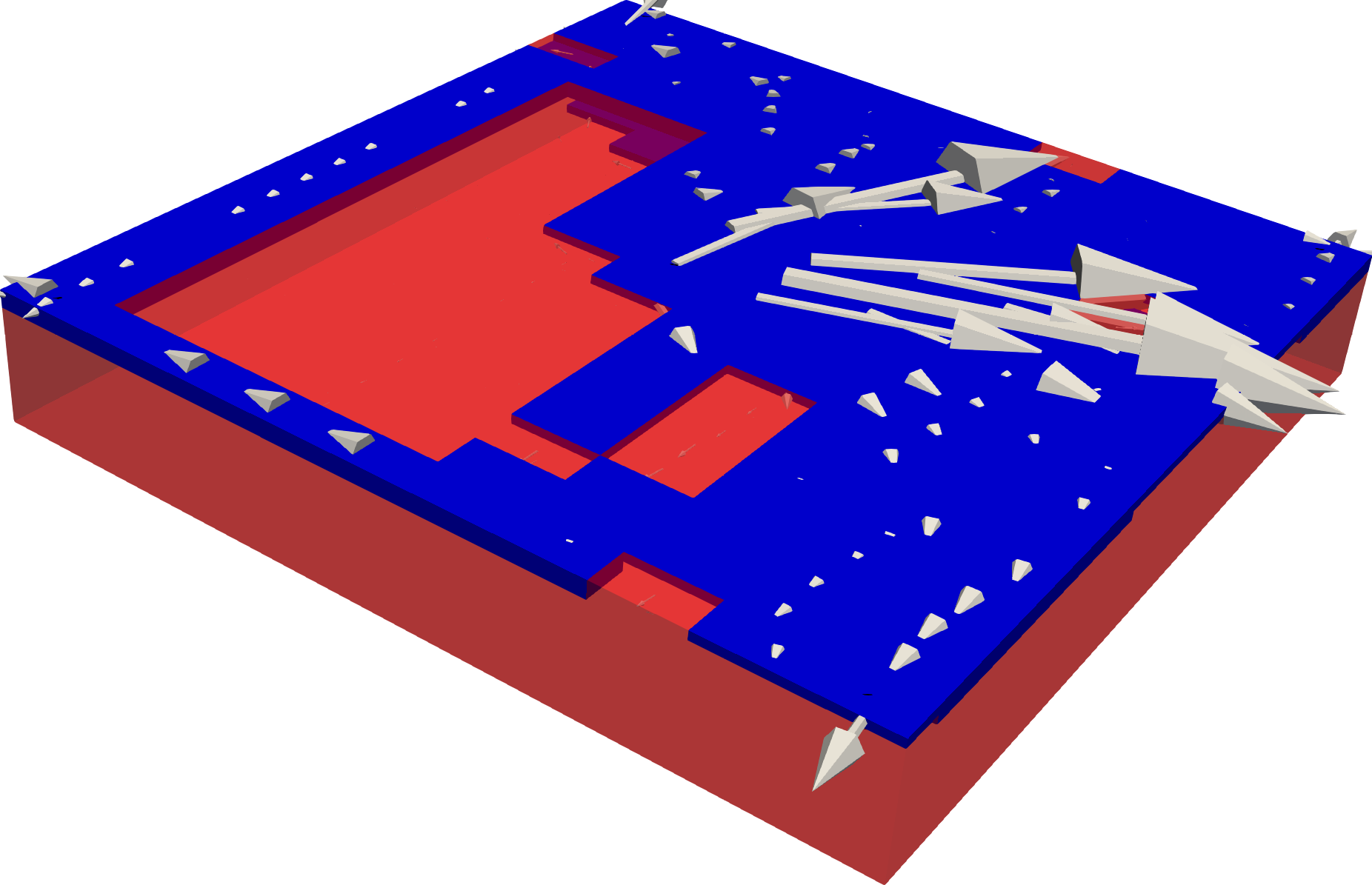}%
    \includegraphics[width=0.33\textwidth]{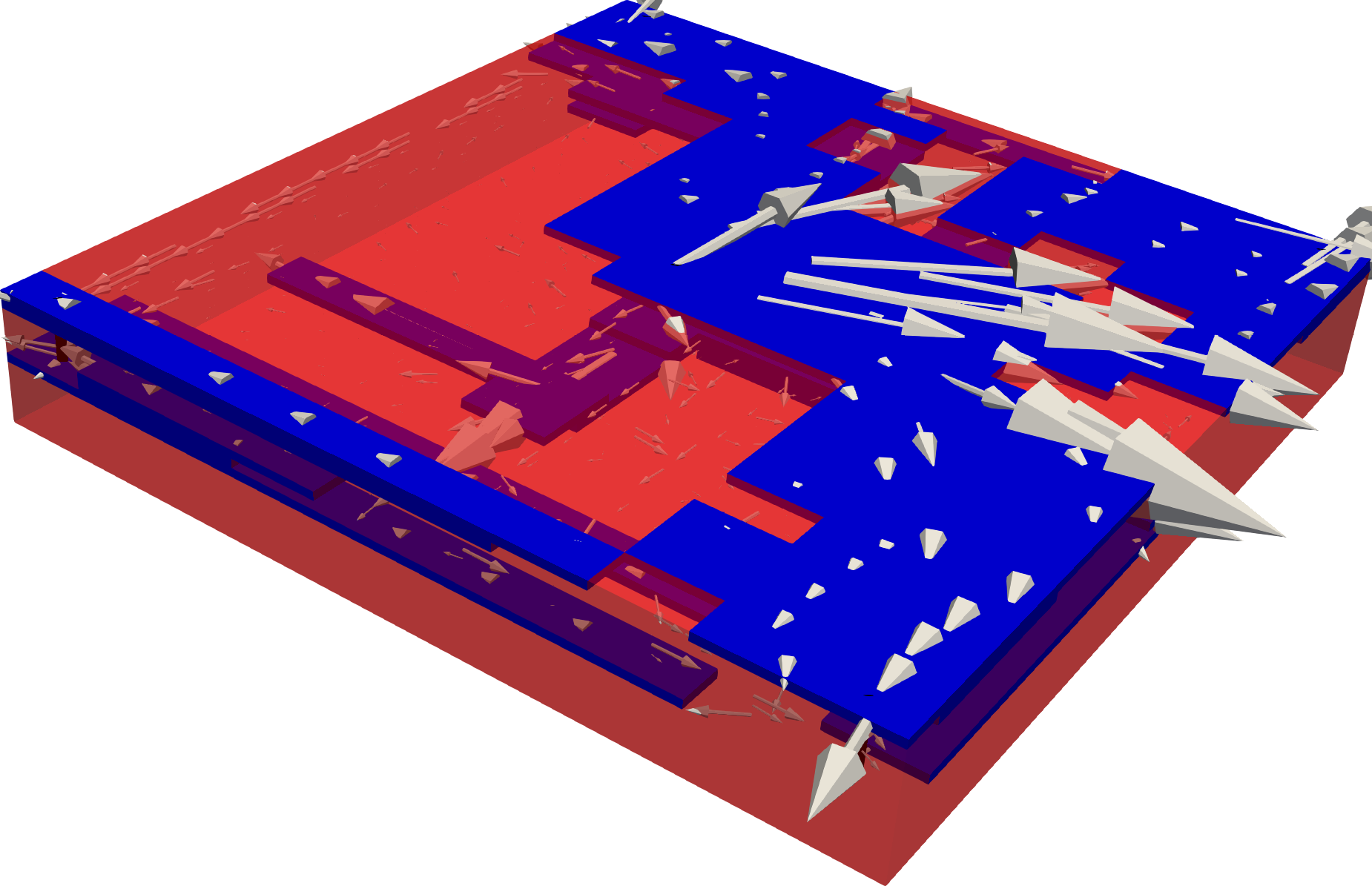}%
    \includegraphics[width=0.33\textwidth]{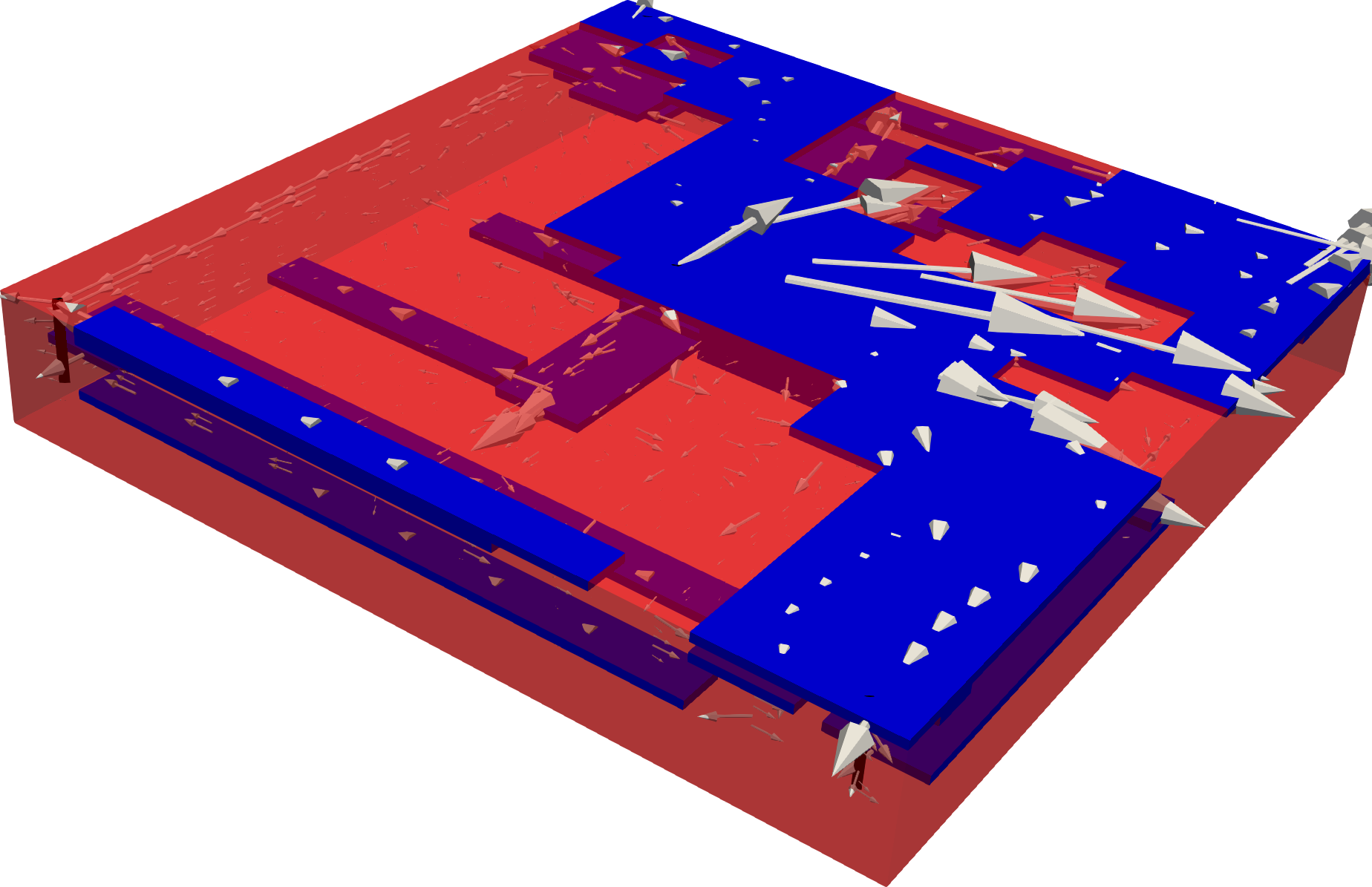}%
    \caption{Darcy (opaqued red) and Darcy--Forchheimer (blue) regions for different well depths for the example in Section \ref{subsec:case3}. From left to right: \textit{Scenario a} (3.5 cells deep), \textit{Scenario b} (7.5) and \textit{Scenario c} (9.5). The velocity arrows for \textit{Scenario a} are halved compared to the other two scenarios.}
    \label{fig:case4_region}
\end{figure}

In Figure \ref{fig:case4_region}, we show a comparison of the obtained adaptive solution by making the wells deeper: 3.5 cells for \textit{Scenario a}, 7.5 cells for \textit{Scenario b}, and 9.5 cells for \textit{Scenario c}. We notice that, by considering the wells deeper, we get more Darcy--Forchheimer cells since indeed more regions of high permeability are thus connected. For the same reason, we observe that the number of Darcy--Forchheimer cells on the top layer reduces in favor of deeper connections. 

%mostrare il risultato dell'adattativo e fare vedere quali sono le regioni
%confrontandole con la permeabilita.
%discussione qualitativa di quello che otteniamo

%sistemare la condizione dei pozzi da metterla come quella del caso 2

%calcolare gli errori ottenuti

%dominio piu grande

%%% Local Variables:
%%% mode: latex
%%% TeX-master: "master"
%%% End:

%% file: conclusion.tex
\section{Conclusion}

In conclusion, our study presents a strategy for an automatic selection of linear- and nonlinear-flow regions in highly heterogeneous porous media, particularly answering the challenges raised by deviations from Darcy's law in highly permeable areas. Through a combination of theoretical derivations, empirical observations, and numerical simulations, we have demonstrated the significance of nonlinear corrections to Darcy's law by accurately capturing complex flow behaviors, especially in regions with high permeability.

Moreover, our approach, by using an adaptive regularized model based on a physically motivated flux threshold, offers a promising strategy for identifying and distinguishing between linear and nonlinear regimes within porous media. A local tolerance on the error, from which the threshold value is derived, has been introduced to switch between the linear and nonlinear laws effectively based on local physical properties. The experimental validation we have conducted across several numerical test cases shows the effectiveness of our methodology, although, as mentioned in the introduction, more investigation is needed to combine our strategy with a domain-decomposition technique and thus yield an complete, accurate and computationally efficient approach for flows in highly heterogeneous porous media.